\documentclass[12pt]{article}
\newtheorem{ass}{Assumption}[section]

\newtheorem{lm}{Lemma}[section]

\newtheorem{prop}{Proposition}[section]

\newtheorem{thm}{Theorem}[section]
\newcommand{\qed}{\hspace*{1mm} \rule{2mm}{3mm}\vspace*{0.3cm}}
\newcommand{\be}{\begin{equation}}
\newcommand{\ee}{\end{equation}}
\newcommand{\bea}{\begin{eqnarray}}
\newcommand{\eea}{\end{eqnarray}}
\newcommand{\ba}{\begin{array}}
\newcommand{\ea}{\end{array}}
\newcommand{\beas}{\begin{eqnarray*}}
\newcommand{\eeas}{\end{eqnarray*}}

\newcommand{\ZZ}{{\bf Z}}
\newcommand{\QQ}{{\bf Q}}
\newcommand{\RR}{{\bf R}}

\newcommand{\nn}{\nonumber}

\def\xitil{\widetilde{\xi}}

\def\xbar{\bar{x}}
\def\ybar{\bar{y}}
\def\tbar{\bar{t}}
\def\ubar{\bar{u}}
\def\zetabar{\bar{\zeta}}

\def\dist{\mbox{dist}}

\def\Var{\mbox{Var}}
\def\e{\varepsilon}
\def\mmGa{{\bf \Gamma}}
\def\mmL{{\bf L}}
\def\Lploc{L^p_{\rm loc}(\RR)}
\def\Hloc{H^{-1}_{\rm loc}(\RR)}

\def\aihpnl{Ann.\ Inst.\ H.\ Poincar\'e Anal.\ Non Lin\'eaire}

\def\ap{Ann.\ Probab.}

\def\ptrf{Probab.\ Theory Related Fields}

\begin{document}
 \title{Diffusive fluctuations for one-dimensional totally 
asymmetric interacting random dynamics}
 \author{Timo Sepp\"al\"ainen\\
 Department of Mathematics\\
 University of Wisconsin\\
 Madison, WI 53706-1388}

 \thispagestyle{empty}
 \maketitle
 \abstract{We study central limit theorems 
for a
 totally asymmetric, one-dimensional interacting random system. 
The models we work with are the Aldous-Diaconis-Hammersley 
process and the related stick model. The A-D-H process
represents a particle configuration on the line, 
or a 1-dimensional interface on the plane
which moves in one fixed direction through random local 
jumps. The stick model is the process of local slopes
of the A-D-H process, and has a conserved quantity. 
The results describe the fluctuations of these systems 
around the deterministic evolution to which the random 
system converges under hydrodynamic scaling. We 
look at diffusive fluctuations, by which we
mean  fluctuations on the scale of the classical 
central limit theorem. In the scaling limit these 
fluctuations obey deterministic equations with
random initial conditions given by the initial
fluctuations. Of particular interest is the effect
of macroscopic shocks, which play a  dominant role because
 dynamical noise is suppressed on the scale we are
working. 
 }

 \thanks{Research  partially supported by NSF grants DMS-9801085
and DMS-0126775.
 \\[.15cm]
 {\sl Abbreviated title}:  Fluctuations for totally asymmetric 
  systems \\[.15cm]
 {\sl AMS (1991) subject classifications}: Primary 60K35; secondary 60F05, 
82C22.}
 \eject

\section{Introduction}
We study fluctuations in the 
scale of the classical central limit theorem for totally asymmetric
interacting random systems in one space dimension. 
The model system for which we prove
 theorems is the Aldous-Diaconis-Hammersley process. 
To summarize this model in one sentence, it consists
of point particles on the real line that  jump to the left,
at rate equal to the distance to the left neighbor, with 
new locations chosen uniformly at random between the jumper
 and its left neighbor. 
The idea for this process appeared in Hammersley's
classical paper \cite{hammersley}, and Aldous and 
Diaconis \cite{alddiac} first defined it
 as an infinite
system of interacting particles. 

We consider the general nonequilibrium hydrodynamic limit
 situation 
where the limiting interface (or tagged particle,
depending on one's point of view)  
is governed by a Hamilton-Jacobi equation
$u_t+f(u_x)=0$. The initial distributions can be 
fairly arbitrary, subject  to a limit assumption
on the fluctuations around the initial macrosopic
profile  and some moment bounds. 
In particular, we do not restrict to product initial 
distributions or particular types of initial
macroscopic profiles.  

The overall picture
is this: the  limiting 
fluctuation field $\zeta(x,t)$ is governed by the linearization
of the hydrodynamic equation: 
$\zeta_t+f'(u_x)\zeta_x=0$ where $u(x,t)$ is the deterministic
limit around which the random interface fluctuates. 
This is a deterministic equation, and 
all the randomness is confined to the initial condition. 
The dynamics transports the initial fluctuations 
along the characteristics and shocks of the hydrodynamic 
equation. This  picture of characteristics rigidly 
transporting fluctuations has been understood to some degree
for quite a while, and has been proved in some special 
cases. What our paper furnishes are proofs 
 in a general setting (but for the particular model). 
In addition we clarify some 
interesting details of this picture that are produced
by the shocks of the hydrodynamic equation, such as 
the definition of the limiting fluctuation variable  $\zeta(x,t)$
at a shock location $(x,t)$.

The mathematical reason for the 
suppression of dynamical noise lies in two facts: (i) The
most general evolution of Hammersley's process can be realized
as an envelope of  an infinite family of simpler processes with 
deterministic initial conditions. This is the microscopic
variational representation of the process. (ii) The 
results of Baik-Deift-Johansson \cite{BDJ} imply that these simpler
processes have fluctuations of order $n^{1/3}$ which
are then swamped by the initial diffusive fluctuations 
of order  $n^{1/2}$. 

It has been more common to use the exclusion process and its
variants  for mathematical theory
of large scale behavior. For our purposes 
 Hammersley's process has one advantage
over the exclusion process. The totally asymmetric versions
of both processes can be 
conveniently coupled with simple growth models, and
both processes possess
particle-level variational formulations in this
coupling. For Hammersley's process this growth model
is the increasing sequences model on a planar Poisson point
process \cite{alddiac}\cite{Seejp}. 
For exclusion it is  the last-passage percolation
model with weakly increasing paths on the 
two-dimensional square lattice \cite{Semprf}. 
 The advantage of Hammersley's process comes from the fact
that presently better probability estimates are available
for the planar increasing sequences model than for the lattice
last-passage model. In particular, Lemma \ref{gammalm}(b)
in our proof has not yet been proved for the exclusion
model. This estimate for the increasing
sequences model was proved by Baik, Deift and Johansson
\cite{BDJ} with Riemann-Hilbert techniques. Obtaining
this estimate for exclusion is not simply a matter of
repeating the argument, but  has turned out to be a somewhat 
tricky problem (personal communication from J.\ Baik). 
 But once this estimate for exclusion becomes available, 
we believe that the results of this paper can be repeated 
for totally asymmetric simple exclusion.  

The reader can find comprehensive   overviews
of fluctuation results for interacting systems in \cite{gialebpre}
and in
 Chapter 11 of \cite{KipLan}. So we make only a few remarks here. 
Past work on the  fluctuations of asymmetric systems
has concentrated  on the exclusion process. 
The proofs  use couplings and monotonicity
arguments and necessitate special 
initial distributions such as i.i.d.\ distributions or
product measures with piecewise constant densities. 

The deepest and most important  work on the
fluctuations of the asymmetric 
exclusion process is undoubtedly   by
Ferrari and Fontes. In a series of papers
(\cite{ferrarifontes1},
 \cite{ferrarifontes2}, \cite{ferrarifontes3})
they study the fluctuations of the current and the 
tagged particle in equilibrium, and the fluctuations
of a second class particle with shock initial
conditions, given by 
a product measure  with different densities 
 to the left and right of the origin. In this situation the authors prove
the basic feature of asymmetric fluctuations, namely the
 rigid transport along the characteristics. Our paper complements their work
 on some questions,
by going into more general nonequilibrium profiles and 
initial distributions, and by giving more complete results
on the convergence of the entire interface and the 
distribution-valued density fluctuation field. 
Our results cover tagged particles for Hammersley's
process, and thereby also the current for the stick 
process. The stick process  is the process of increments for
Hammersley's process, and hence the displacements
of Hammersley's particles are the currents of the stick process. 

Let us  contrast our methods and results with those of
 symmetric, reversible processes. In the one-dimensional
setting symmetric means that particles 
are equally likely to jump both left and right. 
 The fluctuation theory
of reversible interacting processes relies on methods of martingales and
the Holley-Stroock theory of generalized
Ornstein-Uhlenbeck processes. 
The limiting fluctuation fields of reversible processes
 are governed by
equations driven by white noise. 
Both our methods and the qualitative results are different. We use
no martingale theory. Instead our methods rely on sharp 
control of the paths of individual particles, and on the theory of
shocks and characteristics of one-dimensional conservation
laws and Hamilton-Jacobi equations. And we already highlighted 
the main difference, that for asymmetric systems no dynamical
noise is visible on the diffusive scale. 

Between  symmetric and 
asymmetric systems are 
the weakly asymmetric systems where the 
asymmetry vanishes in the hydrodynamic limit. The central limit behavior
of weakly asymmetric systems  is qualitatively
the same as that of symmetric
systems, governed by a linear stochastic partial differential
equation whose drift term is the linearization of the 
hydrodynamic equation. But the  weakly asymmetric systems 
 have an additional interesting feature
proved by Bertini and Giacomin \cite{bertinigiacomin}:
A small perturbation of a flat profile obeys, on larger
space and time scales, a nonlinear stochastic equation
of KPZ type. 
This raises the question whether such a result could be 
obtained for  asymmetric systems at some suitable 
scaling. 

As mentioned above,   the 
 shocks of the hydrodynamic equation  
turn out to have interesting effects on 
fluctuations. For example, 
a basic result one would expect is that
the motion of a tagged particle converges to something
related to Brownian
motion. 
But we find that in the presence of shocks the
fluctuation processes of a
tagged particle are not tight in the Skorokhod space
$D([0,\infty),\RR)$. We can still prove the tagged
particle's convergence to a function of Brownian motion
uniformly on compact time intervals away from shocks, and 
even pointwise at the shocks. The limiting path
has discontinuities at the shock times, and   is not 
right-continuous (in time), but instead lower semicontinuous at
the discontinuities.

We prove a distributional limit theorem for the entire interface
in a weaker topology, as an element of $\Lploc$. This 
limiting process is a weak solution of the linearization of 
the hydrodynamic equation, as mentioned above. However, due to the shocks 
 a particular version
 of the limiting process has to be 
chosen (from among  the a.e.\ equal versions)
 to get a weak solution of the linearized 
equation. And this weak solution turns out to disagree with the pointwise
distributional limit at the shocks. 

{\it Organization of the paper.} Section \ref{sectresults} 
describes the particle process  and the results.
The last part of that section gives a rigorous construction
of the process in terms of increasing sequences among Poisson 
points on the Euclidean plane. This is the variational
coupling formulation basic for our approach.
 Section \ref{charsection}
develops properties of the characteristics and shocks of 
a one-dimensional conservation law. The approach here is 
based on the Hopf-Lax and Lax-Oleinik formulas, with the 
results of \cite{Re} as a starting point.  

Sections \ref{gammasect} and \ref{varcoupsect} contain 
probability estimates needed for taking advantage of the
variational coupling formulation. These are based on the
known estimates for increasing sequences (\cite{BDJ}, \cite{Kim},
\cite{Seldp}). The remaining sections 
go through the proofs of the theorems. Some 
technical  measurability proofs are collected in an 
appendix at the end. 

\section{Results}
\label{sectresults}

We study the large scale behavior of 
the  {\it Aldous-Diaconis-Hammersley process}, or
{\it Hammersley's process} for short. The state of  this process
is $z(t)=(z_i(t):i\in\ZZ)$ that represents a 
countable collection of labeled
 point particles 
on $\RR$. The variable $z_i(t)\in\RR$ 
is the location of particle $i$ at time
$t$.  The particles are ordered, so
 that $z_{i-1}(t)\le z_i(t)$ for all 
$i\in\ZZ$ and $t\ge 0$. 
All particles make jumps to the left, 
 according to the following  rule. Suppose
the state at time $t$ is  $z(t)=(z_i(t):i\in\ZZ)$.
To determine the next jump of
particle $i$, let 
$\sigma$ be a random exponential waiting  time 
with expectation $\left(z_i(t)-z_{i-1}(t)\right)^{-1}$. 
At time $t+\sigma$ particle $i$ jumps to its new location
$z_i(t+\sigma)$,  chosen uniformly at random from the
interval $(z_{i-1}(t), z_i(t))$. This type of event happens
independently and simultaneously for all $i$. Of course
this description needs justification because infinitely many
jumps happen in every positive time interval. 
 In Section \ref{construction}
below we give a rigorous construction 
of this infinite-particle dynamics in terms of
increasing sequences on the plane. 

Instead of thinking
about a particle configuration, we can regard Hammersley's
process as a model for a 1-dimensional interface on the 
plane. The interface is represented by the height function
$z(t)$ defined on the integers, so that $z_i(t)$ is the 
height of the interface above site $i$. Through the jumps
of the $z_i$'s the interface moves downward.  

The {\it stick process}  $\eta(t)=(\eta_i(t):i\in\ZZ)$ is the 
process of  increment variables
defined 
by
$$\eta_i(t)=z_i(t)-z_{i-1}(t).$$ The dynamics of $\eta(\cdot)$
can be represented by the following generator $\cal L$ which acts 
on bounded cylinder functions $\psi$ on the product space
$[0,\infty)^\ZZ$: 
$$
{\cal L}\psi(\eta)=\sum_{i\in\ZZ}\int_0^{\eta_i}
[\psi(\eta^{u,i,i+1})-\psi(\eta)] du
$$
where $\eta^{u,i,i+1}$ represents the configuration after 
a piece of size $u$ has been moved from site $i$ to $i+1$:
$\eta^{u,i,i+1}_i=\eta_i-u$, $\eta^{u,i,i+1}_{i+1}=\eta_{i+1}+u$, 
and
$\eta^{u,i,i+1}_j=\eta_j$ for $j\ne i,i+1$. This process can 
be rigorously defined on a certain subspace of the full product space
$[0,\infty)^\ZZ$, see \cite{Seejp} for details. 

Let  $u_0$ be a nondecreasing
locally Lipschitz continuous function on $\RR$. 
It represents the  initial macroscopic  interface. 
The evolving macroscopic
 interface $u(x,t)$, $(x,t)\in\RR\times[0,\infty)$,
 is the unique viscosity
solution of the Hamilton-Jacobi equation 
\be
u_t+f(u_x)=0,\qquad u(x,0)=u_0(x), 
\label{hjeqn}
\ee
with velocity function $f(\rho)=\rho^2$. Equivalently,
 $u$ is defined for $t>0$ by the Hopf-Lax formula 
\be
u(x,t)=\inf_{y:y\le x}\left\{ u_0(y)+tg\left(\frac{x-y}t\right)
\right\}
\label{hopflax}
\ee
where $g(x)=x^2/4$ is the convex dual of $f$. 
For a fixed $t$ the partial $x$-derivative $\rho(x,t)=u_x(x,t)$ exists 
for all but countably many $x$. This function is the unique entropy
solution of the Burgers equation
\be
\rho_t+f(\rho)_x=0,\qquad \rho(x,0)=\rho_0(x),
\label{burgerseqn}
\ee
where $\rho_0=u_0'$ (a.e.\ defined derivative). 
  We cover some properties of these
equations later in
Section \ref{charsection}. See chapters 3, 10, 11 in \cite{Ev}
for basic theory. 
 
Assume we have  a sequence $z^n(\cdot)$ of Hammersley's processes,
with random initial configurations $\{z^n_i(0):i\in\ZZ\}$, and 
$n=1,2,3,\ldots$ is the index of the sequence.
The objective of our paper is to study the fluctuations 
of the random interface $z^n_{[nx]}(nt)$ around the deterministic 
interface $nu(x,t)$ in the diffusive, or central limit theorem,
scale  $n^{1/2}$. The fluctuations are described by the 
stochastic process
$\zeta_n(x,t)$ defined for $(x,t)\in\RR\times[0,\infty)$ by 
\be
 \zeta_n(x,t)=n^{-1/2} \{ z^n_{[nx]}(nt)-nu(x,t)\}. 
\label{zetanxtdef}
\ee

Think of the initial process $\{\zeta_n(y,0):y\in\RR\}$ as  a 
random function with values in the Skorokhod space $D(\RR)$
of right-continuous functions on $\RR$ with left limits (RCLL functions). 
This space is metrized as follows. Let $\Lambda$ be the 
collection of strictly increasing,
 bijective Lipschitz functions $\lambda:\RR\to\RR$
such that
\be
\|\lambda\|=|\lambda(0)|+\sup_{x\ne y}
\left|\log\frac{\lambda(x)-\lambda(y)}{x-y}\right|<\infty.
\label{deflambdanorm}
\ee
For $\alpha,\beta \in D(\RR)$ and $u>0$ let
$$
d(\alpha,\beta ,\lambda,u)=\sup_{x\in\RR}\left|
\alpha\left((x\wedge u)\vee(-u)\right)-
\beta\left((\lambda(x)\wedge u)\vee(-u)\right)
\right|\wedge 1
$$
and then 
\be
d_S(\alpha,\beta )=\inf_{\lambda\in\Lambda}\left[ \|\lambda\|  
+  \int_0^\infty e^{-u} d(\alpha,\beta ,\lambda,u)du  
 \right].
\label{dSdef}
\ee
The metric $d_S$ is complete and separable. 
Convergence $d_S(\alpha_j,\alpha)\to 0$ is equivalent to the 
existence of a sequence $\lambda_j\in\Lambda$ such that 
$\lambda_j$ converges to the identity function uniformly on
compacts, and $|\alpha_j-\alpha\circ\lambda_j|\to 0$ uniformly on compacts. 
Let $C(\RR)$ denote the subspace of continuous functions. 

Our basic hypothesis is   weak convergence
at time 0 to a continuous limit function:
\be
\begin{array}{rl}
&\mbox{There exists a $C(\RR)$-valued  random function
 $ \zeta_0$ such that 
 }\nn\\ 
&\mbox{$\zeta_n(\cdot,0)\to \zeta_0(\cdot)$ 
in distribution as $n\to\infty$, on the space $D(\RR)$.}
\end{array}
\label{ass1}
\ee

Assumption (\ref{ass1}) is in fact equivalent to a stronger
assumption, which is important for us so we clarify it 
right away.  
Let $D_u(\RR)$ be the space of 
RCLL functions endowed with the $d_u$-metric of uniform 
convergence on compact sets:
\be
d_u(\alpha,\beta )=\sum_{j=1}^\infty 2^{-j} \left\{\sup_{-j\le r\le j}
|\alpha(r)-\beta(r)|\wedge 1\right\}
\quad
\mbox{for $\alpha,\beta \in D_u(\RR)$.}
\label{dudef}
\ee
The  metric $d_u$ is much stronger than the 
Skorokhod metric $d_S$, and in fact $D_u(\RR)$ 
is not even separable. But on $C(\RR)$ the two metrics induce the
same topologies.   Because the jumps of $\zeta_n(\cdot,0)$ occur 
at deterministic  
locations and because the limit process $\zeta_0(\cdot)$ is
continuous, it follows  that $\zeta_n(\cdot,0)$ is measurable as a 
$D_u(\RR)$-valued random function, and assumption (\ref{ass1})
is equivalent to this stronger assumption:
\be
\begin{array}{rl}
&\mbox{There exists a $C(\RR)$-valued  random function
 $ \zeta_0$ such that }\nn\\ 
&\mbox{$\zeta_n(\cdot,0)\to \zeta_0(\cdot)$ 
in distribution as $n\to\infty$, on the space $D_u(\RR)$.}
\end{array}
\label{ass1d}
\ee
We shall not go through the details of this 
point, and refer the reader to
section 18 in \cite{billingsley}.

Since the state space is large, 
 we need a uniformity assumption. But only on one side since the 
dynamics is totally asymmetric. 
\be\begin{array}{rl}
&\mbox{There exists a fixed $b\in\RR$ such that for every $\e>0$ 
one can find $q$ and $n_0$ }\\
&\mbox{such that }\quad \displaystyle
\sup_{n\ge n_0} P\left\{ \sup_{k: k\le nq}{n}{k^{-2}}
\left( z^n_{[nb]}(0)-z^n_k(0)\right) \ge \e\right\} \le \e.
\end{array}
\label{unifass}
\ee
Note that if (\ref{unifass}) holds for some $b$, it holds for all $b$. 
It forces $u_0$ to satisfy 
\be
\lim_{y\searrow -\infty}|y|^{-2}u_0(y)= 0. 
\label{assu0}
\ee
A consequence of assumption (\ref{ass1d}) is convergence in 
probability to the macroscopic interface $u_0$: 
\be
\lim_{n\to\infty}P\left(\sup_{y\in[a,b]}
 |n^{-1}z^n_{[ny]}(0)-u_0(y)|\ge\e\right)=0.
\label{hydrolim0}
\ee
This  and (\ref{unifass}) are sufficient for
a hydrodynamic limit: 
$n^{-1}z^n_{[nx]}(nt)\to u(x,t)$ in probability as $n\to\infty$,
uniformly over $(x,t)$ in compact sets. 
See \cite{Seejp}. 

Property (\ref{assu0}) guarantees
that there exists a nonempty compact set 
$I(x,t)\subseteq(-\infty,x]$
 on which the  infimum in (\ref{hopflax}) is achieved:
$$I(x,t)=\left\{y\le x: u(x,t) =u_0(y)+tg\left(\frac{x-y}{t}\right)
\right\}.
$$
 For $t=0$ it is convenient to have 
the convention $I(x,0)=\{x\}$. 
The minimal and maximal Hopf-Lax minimizers are 
\be
y^-(x,t)=\inf I(x,t)\qquad\mbox{and}\qquad y^+(x,t)=\sup I(x,t).
\label{earlyypmdef}
\ee
Define 
\be
\rho^\pm(x,t)=g'\left( \frac{x-y^\pm(x,t)}t\right)\qquad
\mbox{for $(x,t)\in\RR\times(0,\infty)$.}
\label{earlylaxoleinik}
\ee
 It  turns out that,  for a fixed $t$,   $y^-(x,t)=y^+(x,t)$
for all except at most countably many $x$. At all such points
the function
$\rho(x,t)=\rho^\pm(x,t)$ is defined   and  continuous, 
and is  the $x$-derivative $\rho(x,t)=u_x(x,t)$ of the 
viscosity solution of (\ref{hjeqn}). Definition
(\ref{earlylaxoleinik}) is called the Lax-Oleinik formula. 
 We say
that $(x,t)\in\RR\times(0,\infty)$ is the location of a {\it shock} if
$y^-(x,t)<y^+(x,t)$. We will  not call $(x,0)$ a shock even if the
initial function $u_0$ is nondifferentiable at $x$. 

Our first result shows that later fluctuations are close
to a deterministic transformation of the initial fluctuations.

\begin{thm} Suppose $u_0$ is a locally
Lipschitz continuous function. 
 Assume  {\rm (\ref{ass1d})} and {\rm (\ref{unifass})}.

{\rm (i)}
Let $A\subseteq\RR\times[0,\infty)$ be a compact set such that
either (a) $A$ is finite, or (b) there are no shocks in $A$, in other
words $y^-(x,t)=y^+(x,t)$ for all $(x,t)\in A$. Then 
\be
\lim_{n\to\infty} \sup_{(x,t)\in A}
\left| \zeta_n(x,t)-\inf_{y\in I(x,t)} 
\zeta_n(y,0) \right| =0
\qquad\mbox{in probability.}
\label{problim}
\ee

{\rm (ii)} For all $-\infty<a<b<\infty$, 
 $0<\tau<\infty$, and $1\le p<\infty$, 
\be 
\lim_{n\to\infty} \sup_{0\le t\le \tau}\int_a^b 
\left| \zeta_n(x,t)-\inf_{y\in I(x,t)} 
\zeta_n(y,0) \right|^p dx =0
\qquad\mbox{in probability.}
\label{problim2}
\ee
\label{probthm}
\end{thm}

 From this theorem  we  deduce distributional limits for the 
interface and the stick profile.

\subsection{Weak limits and the linearized equation}

In assumption (\ref{ass1d}) we assumed the existence
of a $C(\RR)$-valued random function $\zeta_0$. On the 
probability space 
of $\zeta_0$ define 
 random variables  $\zeta(x,t)$,  $(x,t)\in\RR\times[0,\infty)$, by
\be
\zeta(x,t)=\inf_{y\in I(x,t)} \zeta_0(y)\,.
\label{zetaxtdef}
\ee

To formulate a process-level weak convergence result, we
consider, for a fixed $t$, the random function 
$x\mapsto \zeta_n(x,t)$ as an element of the space
$\Lploc$ of functions that are locally in $L^p$. By definition,
a measurable function $f$ on $\RR$ lies in $\Lploc$ if
for all $0<k<\infty$, 
$$
\|f\|_{L^p[-k,k]}\equiv \left(\int_{[-k,k]} |f(x)|^p dx\right)^{1/p}<\infty.
$$
$\Lploc$ is a complete separable metric space under the metric
\be
d_p(f,g)=\sum_{k=1}^\infty 2^{-k} \left( \|f-g\|_{L^p[-k,k]}\wedge 1\right),
\label{dpdef}
\ee
and we endow $\Lploc$ with its Borel $\sigma$-algebra. 
We show 
 that for a fixed $t$, $\zeta_n(\cdot,t)$ is 
measurable as an $\Lploc$-valued random element.   
And  that the path $\zeta_n:t\mapsto \zeta_n(\cdot,t)$ is a measurable map 
from the underlying probability space into the 
Skorokhod space $D\left([0,\infty),\Lploc\right)$ of 
right-continuous $\Lploc$-valued paths  with left limits at all time
points $t$. Similarly the random variables $\zeta(x,t)$ defined in 
(\ref{zetaxtdef}) specify an $\Lploc$-valued path 
$\zeta: t\mapsto \zeta(\cdot,t)$. We show that $\zeta$ is 
  a random
element  of the space $C\left([0,\infty),\Lploc\right)$ of 
continuous paths.

\begin{thm} Suppose $u_0$ is a locally
Lipschitz continuous function. 
 Assume {\rm (\ref{ass1d})} and {\rm (\ref{unifass})}.

{\rm (i)} For any finitely many points $(x_i,t_i)\in\RR\times[0,\infty)$,
$1\le i\le k$, we have the limit in distribution
\be
(\zeta_n(x_1,t_1),\ldots,\zeta_n(x_k,t_k))
\stackrel{d}\longrightarrow (\zeta(x_1,t_1),\ldots,\zeta(x_k,t_k))
\qquad\mbox{ as $n\to\infty$}
\label{weaklim1}
\ee
in the space $\RR^k$. 

{\rm (ii)} The process $\zeta_n$ converges in distribution  to the process
$\zeta$ on the path space $D\left([0,\infty), \Lploc\right)$.
\label{weaklimitthm}
\end{thm}

 As one would expect, 
$\zeta(x,t)$ is a solution of the 
linearization of the Hamilton-Jacobi equation (\ref{hjeqn}). For this 
 we must choose the correct version of $\zeta$ in the a.e.\ sense. Let
$\zetabar(x,t)=\frac12\left\{\zeta_0(y^-(x,t))+\zeta_0(y^+(x,t))\right\}$. 
For a fixed $t$, 
$\zetabar(x,t)=\zeta(x,t)$ at all $x$ except 
shock locations.  
 $\zetabar$ is a weak solution of the equation
\be
\zetabar_t(x,t)+f'(\rho(x,t))\zetabar_x(x,t)=0\,,\qquad 
\zetabar(\cdot,0)=\zeta_0(\cdot).
\label{transeqn}
\ee
This is a linear transport equation with a discontinuous coefficient. 
The appropriate definition of a weak solution is that, for all
$\phi\in C_c^\infty(\RR\times[0,\infty))$, $\zetabar$ satisfies
this integral criterion:
\bea
&&\int_0^\infty \int_{\RR} \phi_t(x,t) \zetabar(x,t) dx\,dt
+\int_0^\infty dt \int_{\RR}
\zetabar(x,t) d[\phi(\cdot,t)f'(\rho^+(\cdot,t))](x)\nn\\
&&\qquad\qquad 
+\int_{\RR} \zeta_0(x)\phi(x,0)dx=0.
\label{weaktranseqn}
\eea
For each $t$, the $x$-integral in the second term is with respect to the signed measure
$\mu=\mu(t)$
defined by 
$$
\mu(a,b]=\phi(b,t)f'(\rho^+(b,t))-\phi(a,t)f'(\rho^+(a,t)).
$$
For this to make sense we took the right-continuous version
$\rho^+(\cdot, t)$ of 
$\rho(\cdot, t)$. The definition also requires that $\rho(\cdot, t)$
be locally of bounded variation, which is true by 
the Lax-Oleinik formula (\ref{earlylaxoleinik}).  Equation 
(\ref{weaktranseqn}) shows why the choice of $\zetabar$ matters.
Suppose $(r(t),t)$ is a shock location for $t_0\le t\le t_1$. 
 Then $f'(\rho(\cdot,t))$ jumps
at $r(t)$ and the measure $\mu$ gives nonzero mass to the singleton
$\{r(t)\}$ for each $t$. Clearly the value of the second term 
in (\ref{weaktranseqn}) depends 
on which value $\zetabar(r(t),t)$ takes. 
It is a curious discord that the correct weak solution of 
(\ref{transeqn}) differs from the pointwise limit in (\ref{weaklim1})
at the shocks. 
There is no
dynamically generated noise in  equation (\ref{transeqn}), 
as all the randomness is in the initial
data $\zeta_0$. The equation expresses the point 
that on the diffusive scale the initial noise is transported
along the characteristics, and the noise created by the dynamics 
is not visible because it is of lower order. 

That $\zetabar$ is a weak solution of (\ref{transeqn}) follows
from this more general result.  Given a  convex, differentiable
 flux function 
$f$, let $\Theta(\lambda,\rho)\in[0,1]$ for $\lambda\ne \rho$ be defined by
\be
\frac{f(\lambda)-f(\rho)}{\lambda-\rho}
=\Theta(\lambda,\rho)f'(\rho)+(1-\Theta(\lambda,\rho))f'(\lambda).
\label{Thetadef}
\ee
Let $\rho^\pm(x,t)$ be the functions defined by the Lax-Oleinik
formula (\ref{earlylaxoleinik}). Given a continuous function
$v_0$, set for $(x,t)\in\RR\times(0,\infty)$ first
$$
\theta(x,t)=\Theta(\rho^-(x,t),\rho^+(x,t))
$$
and then 
\be
v(x,t)=\theta(x,t)v_0(y^+(x,t))
+\left(1-\theta(x,t)\right)v_0(y^-(x,t)).
\label{vdef}
\ee

\begin{thm} Suppose $f$ is a convex flux function with convex
conjugate $g$, 
the minimizers $y^\pm(x,t)$ are defined by {\rm (\ref{earlyypmdef})},
and $\rho^\pm(x,t)$ are defined by the Lax-Oleinik 
formula {\rm (\ref{earlylaxoleinik})}. Let $v_0$ be an 
arbitrary continuous function on $\RR$, and define $v$ by
{\rm (\ref{vdef})}. Then $v$ is a weak solution of the linear
transport equation 
\be
v_t+f'(\rho(x,t))v_x=0\,,\quad v|_{t=0}=v_0,
\label{vtranseqn}
\ee
in the sense of the integral criterion {\rm (\ref{weaktranseqn})}.
\label{lineqnthm}
\end{thm}

We would expect $v$ to be the unique
weak solution of (\ref{vtranseqn}) under some natural uniqueness
criterion. Presently a uniqueness theory
 exists for continuous solutions of equations of this
type. See Petrova and Popov
\cite{petrovapopov} and their references. 

 For the special case $f(\rho)=\rho^2$ we get 
$\Theta\equiv\frac12$, which explains why we defined 
$\zetabar$ as the $\frac12,\frac12$ convex combination
of $\zeta_0(y^\pm(x,t))$.
Next some remarks on the hypotheses
and results.

\hbox{} 

\subsubsection{Remark}
\label{contremark}
Above we chose to work with the
$x$-right-continuous function $\zeta_n(x,t)$ defined
by (\ref{zetanxtdef}). 
The reader may prefer to linearly interpolate 
between the point locations $z^n_k$ to define 
an $x$-continuous random interface  
\be
z_n(x,t)=\left(nx-[nx]\right) z^n_{[nx]+1}(nt)
+\left([nx]+1-nx\right) z^n_{[nx]}(nt),
\label{contzndef}
\ee
and then consider the $x$-continuous fluctuation process
$$\zeta_n^{(c)}(x,t)=n^{-1/2}\{ z_n(x,t)-nu(x,t)\}.$$ 
The results would be the same. In particular, assumption
(\ref{ass1d}) is equivalent to $\zeta_n^{(c)}(\cdot,0)\to \zeta_0(\cdot)$
weakly in $C(\RR)$.  
Our estimates imply the following
proposition, which shows that on the scale $n^{1/2}$
large microscopic variations in the index are not
visible. 

\begin{prop} Suppose $0\le \ell=\ell(n)\le Cn^{1/3-\delta}$
for some $C<\infty$ and $\delta>0$. Fix $-\infty<a<b<\infty$ and
$\tau<\infty$. 
 Under assumptions {\rm (\ref{ass1d})} and {\rm (\ref{unifass})},
\be
\lim_{n\to\infty}\sup_{an\le k\le nb\,,\, 0\le t\le \tau}
n^{-1/2}\{ z^n_{k+\ell}(nt)-z^n_k(nt)\}=0
\qquad\mbox{in probability.}
\label{kelllim}
\ee
Under the stronger assumptions {\rm (\ref{assrho})} and 
{\rm  (\ref{assloqeq})}
of the next section, the limit above holds a.s.
\label{kellprop}
\end{prop}

We sketch the proof of this proposition in the Appendix.
 

\subsubsection{Remark} In both theorems part (i) is a sharper statement 
for a restricted set of space-time points, and part (ii) is a weaker statement 
without restriction on space-time points. Let us emphasize that the limits 
in (\ref{problim}) and (\ref{weaklim1}) are valid for
any {\it finite}
 collection of points, including shock locations.  For the global results,
 (\ref{problim2}) and part (ii) of Theorem \ref{weaklimitthm},
 we  integrate over space so that the values of the processes at 
 shocks become  immaterial because the shocks 
are a Lebesgue null set. The same effect could be achieved by integrating
over time. 


\subsubsection{Remark}
\label{unifremark}
 The uniform convergence in (\ref{problim}) cannot be
extended to sets that contain shocks. To see why, suppose 
 $(x,t)$ is a shock, and suppose there are points 
$x_k\nearrow x$ and $x_\ell'\searrow x$ such that 
 $(x_k,t)$ and $(x_\ell',t)$  lie in $A$ but are not shocks. 
Let  $y_k=y^\pm(x_k,t)$ and 
$y_\ell'=y^\pm(x_\ell',t)$ be the 
 (unique) Hopf-Lax minimizers for these points. They 
satisfy  $y_k\nearrow y^-(x,t)$ and 
$y_\ell'\searrow y^+(x,t)$. Consider the $x$-continuous 
version $\zeta_n^{(c)}(x,t)$ defined in Remark \ref{contremark}.
 If 
(\ref{problim}) were to  hold for this set $A$, then with high probability 
$|\zeta_n^{(c)}(x_k,t)-\zeta_n^{(c)}(y_k,0)|$ and 
$|\zeta_n^{(c)}(x_\ell',t)-\zeta_n^{(c)}(y_\ell',0)|$
are small uniformly over $k$ and $\ell$. As we let $k,\ell\to\infty$
and use the continuity of $\zeta_n^{(c)}(\cdot, t)$, we conclude that 
the random variable $\zeta_n^{(c)}(x,t)$ is forced to simultaneously
approximate $\zeta_n^{(c)}(y^-(x,t),0)$ and  $\zeta_n^{(c)}(y^+(x,t),0)$.
This is impossible (except in trivial cases) because
in the shock case the points $y^-(x,t)$ and $y^+(x,t)$ are 
macroscopically separated, and the values $\zeta_n^{(c)}(y^\pm(x,t),0)$
can differ with probability 1 with suitable choice of initial 
distributions.

\subsection{Starting in local equilibrium}

Now we take the point of view of an observer 
 on the initial interface, whose  location 
 is taken as  the origin. Furthermore, 
we assume that at time zero this observer sees the 
interface to his left and right in local equilibrium,
which is an assumption on the local slopes $\eta^n_i(0)=z^n_i(0)-z^n_{i-1}(0)$.
(If we want to think of the $z^n_i$'s as particles, we call
the $\eta^n_i$'s  interparticle 
distances.) Then we can strengthen the distributional limits to almost sure
limits, and give the limiting objects concrete descriptions 
in terms of Brownian motion.

For the precise hypotheses, let $\rho_0$ be 
a nonnegative, locally bounded  measurable function
on $\RR$. It will be  the macroscopic profile of the 
$\eta^n_i(0)$ variables. Assume that for some real number
$b$ (and hence for all $b$), 
\be
\lim_{r\to-\infty}  |r|^{-1}\cdot \sup_{r\le x\le b} \rho_0(x)=0.
\label{assrho}
\ee
Define a locally Lipschitz function $u_0$ by
$$
u_0(0)=0\,,\quad u_0(x)-u_0(y)=\int_y^x \rho_0(r)dr
\quad\mbox{ for all $y<x$.}
$$
Let $u(x,t)$ and $\rho(x,t)$ be again the relevant solutions 
of the macroscopic equations (\ref{hjeqn}) and (\ref{burgerseqn}). 
The assumption on the initial interfaces $z^n(0)$ is as follows.
\be
\begin{array}{ll}
&\mbox{For each $n$, $z^n_0(0)=0$ with probability 1, and 
        the variables $(\eta^n_i(0):i\in\ZZ)$}\\
&\mbox{are mutually independent, exponentially distributed  
with expectations}\\
&\mbox{$\displaystyle E[\eta^n_i(0)]=nu_0({i}/{n})-
nu_0({(i-1)}/{n})=n\int_{(i-1)/n}^{i/n}\rho_0(x)dx.$}
\end{array}
\label{assloqeq}
\ee

Note that now $-z^n_0(t)$ is the cumulative current from
site 0 for the stick process $\eta^n(\cdot)$, in other words
the total stick length that has moved across the bond $(0,1)$
during time interval $(0,t]$.  

Let $B(\cdot)$ denote a two-sided standard Brownian motion.
In other words, 
take two independent 1-dimensional standard Brownian motions
 $B_1(s)$ and  $B_2(s)$ defined for $0\le s<\infty$,
and set 
\be
B(s)=\left\{  \begin{array}{ll}
      B_1(s), &s\ge 0\\
        -B_2(-s), &s<0.
\end{array}
\right.
\ee
The limiting processes are defined in terms of this Brownian motion by  
$$
\zeta_0(y)=B\left( \int_0^y \rho_0^2(s)ds\right)
$$
and
\be
\zeta(x,t)=\inf_{y\in I(x,t)}B\left( \int_0^y \rho_0^2(s)ds\right)
=\inf_{y\in I(x,t)} \zeta_0(y)\,.
\label{zetadefB}
\ee
Note that in the above definitions
the integrals are 
signed, in other words for $y<0$ 
 $\int_0^y \rho_0^2(s) ds=-\int_y^0 \rho_0^2(s) ds \le 0$. 

\begin{thm} Assume {\rm (\ref{assrho})} and {\rm (\ref{assloqeq})}. 
Then we can 
construct the processes $\{z^n(\cdot)\}$ on a common probability
space  with a two-sided Brownian motion $B(\cdot)$ so that 
the following almost sure limits hold.

{\rm (i)} 
Let $A\subseteq\RR\times[0,\infty)$ be a compact set such that
either (a) $A$ is finite, or (b) there are no shocks in $A$, in other
words $y^-(x,t)=y^+(x,t)$ for all $(x,t)\in A$. Then
\be
\lim_{n\to\infty} \sup_{(x,t)\in A}
\left| \zeta_n(x,t) - 
\zeta(x,t)\right|=0\qquad\mbox{a.s.}
\label{limznxtB}
\ee

{\rm (ii)} For all $-\infty<a<b<\infty$,  $0<\tau<\infty$, 
and $1\le p<\infty$, 
\be
\lim_{n\to\infty} \sup_{0\le t\le\tau}
\int_a^b \left| \zeta_n(x,t) - 
\zeta(x,t)\right|^p\,dx=0\qquad\mbox{a.s.}
\label{limznxtB2}
\ee
\label{thmloqeq}
\end{thm}

\subsubsection{Remark}
\label{distremark}
 We assumed the initial increment
variables $\{\eta^n_i(0)\}$ exponentially distributed 
in assumption (\ref{assloqeq}) just to be concrete. It is a
natural choice 
because i.i.d.\ exponential distributions are invariant
 for the $\eta(\cdot)$ process so we can call 
(\ref{assloqeq}) ``local equilibrium.'' But the validity
of Theorem \ref{thmloqeq} does not depend on this special
choice at all. The reader can verify that the proof works
 as long as 
the initial distribution can be embedded in Brownian motion, 
and the moments are sufficiently bounded so that the 
probabilities in (\ref{unifass}) and (\ref{hydrolim0}) 
are summable in $n$. However,  definition (\ref{zetadefB})
of $\zeta(x,t)$ would change with different choices of
initial distributions. The $ \rho_0^2(s)$ inside the integral
$\int_0^y \rho_0^2(s)ds$ appears because the variance of an exponential random
variable is the square of the mean. 

\subsubsection{Moving along a characteristic from the origin} 
If $y^\pm(x,t)=0$, which means that $(x,t)$ is a
point on a genuine characteristic (not a shock) 
 emanating from $(0,0)$, then $\zeta(x,t)=0$ and
 (\ref{limznxtB}) gives 
$\zeta_n(x,t)\to 0$ a.s. This tells us that $n^{-1/2}$ is the wrong
normalization. We might expect the fluctuation to be 
of size  $n^{1/3}$ because the situation studied by 
Baik, Deift and Johansson \cite{BDJ} is of this 
type. Their initial condition corresponds to setting
$z_i(0)=0$ for $i\le 0$ and $z_i(0)=\infty$ for $i>0$. 
And their result can be expressed as the weak limit
of $n^{-1/3}\{z_{[nx]}(nt)-nx^2/(4t)\}$ for $x,t>0$. 
In this situation $u_0(x)=\infty\cdot{\bf 1}_{(0,\infty)}(x)$
and  $y^\pm(x,t)=0$ for all  $x,t>0$. 

On the other hand, suppose $y^-(x,t)\le 0\le y^+(x,t)$
with at least one inequality strict. Then
$(x,t)$ lies on a  characteristic from the origin
that is a shock. Now $\zeta(x,t)\ne 0$ with positive probability,
and with probability 1 if  $y^-(x,t)< 0< y^+(x,t)$.
 (\ref{limznxtB})  says that the current across a shock
has fluctuations of order $n^{1/2}$. 

\subsubsection{Rarefaction fan}
This means that $y^\pm(x,t)=\ybar$ for a nontrivial 
interval of $x$'s. The simplest way to produce this is to take
two densities $\lambda>\rho$ and  the initial profile 
$$
\rho_0(y)=\left\{
\begin{array}{rl}
\rho, &y<\ybar\\
\lambda, &y>\ybar.
\end{array}
\right.
$$
Then  $y^\pm(x,t)=\ybar$ for $(x,t)\in F$
where  $F$ denotes the ``fan'' (cut off at $T<\infty$ to make it compact)
$$
F=\{(x,t): 0\le t\le T\,,\, \ybar+2\rho t\le x\le
\ybar+2\lambda t \}.
$$
As a corollary of (\ref{limznxtB}) we get 
$$
\lim_{n\to\infty}
\sup_{(x,t),(x',t')\in F}\left| \zeta_n(x,t)- \zeta_n(x',t')\right|=0. 
$$
So $n^{-1/2}$ is not the right normalization for fluctuations
inside a rarefaction fan, and further work is called for. 

\subsubsection{Shock} A shock produces discontinuous
fluctuations that jump across segments of the Brownian path
that represents the initial fluctuations. Consider the simplest
shock case, with initial profile 
$$
\rho_0(y)=\left\{
\begin{array}{rl}
\lambda, &y<0\\
\rho, &y>0,
\end{array}
\right. 
$$
where still $\lambda>\rho$. The convex flux
$f(\rho)=\rho^2$ preserves a downward jump. 
(An upward jump is smoothed out into the rarefaction
fan.) At later times $t>0$ the 
shock is located at  $x=(\rho+\lambda)t$, and the profile is given by
$$
\rho(x,t)=\left\{
\begin{array}{rl}
\lambda, &x<(\rho+\lambda)t\\
\rho, &x>(\rho+\lambda)t.
\end{array}
\right.
$$ 
The Hopf-Lax minimizers are 
$y^\pm(x,t)=x-2\lambda t$ for $x<(\rho+\lambda)t$, 
$y^\pm(x,t)=x-2\rho t$ for $x>(\rho+\lambda)t$, 
and $I(x,t)=\{(\rho-\lambda)t, (\lambda-\rho)t\}$ 
for $x=(\rho+\lambda)t$. 
At macroscopic time $t$, the limiting fluctuation process
 is 
\be
\zeta(x,t)=\left\{
\begin{array}{ll}
B\left(\lambda^2(x-2\lambda t)\right), &x<(\rho+\lambda)t\\
\min\{ B\left(\lambda^2t(\rho-\lambda )\right), 
B\left(\rho^2t(\lambda-\rho )\right) \},
&x=(\rho+\lambda)t\\
B\left(\rho^2(x-2\rho t)\right), &x>(\rho+\lambda)t.
\end{array}
\right.
\label{limfluct}
\ee
There is a jump in $\zeta(\cdot,t)$  at 
the shock $x=(\rho+\lambda)t$, and the path may be 
left- or right-continuous, depending on which choice
makes it lower semicontinuous. The initial
fluctuation in the range 
$\{B(s): \lambda^2t(\rho-\lambda)<s<\rho^2t(\lambda-\rho )\}$
 disappeared from (\ref{limfluct}).   Ferrari and Fontes 
\cite{ferrarifontes1} show that in asymmetric exclusion
 this becomes the 
fluctuation of a second class particle.

\subsubsection{A tagged particle fails to be tight in the 
presence of shocks}
A basic question is to ask 
about the fluctuations of the motion of a tagged particle. 
In other words, fix $x$ and consider  the process
$\zeta_n(x,t)=n^{-1/2}\{z^n_{[nx]}(nt)-nu(x,t)\}$
as $t$ varies in $[0,T]$. If there are no shocks 
in $\{x\}\times[0,T]$, (\ref{limznxtB}) gives uniform 
convergence to a time-changed Brownian 
path. 

But if  $(x,\sigma)$ 
is a shock  for some 
$\sigma\in(0,T)$, it turns out that the sequence
of processes $\{\zeta_n(x,\cdot)\}$ is not
even tight in the Skorokhod space $D([0,T],\RR)$. To see
this, recall this condition for tightness: for every 
$\e>0$ there must exist a $\delta>0$ such that 
$P(w'_n(\delta) >\e)<\e$ for all $n$, where $w'_n(\delta)$
is the following modulus of continuity: $w'_n(\delta)=\inf_{\{t_i\}}
w_n(\{t_i\})$ 
where the infimum is over partitions $\{t_i\}$ of 
$[0,T]$ such that $t_i-t_{i-1}>\delta$ for all $i$, and 
$$
w_n(\{t_i\})=  \max_i
\sup\{ |\zeta_n(x,s)-\zeta_n(x,t)| : s,t\in [t_{i-1},t_i)\}.
$$
 (See \cite[Chapter 3]{billingsley} or \cite[Chapter 3]{ethierkurtz}.)

Now fix $y_0<y^-(x,\sigma)$, a constant $\alpha>0$,  the event
$$
A=\{ \mbox{$\zeta_0(y^+(x,\sigma))< \zeta_0(y)-\alpha$ for 
$y\in[y_0,y^-(x,\sigma)]$} \},
$$
and $\beta=P(A)>0$. The probability $P(A)$ is positive 
because $y^-(x,\sigma)<y^+(x,\sigma)$ by the assumption that 
$(x,\sigma)$ is a shock. Let $\e<(\alpha\wedge\beta)/8$, and
suppose there is a 
$\delta>0$ such that 
$P(w'_n(\delta) \ge\e)<\e$ for all $n$. Fix $\tau\in(\sigma,\sigma+\delta/2)$
so that $(x,\tau)$ is not a shock and so that 
$y(x,\tau)\in [y_0, y^-(x,\sigma)]$. This is possible 
because $t\mapsto y^-(x,t)$ is right-continuous and 
nonincreasing.  Let $F_n$ be the event
$$
F_n=\{ w'_n(\delta) < \e\,,\, |\zeta_n(x,t)-\zeta(x,t)|\le\e
\ \mbox{ for $t=\sigma,\tau$} \}.
$$
By (\ref{limznxtB}), $P(F_n)\ge 1-2\e$ for large enough $n$,
and then $P(A\cap F_n)\ge \beta/2>0$. Fix a sample point
$\omega\in A\cap F_n$. Fix a partition $\{t_i\}$ that achieves
$w_n(\{t_i\})<\e$ for this $\omega$. At this $\omega$, 
\beas
\zeta_n(x,\sigma)&\le& \zeta(x,\sigma)+\e\le \zeta_0(y^+(x,\sigma))+\e
\le \zeta_0(y(x,\tau))-\alpha+\e = \zeta(x,\tau)-\alpha+\e\\
&\le& \zeta_n(x,\tau)-3\alpha/4.
\eeas
This implies that $\sigma$ and $\tau$ cannot lie in the
same partition interval $[t_{i-1},t_i)$, so there must be at least
one partition point in $(\sigma,\tau]$. 
Since $\tau-\sigma<\delta/2$, there must be a unique partition 
point $t_k\in\{t_i\}\cap(\sigma,\tau]$. Then $w_n(\{t_i\})<\e$ 
forces $|\zeta_n(x,t)-\zeta_n(x,\sigma)|<\e$ for 
$t\in[t_{k-1}, t_k)$, while   
$|\zeta_n(x,t)-\zeta_n(x,\tau)|<\e$ for 
$t\in[t_k,t_{k+1})$. Combining this with the earlier inequality
gives
$$
\mbox{ $\zeta_n(x,t_k-)\le \zeta_n(x,t_k)-\alpha/2,$}
$$
which by the continuity of $u(x,t)$ implies 
$$
\mbox{ $z^n_{[nx]}(nt_k-)\le z^n_{[nx]}(nt_k)-n^{1/2}\alpha/2$}
$$ 
and contradicts the basic rule that  the particle
$z^n_{[nx]}(\cdot)$ jumps {\it leftward}. 

\subsection{Fluctuations for the conserved quantity}

We continue assuming that the process starts in local 
equilibrium according to assumptions  
(\ref{assrho}) and (\ref{assloqeq}).
In this section we consider the fluctuations
of the empirical density of the stick variables
$\{\eta^n_i(nt)\}$. 
Total stick length is conserved by the dynamics, as each 
jump of particle $z_i$ means that a random portion is subtracted
from  
$\eta_i$ and added on to $\eta_{i+1}$. Under
 assumptions (\ref{assrho}) and (\ref{assloqeq}) the empirical 
measure $n^{-1}\sum_i \eta^n_i(nt)\delta_{i/n}$ 
satisfies a hydrodynamic limit. 
Precisely, for any finite $a<b$, 
\be
\lim_{n\to\infty}\frac1n\sum_{i=[na]+1}^{[nb]}\eta^n_i(nt)
=\int_a^b \rho(x,t)dx
\qquad\mbox{a.s.}
\label{stickhydrolim}
\ee
See \cite{Seejp}. Actually only a limit in probability
is proved in \cite{Seejp}, but the result  can be strengthened
under assumption (\ref{assloqeq}). 

The next theorem is the fluctuation
theorem for this hydrodynamic limit. 
 The result is stated 
for the random distribution $\xi_n(t)$ defined below. 
First set
$$\rho^n_i(t)=n\int_{(i-1)/n}^{i/n}\rho(x,t)dx=n
u(i/n,t)-nu((i-1)/n,t).
$$
Then, 
for compactly supported test functions $\phi$,  define 
$$
\xi_n(t,\phi)=n^{-1/2}\sum_{i\in\ZZ}
\phi(i/n)\left( \eta^n_i(nt)-\rho^n_i(t)\right).
$$
Define another random distribution 
$\xi(t)$ by
\be
\xi(t,\phi)=-\int_{-\infty}^\infty \phi'(x)\zeta(x,t)dx,
\label{xidef}
\ee
where $\zeta(x,t)$ is defined by (\ref{zetadefB}) in terms of
the Brownian motion $B(\cdot)$. 

We want to put $\xi_n(t)$ and $\xi(t)$
 into some reasonable metric space, and a workable choice
turns out to be  the space
$\Hloc$ of distributions that are locally in $H^{-1}(\RR)$. 
To explain this we need some definitions. For the reader
unfamiliar with this,  Chapter 9
in \cite{folland} covers enough of the theory for following
our paper. 
Let $\cal D'$ be the space of distributions in Schwartz's
notation. Elements
$F\in\cal D'$ are linear functionals on the space 
$C_c^\infty(\RR)$ of compactly supported infinitely differentiable
functions, and they are continuous in this sense: 
$F(\phi_j)\to F(\phi)$ if all derivatives of $\phi_j$ converge 
uniformly to the corresponding derivatives of $\phi$, and 
all $\phi_j$ and $\phi$ are supported on a common compact set. 
Distributions can be multiplied by smooth functions: 
if $\chi$ is a $C^\infty$-function then the distribution 
$\chi F$ is defined by $\chi F(\phi)=F(\chi\phi)$. 

The Sobolev space $H^1(\RR)$ contains those $L^2$-functions $v$ that 
possess a weak derivative $v'$ in $L^2$.  It is a separable
 Hilbert space with (one possible) norm 
$$
\|v\|_{H^1(\RR)}=\|v\|_{L^2(\RR)}+\|v'\|_{L^2(\RR)}.
$$
$H^{-1}(\RR)$ is the dual space of $H^{1}(\RR)$, and itself
a separable Hilbert space. A continuous linear functional
on $H^1(\RR)$ also acts continuously on $C^\infty_c(\RR)$, and consequently
the elements of $H^{-1}(\RR)$ are from the space $\cal D'$.
Give $H^{-1}(\RR)$ the operator norm
$$
\|F\|_{H^{-1}(\RR)}=\sup\{ |F(v)|: \|v\|_{H^1(\RR)}\le 1\}.
$$
Now we can define the space of distributions locally in $H^{-1}$:
\be
\Hloc=\{ F\in {\cal D}': \mbox{$\chi F\in H^{-1}(\RR)$ for all 
$\chi\in C^\infty_c(\RR)$}
\}.
\label{Hlocdef}
\ee
Fix once and for all an increasing sequence of $C^\infty_c(\RR)$ functions
$\chi_k$ such that 
$${\bf 1}_{[-k+1,k-1]}\le \chi_k\le 
{\bf 1}_{(-k,k)}.
$$
  Then a distribution $F$ lies in 
$\Hloc$ iff $\chi_kF\in H^{-1}(\RR)$ for all $k$. We metrize
$\Hloc$ by
\be
R(F,G)=\sum_{k=1}^\infty 2^{-k} \left\{ 1\wedge 
\| \chi_kF-\chi_kG\|_{H^{-1}(\RR)}\right\}.
\label{metricdef}
\ee
Under this metric $\Hloc$ is a complete separable metric space.
 
We shall show that the process $t\mapsto \xi_n(t)$ is a 
random element of the  Skorokhod space 
$D([0,\infty),\Hloc)$, and that  $t\mapsto \xi(t)$ is a 
random element of the  space 
\break
$C([0,\infty),\Hloc)$ of continuous $\Hloc$-valued paths.

\begin{thm}  Assume {\rm (\ref{assrho})} and {\rm (\ref{assloqeq})}. 
Construct the processes $\{z^n(\cdot)\}$ on a common probability
space  with a two-sided Brownian motion $B(\cdot)$ so that the
conclusions of Theorem
\ref{thmloqeq} are valid. 
Fix a finite time horizon
 $\tau<\infty$. Then almost surely
\be
\lim_{n\to\infty} \sup_{t\in[0,\tau]}
R(\xi_n(t), \xi(t))=0.
\label{distrlim}
\ee
In particular, $\xi_n(\cdot)$ converges almost surely 
to $\xi(\cdot)$ on the path space $D([0,\infty),\Hloc)$. 
\label{stickthm}
\end{thm}

Theorem \ref{stickthm} is a corollary of Theorem \ref{thmloqeq}
and is valid under any hypotheses that make 
 Theorem \ref{thmloqeq} true. See Remark \ref{distremark}. 

To complement the theorem, we give alternative characterizations 
of the limiting distribution-valued process $\xi(\cdot)$. 
Spohn  \cite[page 260]{Sp} argued that the limiting fluctuations of 
an asymmetric conservative system should be governed by the
equation 
\be
\partial_t\xi+\partial_x[f'(\rho)\xi]=0.
\label{spohneqn}
\ee
By definition (\ref{xidef}), $\xi(t)=\partial_x\zeta(\cdot,t)$
in the 
distribution sense. Formally differentiating
through (\ref{transeqn}) with respect to $x$ then gives exactly
equation (\ref{spohneqn}). Thus we can regard $\xi(\cdot)$ as a 
distribution solution to (\ref{spohneqn}). 
Following (\ref{weaktranseqn}), the correct interpretation
 of the distribution 
$\partial_x[f'(\rho(\cdot,t))\xi(t)]$ is then,  
applied to a test function 
$\psi\in C^\infty_c(\RR)$,  
\be
\partial_x[f'(\rho(\cdot,t))\xi(t)](\psi)=
\int_{\RR}\zetabar(x,t)d[f'(\rho(\cdot,t))\psi'](x).
\label{weakspohneqn}
\ee

We can also consider $\xi$ as a Gaussian process indexed
by time and compactly supported test functions. For this 
we briefly introduce forward characteristics $w^\pm(a,t)$.
These are  inverse functions of $y^\pm(x,t)$ defined in 
(\ref{earlyypmdef}), themselves defined by 
$$
w^-(a,t)=\inf\{x: y^\pm(x,t)\ge a\}
\quad\mbox{and}\quad
w^+(a,t)=\sup\{x: y^\pm(x,t)\le a\}.
$$
We discuss these characteristics in Section \ref{charsection}. For now, 
we note that for a fixed $t$, $w^-(a,t)=w^+(a,t)$ for all but countably many points $a\in\RR$.
As functions of $t$, $w^\pm(a,\cdot)$ are the minimal and maximal 
Filippov solutions of the initial value problem 
\be
\frac{dx}{dt}=f'(\rho(x,t))\,,\quad x(0)=a.
\label{charode}
\ee
See \cite{dafermos} and \cite{Re} for more about this.

Ignoring the Lebesgue null set of shocks, we can write
$$
\xi(t,\phi)=-\int_{\RR}\phi'(x)
B\left(\int_0^{y^\pm(x,t)}\rho_0^2(r)dr \right)dx
$$
which shows that $\xi=\{\xi(t,\phi): t\in[0,\infty), \phi\in C^\infty_c(\RR)\}$
is a mean zero Gaussian process. Its distribution  is determined
by the correlations $E[\xi(s,\psi)\xi(t,\phi)]$, which we will show
in Section \ref{distrsect} to equal 
\be
E[\xi(s,\psi)\xi(t,\phi)]=\int_{\RR}\psi(w(r,s))\phi(w(r,t))\rho^2_0(r)dr. 
\label{corr1}
\ee
Here we wrote $w(r,t)$ for the a.e.\ defined function that agrees
with both $w^-(r,t)$ and $w^+(r,t)$ at a.e.\ $r$, for any fixed $t$. 

Correlations (\ref{corr1}) show that $\xi(t,\phi)$ can be 
equivalently described as follows. Fix a single two-sided Brownian
motion $W(\cdot)$. For $t\in[0,\infty)$ and $\phi\in C^\infty_c(\RR)$, define 
the random variables $\xitil(t,\phi)$ by the It\^o integrals
\be
\xitil(t,\phi)=\int_{\RR}\phi(w(r,t))\rho_0(r)dW(r).
\label{xitildef}
\ee
The function $\phi(w(r,t))$ is supported on some compact interval 
$a\le r\le b$ so there is no problem in defining the stochastic integral
(\ref{xitildef}) as a function of the increments $\{W(r)-W(a): a\le r\le b\}$.
The process $\xitil=\{\xitil(t,\phi): t\in[0,\infty), \phi\in C^\infty_c(\RR)\}$
 has the correlations given in (\ref{corr1}). From this we 
conclude that on the product space $\RR^{[0,\infty)\times C^\infty_c(\RR)}$
the distributions of $\xi$ and $\xitil$ are identical.

\subsection{Construction of the process and the variational coupling}
\label{construction}

The purpose of this section is mainly to establish the notation.
For more explanation and justification of this construction
we refer to \cite{alddiac}, 
\cite{Seejp}, \cite{Seldp}, \cite{Seperturb}. 
Consider a rate one, homogeneous
 Poisson
point process on $\RR\times(0,\infty)$.
 A sequence $(x_1,t_1)$, $(x_2,t_2)$, $\ldots$, 
$(x_m,t_m)$ of Poisson points is {\it increasing} if 
$$x_1<x_2<\cdots<x_m\qquad\mbox{and}\qquad t_1<t_2<\cdots<t_m\,.
$$
For $(a,s)$, $(b,t)\in\RR\times[0,\infty)$, let 
 $\mmL((a,s), (b,t))$ be the 
maximal number of  Poisson points on an increasing sequence
contained in
 $(a,b]\times(s,t]$. Abbreviate
$\mmL(b,t)=\mmL((0,0), (b,t))$. 

 Define an inverse to $\mmL$ by
$$\mmGa((a,s), m, \tau)=\inf\{ h>0: 
\mmL((a,s), (a+h,s+\tau)) \ge m\}\,.
$$
Again abbreviate
$\mmGa(m,\tau)=\mmGa((0,0), m,\tau)$. 
The well-known laws of large numbers are 
$$\lim_{s\to\infty} 
\frac1s\mmL(sb,st)= 2\sqrt{bt\,}\qquad\mbox{ and }\qquad
\lim_{s\to\infty} \frac1s\mmGa([sa],st)= \frac{a^2}{4t}\qquad\mbox{a.s.}
$$

Assume given a probability space $(\Omega, {\cal F}, P)$
 on which are defined the homogeneous 
Poisson point process on  $\RR\times(0,\infty)$ and  an initial 
configuration $(z_i(0):i\in\ZZ)$ for Hammersley's process. 
The process $z(t)=(z_k(t): k\in\ZZ)$ 
is defined by
\be
z_k(t)=\inf_{i:i\le k}\left\{ z_i(0)+
\mmGa((z_i(0) ,0), k-i, t)\right\}
\label{varcoup1}
\ee
for all $k\in\ZZ$ and $t>0$. 
Define the state space 
$${\cal Z}=\left\{ z=(z_i)\in\RR^{\ZZ}: \mbox{$z_{i-1}\le z_i$
for all $i$, and }
\lim_{i\to-\infty} i^{-2}z_i=0\right\}\,.
$$
If $(z_i(0))\in \cal Z$ a.s., then 
the infimum in (\ref{varcoup1}) is attained at some finite $i$
and $z(t)\in \cal Z$ for all
$t$ a.s. Thus  (\ref{varcoup1}) 
 defines a time-homogeneous Markov process $z(\cdot)$
with state space $\cal Z$. 

In this paper we work with a family of processes $\{z^n(\cdot)\}$. 
For each $n$ we assume the existence of some
probability space  $(\Omega,{\cal F},P)$
that supports the initial configuration 
 $z^n(0)=(z^n_i(0):i\in\ZZ)$ in addition to  
 the space-time Poisson
point process.
On this probability space  define the random variables
\be
\Gamma^{n,i}_m(t)=\mmGa((z^n_i(0),0), m,t)\,.
\label{Gandef}
\ee
Then, following (\ref{varcoup1}), the processes $\{z^n(t)\}$
are defined by 
$$z^n_k(t)=\inf_{i:i\le k}\left\{ z^n_i(0)+
\Gamma^{n,i}_{k-i}(t)\right\}\,.
$$

\section{Characteristics and the Hopf-Lax formula} 
\label{charsection}

The proofs of this paper take advantage of the 
correspondence between the macroscopic and microscopic
situations, and estimates on the probability that the
microscopic situation deviates from the macroscopic one.
First we  study the 
 macroscopic situation.
Without any additional trouble, we can relax the 
regularity assumption on 
 the initial interface $u_0$.  We adopt this
 standing assumption for this section:

\begin{ass}  $u_0$ is a nondecreasing, left-continuous
real-valued function on $\RR$ that satisfies the left growth
bound  {\rm (\ref{assu0})}.
\label{assu01}
\end{ass}

 We work throughout with the flux
 $f(\rho)=\rho^2$ with convex conjugate 
$g(x)=x^2/4$. Same results 
can be derived for any strictly convex, differentiable
conjugate pair $(f,g)$. The growth bound (\ref{assu0}) would
need to be
tailored to the $g$ in question. 

Under Assumption \ref{assu01} 
the function $\Phi(y)=u_0(y)+tg((x-y)/t)$, minimized
in the Hopf-Lax formula $(\ref{hopflax})$ over $y\in(-\infty,x]$, is 
lower semicontinuous and satisfies $\lim_{y\to-\infty}
\Phi(y)=\infty$. 
Consequently the minimum 
in  $(\ref{hopflax})$ is achieved at some point $y$, and the
set of minimizers is compact. 
 Define the function  $u(x,t)$ by 
the initial condition $u(x,0)=u_0(x)$ and by
(\ref{hopflax}).  Then  for a fixed $t>0$, $u(\cdot,t)$  
 is locally Lipschitz in $x$ (we check this
below), and  $x^{-2}u(x,t)\to 0$ as 
$x\to-\infty$. 
The Hopf-Lax formula can be iterated as a semigroup: 
\be
u(x,t)=\inf_{y:y\le x}\left\{ u(y,s)+(t-s)g\left(\frac{x-y}{t-s}\right)
\right\}
\label{hopflaxst}
\ee
for all $0<s<t$ and $x\in\RR$.  Define the set of minimizers
 in (\ref{hopflaxst}) by 
\be
I(x;s,t)=\left\{ y\le x: 
u(x,t)=u(y,s)+(t-s)g\left(\frac{x-y}{t-s}\right)
\right\}. 
\label{Ixtdef}
\ee
$I(x;s,t)$ is nonempty and compact. Define 
 minimal and maximal
minimizers by
\be
y^-(x;s,t)=\inf I(x;s,t) \quad
\mbox{ and }\quad
y^+(x;s,t)=\sup I(x;s,t).
\label{ypmdef}
\ee
The following properties can be checked: If $x_1<x_2$ then
 $y^+(x_1;s,t)\le y^-(x_2;s,t)$,
while if $t_1<t_2$ then $y^+(x;s,t_2)\le y^-(x;s,t_1)$.
 $y^\pm(x;s,t)$
is nondecreasing in $x$ and nonincreasing in $t$. 
$y^+$ is right- and $y^-$ left-continuous in $x$, 
while $y^+$ is left- and $y^-$ right-continuous in $t$. 
Consequently, for fixed $s<t$,  $y^\pm(\cdot\,;s,t)$
have the same continuity points, they coincide
on these continuity points, and  
$y^-(x;s,t)<y^+(x;s,t)$ iff $x$ is a discontinuity point. 
A  similar statement holds for $y^\pm(x;s,\cdot)$
as a function of $t$, for fixed $x,s$. 

Next define minimal and maximal forward characteristics by 
\be
w^-(a;s,t)=\sup\{ x: y^\pm(x;s,t)< a\}
= \inf\{ x: y^\pm(x;s,t)\ge a\} 
\label{x-def}
\ee
and
\be
w^+(a;s,t)=\sup\{ x: y^\pm(x;s,t)\le a\}
= \inf\{ x: y^\pm(x;s,t)> a\}. 
\label{x+def}
\ee
The equalities between the alternative definitions
follow from the properties of \break 
$y^\pm(x;s,t)$. 
For $w^\pm(a;s,t)$ we have these properties:
nondecreasing in $a$, nondecreasing in $t$, 
$w^+$ is right- and $w^-$ left-continuous in $a$, 
$w^+(a_1;s,t)\le w^-(a_2;s,t)$ for $a_1<a_2$. As above,
for fixed $s<t$,  $w^\pm(\cdot\,;s,t)$
have the same points of continuity, coincide
on  continuity points, and  
$w^-(a;s,t)<w^+(a;s,t)$ iff $a$ is a discontinuity point.
Note the equivalence
\be
y^-(x;s,t)\le a\le y^+(x;s,t)\Longleftrightarrow
w^-(a;s,t)\le x\le w^+(a;s,t).
\label{xyequiv}
\ee
Note also that as a trivial consequence of the definitions, 
$y^-(x;s,t)\le y^+(x;s,t)\le x$, and 
$a\le w^-(a;s,t)\le w^+(a;s,t)$. 

We adopt the following notational conventions. 
When the $\pm$ functions coincide we write 
$y^\pm(x;s,t)=y(x;s,t)$ and $w^\pm(a;s,t)=w(a;s,t)$. 
When $s=0$ abbreviate $y^\pm(x;0,t)=y^\pm(x,t)$ and similarly
$y(x,t)$, $w^\pm(a,t)$, $w(a,t)$. 

As mentioned earlier, $u(x,t)$ is the unique viscosity
solution of the Hamilton-Jacobi equation 
$u_t+f(u_x)=0$ with $f(\rho)=\rho^2$ and initial data
$u|_{t=0}=u_0$. Set 
$b(x)=g'(x)=x/2$, and define two functions $\rho^\pm(x,t)$ by
\be
\rho^\pm(x,t)=b\left( \frac{x-y^\pm(x,t)}t\right)\qquad
\mbox{for $(x,t)\in\RR\times(0,\infty)$.}
\label{laxoleinik}
\ee
For a fixed $t>0$, $\rho^\pm$ give the one-sided
$x$-derivatives of $u$:
$$
\rho^\pm(x,t)=\lim_{\e\to 0^\pm}\frac{u(x+\e,t)-u(x,t)}\e.
$$
There is a 
 function $\rho$ such that 
 $\rho(x,t)=\rho^\pm(x,t)$ for all but
countably many $x$, because $y^-(x,t)=y^+(x,t)$ for all but
countably many $x$ (for fixed $t$ again). 
The a.e.\ defined function $\rho(x,t)$ 
is the unique entropy
solution of the Burgers equation 
$\rho_t+f(\rho)_x=0$ with initial condition given
by the Radon measure $du_0(x)$. More precisely, we mean that
$\rho(x,t)$ is a weak solution in this integral sense:
for all $\phi\in C^\infty_c(\RR)$, 
\be
\int_\RR \phi(x)\rho(x,t)dx-\int_\RR \phi(x) du_0(x) =
\int_0^t \int_\RR \phi'(x)f(\rho(x,s))dx ds.
\label{weakburgers}
\ee
Formula (\ref{laxoleinik})
 is known as the  Lax-Oleinik formula.
See  \cite{Ev} for the textbook
p.d.e.\ theory. The appendix in \cite{Seejp} develops a 
uniqueness theory for $\rho(x,t)$ when the initial 
condition $du_0$ is a measure with singularities.

The next lemma collects some properties proved in 
Section 3 of Rezakhanlou \cite{Re}.
 In that paper the initial density profile $\rho_0=u_0'$ is assumed 
bounded and integrable, but the proofs
work with at most minor modifications  under our
Assumption \ref{assu01}. 
 It is worthwile to 
note that strict convexity 
and continuous differentiability 
of $g(x)=x^2/4$ are critical for many of the good
properties of the characteristics utilized in this
section. The reader can compare with \cite{SepII} where
the $g$ function corresponding to the $K$-exclusion process 
  is not known to possess these properties. 

\begin{lm} {\rm (a)} Suppose $0\le t_1<t_2<t_3$, $y_2\in I(x;t_2,t_3)$,
and $y_1\in I(y_2;t_1,t_2)$. Then $y_1\in I(x;t_1,t_3)$ 
and  
$$\frac{x-y_1}{t_3-t_1}=\frac{x-y_2}{t_3-t_2}=\frac{y_2-y_1}{t_2-t_1}.
$$
In other words, the points $(y_1,t_1)$, $(y_2,t_2)$, and 
$(x,t_3)$ lie on a line segment. 

{\rm (b)} For $0\le s< s_1<t$, 
$$y^\pm(x;s_1,t)=\frac{s_1-s}{t-s}x+
\frac{t-s_1}{t-s} y^\pm(x;s,t).$$

{\rm (c)} For $0<s<t$ and all $a\in\RR$, $w^\pm(a;s,t)=w(a;s,t)$. For $s=0$
and $a\in\RR$, 
 $w^\pm(a,t)=w(a,t)$ is guaranteed by 
\be
\liminf_{\e\searrow 0}\frac{u_0(a+\e)-u_0(a)}\e
\le
\limsup_{\e\searrow 0}\frac{u_0(a)-u_0(a-\e)}\e. 
\label{fancond1}
\ee

{\rm (d)} Suppose $0\le t_1<t_2<t_3$. Then $w^\pm(a;t_1,t_3)
=w(w^\pm(a;t_1,t_2);t_2,t_3)$.
\label{charlm1}
\end{lm}

Statements (a) and (b) can be  augmented as follows.

\begin{lm} Let $y_1\in I(x_1,t_1)$. Let
$z(t)=(t/t_1)x_1+\left(1-(t/t_1)\right)y_1$, $t\in[0,t_1]$,
be the line segment from $(y_1,0)$ to $(x_1,t_1)$.
Then 

{\rm (i)} $z(t)\in I(x_1;t,t_1)$ for each $t\in[0,t_1)$, 
and 

{\rm (ii)} $I(z(t);s,t)=\{z(s)\}$ for all $0\le s<t<t_1$. 
\label{charlm2}
\end{lm}

{\it Proof.} {\bf Step 1:} we show that $I(z(t);0,t)=\{y_1\}$
 for all $0<t<t_1$. The assumption $y_1\in I(x_1,t_1)$
implies that 
$$
u_0(y_1)+t_1g\left(\frac{x_1-y_1}{t_1}\right)
\le 
u_0(y)+t_1g\left(\frac{x_1-y}{t_1}\right)
\quad\mbox{ for all $y\le x_1$.}
$$
 This rearranges to give 
\be
\frac{u_0(y_1)-u_0(y)}{y_1-y}\le 
\frac{g\left(\frac{x_1-y}{t_1}\right)
-g\left(\frac{x_1-y_1}{t_1}\right)}
{\frac{x_1-y}{t_1}-\frac{x_1-y_1}{t_1}}
\qquad \mbox{for all $y<y_1$}
\label{ylty1ineq}
\ee
and 
\be
\frac{u_0(y)-u_0(y_1)}{y-y_1}\ge 
\frac{g\left(\frac{x_1-y_1}{t_1}\right)
-g\left(\frac{x_1-y}{t_1}\right)}
{\frac{x_1-y_1}{t_1}-\frac{x_1-y}{t_1}}
\qquad \mbox{for all $y_1<y\le x$.}
\label{ygty1ineq}
\ee
Now let $y\in(y_1,z(t)]$. By the definition of $z(t)$, 
$$
\frac{z(t)-y}{t}<\frac{x_1-y}{t_1}<\frac{z(t)-y_1}{t}=\frac{x_1-y_1}{t_1}.
$$
By the strict convexity of 
$g(x)=x^2/4$, 
$$
\frac{g\left(\frac{x_1-y_1}{t_1}\right)
-g\left(\frac{x_1-y}{t_1}\right)}
{\frac{x_1-y_1}{t_1}-\frac{x_1-y}{t_1}}
>
\frac{g\left(\frac{z(t)-y_1}{t}\right)
-g\left(\frac{z(t)-y}{t}\right)}
{\frac{z(t)-y_1}{t}-\frac{z(t)-y}{t} }.
$$
This combined with (\ref{ygty1ineq}) gives
$$
u_0(y)-u_0(y_1)>g\left(\frac{z(t)-y_1}{t}\right)
-g\left(\frac{z(t)-y}{t}\right)
$$
which implies that no $y>y_1$ can be in the minimizing set
$I(z(t);0,t)$. A similar argument that utilizes (\ref{ylty1ineq})
rules out $y<y_1$, and Step 1 is complete. 

 {\bf Step 2:} we show $z(t)\in I(x_1;t,t_1)$. 
\beas
u(x_1,t_1)&=&u_0(y_1)+t_1g\left(\frac{x_1-y_1}{t_1}\right)\\
&=&u_0(y_1)+ tg\left(\frac{z(t)-y_1}{t}\right)
+ (t_1-t)g\left(\frac{x_1-z(t)}{t_1-t}\right)\\
&=& u(z(t),t)+ (t_1-t)g\left(\frac{x_1-z(t)}{t_1-t}\right)
\eeas
which implies the conclusion. Above we used the 
line segment assumption in the form 
$$
\frac{x_1-y_1}{t_1}=
\frac{z(t)-y_1}{t}=\frac{x_1-z(t)}{t_1-t}
$$
and then Step 1. 
  
 {\bf Step 3:} It remains to show  $I(z(t);s,t)=\{z(s)\}$ 
for $0< s<t<t_1$. Now we know $z(s)\in I(x_1;s,t_1)$ 
by Step 2, so we can simply repeat Step 1 for $s>0$ 
in place of $s=0$. 
\qed

We emphasize the conclusion of part (ii) of the last lemma:
Along the line segment $z(t)$, $0\le t<t_1$, Hopf-Lax
minimizers are unique.  

\begin{lm} 
{\rm (i)} Let $0<s<t$ and $x_1=w(x_0;s,t)$. Then 
\be
y^-(x_1,t)\le y^-(x_0,s)\le y^+(x_0,s)\le y^+(x_1,t).
\label{yminusineq1}
\ee
Conversely, if {\rm (\ref{yminusineq1})}
 holds and the middle inequality is strict, then $x_1=w(x_0;s,t)$.

{\rm (ii)} Let $(x,t)$ and $(x_1,t_1)$ be arbitrary points
in $\RR\times(0,\infty)$ with $t\le t_1$. Suppose  the open intervals
$J_{x,t}=(y^-(x,t),y^+(x,t))$ and $J_{x_1,t_1}=(y^-(x_1,t_1),y^+(x_1,t_1))$ 
are nonempty. Then one of two cases happens: either the intervals
 are disjoint, which happens if $t=t_1$ and $x\ne x_1$, or if $t<t_1$ and 
$x_1\ne w(x;t,t_1)$. Or $J_{x,t}\subseteq J_{x_1,t_1}$ which happens 
if  $(x,t)=(x_1,t_1)$ or if $t<t_1$ and 
$x_1= w(x;t,t_1)$.

{\rm (iii)} Let $y\in I(x,t)$ and 
$z(s)=(s/t)x+\left(1-(s/t)\right)y$, $s\in[0,t]$,
be the line segment from $(y,0)$ to $(x,t)$.
Suppose $0\le s_1<s_2\le t$ and 
$w^-(x_1;s_1,s_2)\le x_2\le w^+(x_1;s_1,s_2).$
Then the points $(x_1,s_1)$ and 
$(x_2,s_2)$ must be on the same side of the line
segment $z(\cdot)$. In other words, 
$x_1< z(s_1)$ implies $x_2\le z(s_2)$, and 
$x_1> z(s_1)$ implies $x_2\ge z(s_2)$. 
\label{charlm3}
\end{lm}

{\it Proof.} (i) By (\ref{xyequiv}),
$$y'\equiv y^-(x_1;s,t)\le x_0\le y^+(x_1;s,t)\equiv y''.$$
By Lemma \ref{charlm1}(a), $y^-(y';0,s), y^+(y'';0,s)\in
I(x_1,t)$. By the monotonicity of $y^\pm(\cdot\,,s)$,
$$y^-(x_1;0,t)\le y^-(y';0,s)\le  y^-(x_0;0,s)
\le y^+(x_0;0,s)\le y^+(y'';0,s)\le y^+(x_1;0,t).$$

For the converse part, if $x_1>w(x_0;s,t)$ then by monotonicity
and the part already proved,
$$y^-(x_1;0,t)\ge  y^+(w(x_0;s,t) ;0,t) \ge y^+(x_0; 0,s).$$
This contradicts  (\ref{yminusineq1}) if the middle
inequality of (\ref{yminusineq1}) 
is strict. Similarly rule out the case $x_1<w(x_0;s,t)$.

(ii) 
  If $t=t_1$ then either $x=x_1$ or 
the intervals must be disjoint, because $x<x_1$ implies
$y^+(x,t)\le y^-(x_1,t)$. 
Suppose $t_1 >t$. If $x_1=w(x;t,t_1)$ then by part (i) 
$(y^-(x,t),y^+(x,t))$ is contained in $(y^-(x_1,t_1),y^+(x_1,t_1))$.
 On the other hand, if $x_1>w(x;t,t_1)$ then we have 
disjointness by the argument already used in part (i): 
$$y^-(x_1;0,t_1)\ge  y^+(w(x;t,t_1) ;0,t_1) \ge y^+(x; 0,t).$$
Similar for the remaining cases.

(iii) Let us show $x_1> z(s_1)$ implies $x_2\ge z(s_2)$. 
By Lemma \ref{charlm2}, $y^-(z(s_2);s_1,s_2)\le z(s_1)<x_1$, 
so
$$x_2\ge w^-(x_1;s_1,s_2)=\sup\{\xi: y^-(\xi;s_1,s_2) <x_1\}
\ge z(s_2). 
\qed
$$

Part (i) of the previous lemma has the following 
meaning. We say that $(x,t)$ is a {\it shock} if 
$y^-(x,t)<y^+(x,t)$. By (\ref{laxoleinik}), this is 
the same as saying that the Lax-Oleinik solution
 $\rho(\cdot,t)$ is discontinuous at
$x$. In fact this is the same as saying that $\rho$ is 
continuous at $(x,t)$. For  when the minimizer
$y(x,t)=y^\pm(x,t)$ is unique and $(x_j,t_j)\to (x,t)$, 
then any choice $y_j\in I(x_j,t_j)$ satisfies
$y_j\to y(x,t)$.  Inequalities (\ref{yminusineq1}) imply that
once a shock is created, it moves along a forward
characteristic and never disappears. Note though that 
shocks merge when characteristics merge, so the number
of shocks may decrease. 

Next we look at the continuity of characteristics and 
$u(x,t)$. 

\begin{lm} {\rm  (a)} Given $a<b$ and $T>0$, there exists 
a constant $C$ such that $I(x,t)\subseteq[x-Ct^{1/2},x]$
 for all $x\in[a,b]$ and $t\in(0,T]$. 

{\rm  (b)} Fix $0\le s<T<\infty$ and suppose 
the function $u(\cdot,s)$ on the
right-hand side of {\rm (\ref{hopflaxst})} is locally Lipschitz in
the  $x$-variable. Fix $a<b$, and let $K$ be the Lipschitz
constant of $u(\cdot,s)$ on the interval $[y^-(a;s,T),b]$. 
Then for all $x\in[a,b]$ and $t\in(s,T]$, 
$I(x; s,t)\subseteq[x-4K(t-s),x]$. 

{\rm (c)} Assume $u_0$ is locally Lipschitz, in addition to
assumption
{\rm (\ref{assu0})}. Then $u(x,t)$ is locally 
Lipschitz on $\RR\times[0,\infty)$. 
For any $a<b$ and $T<\infty$ there exists a constant
$L=L(a,b,T)$ such that 
$I(x; s,t)\subseteq[x-L(t-s),x]$ for all $x\in[a,b]$ and
$0\le s<t\le T$.

{\rm (d)} Under the original assumptions of left-continuity
and {\rm (\ref{assu0})}, $u(x,t)$ is locally Lipschitz on
$\RR\times(0,\infty)$, and $\lim_{t\to 0}u(x,t)=u_0(x)$
for all $x\in\RR$.
\label{Ixtlm}
\end{lm}

{\it Proof.} (a) Pick $\e>0$ small enough so that
$\sqrt{T\e\,}<1/2$. Pick $M$ so that
$|u_0(y)|\le \e y^2$ for  $y\le M$.  
Let $(x,t)\in [a,b]\times(0,T]$. Any 
$y\in I(x,t)$ must satisfy $y\le x$ and 
$u_0(x)\ge u_0(y)+(x-y)^2/(4t)$, from which 
follows 
$$(x-y)^2\le 4t(u_0(x)-u_0(y))\le 4t( C_1+\e y^2)
$$
where we picked $C_1\ge 0$ so that $2|u_0|\le C_1$
on  $[M,b]$. Squareroots and 
algebra give
$$x-2\sqrt{C_1t\,} \le y+2|y|\sqrt{t\e\,}
 \le y+2\sqrt{t\e\,}(|b|+|y^-(a,T)|), $$ 
from which the conclusion follows.

 (b) Let $c=y^-(a;s,T)$. Since $y^\pm(x;s,t)$ is
nonincreasing in $t$, 
$I(x;s,t)\subseteq [c,b]$  for all 
$(x,t)\in [a,b]\times(s,T]$.  Any 
$y\in I(x;s,t)$ must satisfy $y\le x$ and 
$u(x,s)\ge u(y,s)+(x-y)^2/(4(t-s))$, from which 
follows 
$$(x-y)^2\le 4(t-s)(u(x,s)-u(y,s))\le  4(t-s)K(x-y).
$$

(c) First step: to show that on a bounded interval
$[a,b]$,  
$u(\cdot, t)$ is locally Lipschitz in the 
$x$-variable with Lipschitz constant independent
of $t\in[0,T]$. 
Let $K$ be the Lipschitz constant of $u_0$
on $[y^-(a,T),b]$. Let $x_1<x_2$ in $[a,b]$ and pick
$y_1\in I(x_1,t)$.  Then
$y_1\in [y^-(x_1,t),x_1]\subseteq [y^-(a,T),b]$. 
Set $y_2=y_1+x_2-x_1\in [y_1,x_2]\subseteq [y^-(a,T),b]$. 
Then
\beas
0&\le& u(x_2,t)-u(x_1,t)\\
&\le& u_0(y_2)+tg\left(\frac{x_2-y_2}t\right)
-u_0(y_1)-tg\left(\frac{x_1-y_1}t\right)\\
&=&u_0(y_2)-u_0(y_1)\le K|y_2-y_1|=K|x_2-x_1|.
\eeas

Second step: to show the existence of 
a constant $C$ such that  for all $x\in[a,b]$ 
and $0\le s<t\le T$,  
$|u(x,t)-u(x,s)|\le C|t-s|$. The two steps together 
imply that $u$ is locally Lipschitz. 

For the second step, apply the first step
to let $K$ be the common Lipschitz
constant for the functions $\{u(\cdot,t):0\le t\le T\}$ 
on the $x$-interval
$y^-(a,T)\le x\le b$. 
Let $y= y^-(x;s,t)\in I(x;s,t)$. By part (b) of this lemma, 
$|x-y|\le 4K(t-s)$. Furthermore, by Lemma \ref{charlm1}(b),
$$y= y^-(x;s,t)=\frac{s}{t}x+\left(1-\frac{s}{t}\right) y^-(x,t)
\in [y^-(a,T), b].$$
Now we may reason as follows for $x\in[a,b]$ 
and $0\le s<t\le T$: 
\beas
0&\le& u(x,s)-u(x,t)=u(x,s)-u(y,s)-\frac{(x-y)^2}{4(t-s)}\\
&\le& u(x,s)-u(y,s)\le K|x-y|\le 4K^2(t-s).  
\eeas
We may take $C=4K^2$ and the proof  of Lipshcitz continuity  is complete. 

  By Lemma \ref{charlm1}(b), $[y^-(a;s,T),b]\subseteq
[y^-(a,T), b]$ for all $0\le s<T$, so 
a single Lipschitz constant 
 $K$ works for all $s$ in part (b). 

(d) For local Lipschitz continuity of
$u$ on $\RR\times(0,\infty)$
 it suffices to show that $u(\cdot,s)$ is Lipschitz in $x$
for any fixed $s>0$, for then we can apply part (c) to the 
solution obtained for $t\ge s$. Let $x_1<x_2$ in $[a,b]$ and  
$y\in I(x_1,s)$. 
\beas
0&\le& u(x_2,s)-u(x_1,s)\le u_0(y)+sg((x_2-y)/s)-u_0(y)-sg((x_1-y)/s)\\
&=&sg((x_2-y)/s)-sg((x_1-y)/s)=g'(\xi)(x_2-x_1)\le C(x_2-x_1)
\eeas
where we used the mean value theorem, and chose $C$ as an upper
bound for $g'$ on the interval $[0, s^{-1}b-s^{-1}y^-(a,s)]$.

The limit $u(x,0+)=u_0(x)$ follows from 
$u_0(y^+(x,t))\le u(x,t)\le u_0(x)$, part (a), and left-continuity
of $u_0$. 
\qed

Let $0\le t_0<T$. 
A {\it forward characteristic} emanating from $(x_0,t_0)$
 is any
function $r(t)$, $t_0\le t\le T$, that satisfies $r(t_0)=x_0$ and 
\be
\mbox{$w^-(r(s);s,t)\le r(t)\le w^+(r(s);s,t)$ 
for all $t_0\le s<t\le T$.}
\label{rtdef}
\ee
Notice that this implies  $r(t)=w(r(s);s,t)$ for all $0<s<t$.
Multiple forward characteristics can emanate
only from a point $(a,0)$ on the $t=0$ line
[and only if (\ref{fancond1}) fails]. For example, 
$w^\pm(a,t)$ are forward characteristics that emanate 
from $(a,0)$. 

\begin{lm} {\rm (i)} A forward characteristic $r(\cdot)$ is 
absolutely continuous on $[0,T]$ for any $T<\infty$, and locally
Lipschitz continuous on $(0,\infty)$. If $u_0$ is locally
Lipschitz, then  $r(\cdot)$ is 
 locally Lipschitz on $[0,\infty)$.

{\rm (ii)} Let $A\subseteq\RR\times[0,\infty)$
be any set of points such that no two points of $A$ have 
the same $t$-coordinate, and for every 
pair $(x_1,t_1),(x_2,t_2)\in A$, if $t_1<t_2$ then 
$w^-(x_1;t_1,t_2)\le x_2\le w^+(x_1;t_1,t_2)$. (If $t_1>0$
it follows  
we must have equality $x_2= w(x_1;t_1,t_2)$.) 
Then there exists a forward characteristic $r(t)$, $t\ge 0$,
such that all points of $A$ lie on $r(\cdot)$. 
\label{rlm1}
\end{lm}

{\it Proof.} Part (i) is a consequence of Lemma \ref{Ixtlm}. 

For part (ii), let 
$$
\tau=\inf\{ t: \mbox{there exists $x\in\RR$ such that $(x,t)\in A$}\}.
$$
Let $\xi$ be such that $(\xi,\tau)\in A$
if such a point exists. If not, pick a sequence
$(x_n,t_n)\in A$ such that $t_n\searrow \tau$. 
By assumption $x_n=w(x_{n+1};t_{n+1},t_n)$ for all $n$
(note that now $t_n>0$ for all $n$) so 
that $x_n\ge x_{n+1}$. Thus there is a limit 
$x_n\searrow\xi$. This limit must be finite because
$x_1=w(x_{n};t_{n},t_1)$ implies 
$x_n\ge y^-(x_1;t_n,t_1)\ge y^-(x_1,t_1)$. 

Let us first show that for all $(x,t)\in A$ with $t>\tau$, 
\be
w^-(\xi;\tau,t)\le x\le w^+(\xi;\tau,t).
\label{x+-xiineq}
\ee
If $(\xi,\tau)\in A$ then (\ref{x+-xiineq}) is part of the assumption.
Otherwise, for large $n$ so that $t_n<t$, the assumption gives
$x=w(x_{n};t_{n},t)$ from which follows
$$y^-(x;t_n,t)\le x_n\le y^+(x;t_n,t).$$
Let $n\to\infty$ and apply Lemma \ref{charlm1}(b) to get
$$y^-(x;\tau,t)\le \xi\le y^+(x;\tau,t)$$
which implies (\ref{x+-xiineq}). 

Suppose first $\tau>0$. Pick any $y\in I(\xi,\tau)$ and 
define $r(\cdot)$ by
$$
r(t)=
\left\{
\begin{array}{ll}
(t/\tau)\xi +(1-(t/\tau))y, & t\in[0,\tau]\\
w(\xi;\tau,t),& t\in(\tau,\infty).
\end{array}
\right.
$$
Properties (\ref{rtdef}) are satisfied by Lemmas \ref{charlm1}(d)
and  \ref{charlm2}.

If $\tau=0$ set $r(0)=\xi$ and for every $s>0$
set $r(s)=w(x;t,s)$ where $(x,t)$ is an arbitrary
point in $A$ such that $0<t<s$. Such points exist
since $\tau=0$.  $r(s)$ is well
defined because if $t_1<t_2<s$ and 
 $(x_1,t_1),(x_2,t_2)\in A$, then 
 $w(x_2;t_2,s)= w(w(x_1;t_1,t_2);t_2,s)=w(x_1;t_1,s) $
by assumption and by Lemma \ref{charlm1}(d).
This same Lemma \ref{charlm1}(d) implies  that 
$r(\cdot)$ satisfies the defining properties
(\ref{rtdef}). 
\qed

By Lemma \ref{rlm1}(i), $r'(t)$ exists at a.e.\ $t$, and 
$r$ is the integral of its derivative. Next, a formula
for the derivative. Let
\be
h(x,t)=\left\{
\begin{array}{ll}
f'(\rho(x,t)), &\mbox{if $y^-(x,t)=y^+(x,t)$} \smallskip\\
\displaystyle\frac{f(\rho^+(x,t))-f(\rho^-(x,t))}{\rho^+(x,t)-\rho^-(x,t)}, 
&\mbox{if $y^-(x,t)<y^+(x,t)$.}
\end{array}\right. 
\label{hdef}
\ee

\begin{thm} For any forward characteristic $r(\cdot)$, 
$r'(t)=h(r(t),t)$ for Lebesgue a.e.\ $t$. 
\label{odethm}
\end{thm}

{\it Proof.} Let $t_1>t_0>0$ and $x_1=w(x_0;t_0,t_1)$. Let
$r(t)=w(x_0;t_0,t)$ be the forward characteristic emanating from $(x_0,t_0)$. 

{\it Case 1.} If $(x_0,t_0)$ and
$(x_1,t_1)$ are not shocks, the proof on p.\ 136--137 of \cite{Re}
shows that $(x_1-x_0)/(t_1-t_0)=f'(\rho(x_0,t_0))=f'(\rho(x_1,t_1))$. 

{\it Case 2.}
Suppose $(x_0,t_0)$ is a shock. Abbreviate $\rho^\pm=\rho^\pm(x_0,t_0)$.
 Let $\e>0$.  Take $t_1$ close enough to
$t_0$ and $a<x_0<b$ close enough to $x_0$ so that, for
$(x,t)\in[a,b]\times[t_0,t_1]$, 
\beas
|\rho^\pm(x,t)-\rho^-|<\e &&\mbox{if $x<r(t)$, and }\\
|\rho^\pm(x,t)-\rho^+|<\e &&\mbox{if $x>r(t)$.}
\eeas
This can be achieved because $\rho^\pm(x,t)=b((x-y^\pm(x,t))/t)$,
$b=g'$ is continuous, 
and $y^\pm(x,t)\to y^-(x_0,t_0)$ as $(x,t)$ approaches $(x_0,t_0)$
with $t\ge t_0$ and  $x<r(t)$. Decrease $t_1$ further towards 
$t_0$ so that $a<x_0\le r(t)\le x_1<b$ for $t\in[t_0,t_1]$. 

Now apply (\ref{weakburgers}) to a test function $\phi\in C^\infty_c(\RR)$ compactly
supported inside $(a,b)$ and such that $\phi\equiv 1$ on $[x_0,x_1]$,
$\phi'\ge 0$ on $(a,x_1]$ and $\phi'\le 0$ on $[x_1,b)$.  
After a calculation this gives
$$
\frac{x_1-x_0}{t_1-t_0}=\frac{f(\rho^+)-f(\rho^-)}{\rho^+-\rho^-}+O(\e)
$$
provided the ratio $(b-a)/(t_1-t_0)$ is kept bounded. Letting
$t_1\searrow t_0$  shows 
that the right derivative 
$r'(t_0+)$ at $t_0$ is given by $h(r(t_0),t_0)$. 

The only type of point $(x_0,t_0)$ not covered by Cases 1 and 2 
is such that $(x_0,t_0)$ is not a shock but $(x_1,t_1)$ is a shock 
for every $t_1>t_0$. There can be 
only one such point on any forward characteristic. 
These two calculations  suffice to prove the theorem, because 
from absolute continuity we know $r'(t)$ exists almost
everywhere. 
\qed

From this theorem one concludes
 a result of  Dafermos \cite{dafermos} and 
Rezakhanlou \cite{Re}: A forward characteristic is a Filippov solution
of the initial value problem (\ref{charode}).
  We shall not discuss this
point further as we make no use of it. But we will use the formula 
$(d/dt)w^\pm(q,t)=h(w^\pm(q,t),t)$.

At certain stages in the proofs we need
 to know  that an expression such as 
$$\sup_{x\in[a,b], t\in[0,\tau]}\left| \zeta_n(x,t)-\inf_{y\in I(x,t)}
\zeta_n(y,0)\right|$$
 is a measurable   random variable. 
We want  to argue that the supremum can 
be taken over a countable set $\{(x_k,t_k)\}$. 
Such an attempt will
reveal 
as troublesome those points $(x,t)$ for 
which the Hopf-Lax formula (\ref{hopflax}) possesses
more than two minimizers. [In other words, 
$I(x,t)$ does not coincide with 
$\{y^\pm(x,t)\}$, which is either a singleton or a two-point
set.] 
We show here that such points are at most 
 countable.

\begin{thm} Let $T>0$, and 
$${\cal U}=\{(x,t)\in\RR\times(0,T]: 
I(x,t)\ne \{y^\pm(x,t)\} \}.
$$
The set $\cal U$ is countable. 
\label{Uthm}
\end{thm}

%
%

{\it Proof.} If $(x,t)\in{\cal U}$, then necessarily 
$y^+(x,t)-y^-(x,t)>0$.
 Thus it suffices to prove
the countability of the set
$${\cal U}_\alpha=\{(x,t)\in {\cal U}: 
y^+(x,t)-y^-(x,t)\ge\alpha   \}
$$
for an arbitrary $\alpha>0$.
The countability of ${\cal U}_\alpha$ will be achieved by
showing that on any bounded set $[a,b]\times(0,T]$,
the points of ${\cal U}_\alpha$ lie on a finite collection
of forward characteristics, and each characteristic 
contains at most countably many $\cal U_\alpha$-points.

\hbox{}

{\bf Step 1.} Fix a forward characteristic
$r(\cdot)$. We show that 
$r(\cdot)$ contains  at most countably many points from
$\cal U$, so in particular, at most countably many points from
$\cal U_\alpha$. We do this by associating to each
point $(r(t),t)\in\cal U$ a nonempty open interval $J_t$
so that these $J_t$'s are pairwise disjoint. 

Let $t_1>0$ be such that $(r(t_1),t_1)\in\cal U$. This
implies that the minimizing set $I(r(t_1),t_1)$ contains
at least three points
$y^-(r(t_1),t_1)<y_1<y^+(r(t_1),t_1)$. Let 
$\zeta=\{(z(t),t): 0\le t\le t_1\}$ be the line segment 
from $(y_1,0)$ to $(r(t_1),t_1)$. By Lemma \ref{charlm3},
the curve $\{(r(t),t): 0\le t\le t_1\}$ must lie
entirely on one side of $\zeta$. 
Also, no $(z(t),t)\in\zeta$ lies in $\cal U$
for $t\in(0,t_1)$. This is because by Lemma \ref{charlm2}(ii)
the minimizer set $I(z(t),t)$ is the singleton $\{y_1\}$,
while every $\cal U$-point has multiple minimizers. 
Thus one of these cases holds: 
$$
\begin{array}{ll}
&\mbox{{\it Case I:} $r(t)<z(t)$ for all $t\in(0,t_1)$ such that
$(r(t),t)\in\cal U$.}\\
&\mbox{{\it Case II:} $r(t)>z(t)$ for all $t\in(0,t_1)$ such that
$(r(t),t)\in\cal U$.} 
\end{array}
$$
In {\it Case I} set $J_{t_1}=(y_1,y^+(r(t_1),t_1))$, and in
{\it Case II} set $J_{t_1}=(y^-(r(t_1),t_1), y_1)$. In the special
case where $(r(t),t)\notin\cal U$ for all $t\in(0,t_1)$ we may
set $J_{t_1}$ either one of the two alternatives. 

For the pairwise disjointness it suffices to show that
$J_{t_1}\cap J_t=\emptyset$ for any  $t<t_1$ such 
that $(r(t),t),(r(t_1),t_1)\in\cal U$. In Case I $r(t)<z(t)$, 
and Lemma \ref{charlm2}(ii) implies
$$y^-(r(t),t)<y^+(r(t),t)\le y^-(z(t),t)=y_1. $$
Consequently $J_{t_1}\cap J_t=\emptyset$
because  $J_t\subseteq (y^-(r(t),t),y^+(r(t),t))$ and
$J_{t_1}$ $=$ $(y_1,y^+(r(t_1),t_1))$. Similarly in Case II. 

\hbox{}

{\bf Step 2.} We show that in 
 a bounded set $[a,b]\times(0,T]$,  
all the $\cal U_\alpha$-points lie on 
 a finite collection $\{r_j(\cdot):1\le j\le M\}$
of forward characteristics.  

Let
\beas
&&\Lambda=\{U\subseteq{\cal U}_\alpha: \mbox{
$U\subseteq[a,b]\times(0,T]$, and for any
two distinct points $(x,t)$, $(x',t')$} \\
&&\mbox{$\in U$, the intervals $(y^-(x,t),y^+(x,t))$ 
and $(y^-(x',t'),y^+(x',t'))$ are disjoint}\}.
\eeas
Note that for each $(x,t)\in \cal U_\alpha$,
the interval $(y^-(x,t),y^+(x,t))$ is a subinterval of
$[y^-(a,T),b]$ of length at least $\alpha$.
Consequently no set $U\in\Lambda$ 
contains more than $(b-y^-(a,T))/\alpha$ points.
Let $M=\max\{|U|:U\in\Lambda\}$ be  the maximal number 
of points in any element of $\Lambda$.  
 Fix $U_0\in\Lambda$ that has
$M$ points, 
$$U_0=\{(x_j,t_j): 1\le j\le M\}.$$
For each $j$, define the forward characteristic emanating
from $(x_j,t_j)$ by $r_j(t_j)=x_j$
and  $r_j(t)=w(x_j;t_j,t)$ for $t\in(t_j,\infty)$. 
And let $A_j$ be the open triangle with vertices 
$(x_j,t_j)$, $(y^-(x_j,t_j), 0)$, and $(y^+(x_j,t_j), 0)$.
By Lemma \ref{charlm1}(b) the sides
of this triangle are formed by the line segments
$\{ (y^\pm(x_j;t,t_j),t): 0<t<t_j\}$. Hence
$$A_j=\{(x,t): 0<t<t_j, y^-(x_j;t,t_j)<x<y^+(x_j;t,t_j)\}.$$

As an intermediate conclusion we claim that each 
$\cal U_\alpha$-point in  $[a,b]\times(0,T]$ 
lies either on the forward characteristic from 
some $(x_j,t_j)$, or in some $A_j$. To justify this, note
first that no $\cal U_\alpha$-point can lie on a
side $\{ (y^\pm(x_j;t,t_j),t): 0<t<t_j\}$ of a triangle $A_j$ 
because by Lemma \ref{charlm2}(ii) such a point has 
a unique minimizer. Secondly, if a $\cal U_\alpha$-point
$(x,t)$ lies outside the closure of  
the union of $\{r_j(\cdot), A_j\}$, then by Lemma
\ref{charlm3}(i)--(ii) the open interval 
$(y^-(x,t),y^+(x,t))$ is disjoint from all 
$(y^-(x_j,t_j),y^+(x_j,t_j))$. Then  we can add $(x,t)$
to $U_0$, thereby contradicting the definition of 
$M=|U_0|$ as the maximal size of an element of $\Lambda$. 

To complete Step 2 it remains to argue that we can
extend the definition of $r_j(\cdot)$ to $[0,t_j)$
so that all $\cal U_\alpha$-points in $A_j$ lie on
 $r_j(\cdot)$. 

\hbox{}

{\it Final claim.} There cannot exist two 
$\cal U_\alpha$-points $(x,t)$ and $(x',t')$ in $A_j$ 
such that $t'>t$ but $x'\ne w(x;t,t')$. 

Suppose such points did exist. But then  
$(y^-(x,t),y^+(x,t))$ and $(y^-(x',t'),y^+(x',t'))$
are disjoint. Furthermore, as subintervals of
$(y^-(x_j,t_j),y^+(x_j,t_j))$ they are both disjoint
from all the other intervals $(y^-(x_i,t_i),y^+(x_i,t_i))$, $i\ne j$.
We can contradict the maximality
of $M=|U_0|$ by replacing $(x_j,t_j)$ 
with $(x,t)$ and $(x',t')$. This proves the final claim. 

By the final claim, we can apply Lemma \ref{rlm1}(ii) to
get  a single characteristic $r_j(\cdot)$ that contains 
 all $\cal U_\alpha$-points in 
$A_j$ and those on the forward characteristic from $(x_j,t_j)$.
This completes the proof of Step 2 and thereby the proof
of Theorem \ref{Uthm}. 
\qed

Another technical result we need is that all the shocks in a 
compact set can be enclosed in an open set with small $t$-sections.

\begin{prop} Fix $-\infty<a<b<\infty$, $0<\tau<\infty$, and 
$\e>0$. Then there exists an open set $G\subseteq\RR\times(0,\infty)$
such that $G$ contains all the shocks in $[a,b]\times(0,\tau]$,
and for each $t\in(0,\tau]$, the $t$-section 
$G_t=\{x: (x,t)\in G\}$ has 1-dimensional Lebesgue measure
$|G_t|\le \e$.  
\label{Gprop}
\end{prop}

The following notion will be helpful for the proof: Say a shock 
$(x,t)$ is a {\it new shock} if there does not exist a shock
$(x_0,t_0)$ such that $t_0<t$ and $x=w(x_0;t_0,t)$. There can be 
at most countably many new shocks because if 
$(x,t)$ and $(x',t')$ are new shocks, the open intervals 
$(y^-(x,t),y^+(x,t))$ and $(y^-(x',t'),y^+(x',t'))$  must be 
disjoint by Lemma \ref{charlm3}(ii). 

\hbox{}

{\it Proof of Proposition \ref{Gprop}.} Fix a point $c<y^-(a,\tau)$. 
Let $\{(x_i,t_i):i\ge 1\}$ be the (at most countably many)
 new shocks in $[a,b]\times(0,\tau]$. 
Let $S$ be the set of all shocks in $[c,b+1]\times(0,\tau+1]$.
Define $\Delta y(x,t)=y^+(x,t)-y^-(x,t)$, so that
$\Delta y(x,t)>0$ iff $(x,t)$ is
a shock. Let $\e_1<(\e/2)\cdot(b-y^-(c,\tau))^{-1}$. 
Write $B({\bf x},\delta)$ for the Euclidean ball in $\RR^2$ centered at
${\bf x}$ with  radius $\delta$. 
 Define the following subset of $\RR^2$: 
$$
H=\left\{\bigcup_{(x,t)\in S}(x-\e_1\Delta y(x,t), x+\e_1\Delta y(x,t))\times\{t\}\right\}
\cup
\left\{ \bigcup_{i\ge 1} B((x_i,t_i),2^{-i-2}\e)\right\}.
$$ 

We claim that every shock in $[a,b]\times (0,\tau]$ is an interior point
of $H$. This is clear for new shocks. If $(x_1,t_1)$ is a non-new
shock, we can find a shock 
$(x_0,t_0)$ such that $0<t_0<t$ and $x_1=w(x_0;t_0,t_1)$. Since we chose
$c<y^-(a,\tau)\le y^-(x_1,t_1)$, the shock $(x_0,t_0)$ and the forward
characteristic $w(x_0;t_0,t)$ for $t_0\le t\le t_1$ lie in $S$. 
Furthermore, since $S$ contains the shocks in $[x_1,b+1]\times[t_1,\tau+1]$,
we can choose $t_2>t_1$ so that $S$ contains the forward characteristic
$r(t)\equiv w(x_0;t_0,t)$ for $t_0\le t\le t_2$. Let
$h=\e_1\Delta y(x_0,t_0)>0$. By (\ref{yminusineq1}), 
$\e_1\Delta y(r(t),t)\ge h$ for  $t_0\le t\le t_2$. Consequently $H$ contains
the set $\bigcup_{t_0<t<t_2}(r(t)-h,r(t)+h)\times\{t\}$. 
This latter contains an 
open neighborhood of $(x_1,t_1)$ because $r(\cdot)$ is  a Lipschitz curve
by Lemma \ref{rlm1}(i). 

Let $G$ be the interior of $H$. Then for $t\in(0,\tau]$
\beas
G_t&\subseteq&
\left\{\bigcup_{x: (x,t)\in S}(x-\e_1\Delta y(x,t), x+\e_1\Delta y(x,t))
\right\}\\
&&\qquad\qquad \cup
\left\{ \bigcup_{i\ge 1} (x_i-2^{-i-2}\e,x_i+2^{-i-2}\e)\right\}, 
\eeas
and consequently 
\beas
|G_t|&\le& \sum_{c\le x\le b} 2\e_1\Delta y(x,t)+\sum_{i\ge 1}2^{-i-1}\e\\
&\le& 2\e_1(y^+(b,t)-y^-(c,t))+\e/2 < \e.
\eeas
The inequalities above follow because, as $x\in[c,b]$ ranges over the shock locations
with time coordinate $t$,  the open intervals
$(y^-(x,t),y^+(x,t))$ are disjoint subintervals of $(y^-(c,t),y^+(b,t))$, which 
itself is a subinterval of $(y^-(c,\tau), b)$. 
\qed

\section{Estimates for increasing sequences }
\label{gammasect}

We have the following bounds on ${\bf L}$ and $\mmGa$. 

\begin{lm} Suppose $a$, $s$ and $h$
are positive real numbers.

{\rm (a)} For $x\ge 2$, define 
$$I(x)=2x\cosh^{-1}(x/2)-2\sqrt{x^2-4\,}\,.
$$
When $x>0$ is small enough, there is a constant $C$ such that
$I(2+x)\ge Cx^{3/2}$. 
For any $C$, $I(x)\ge Cx$ for large enough $x$.
For all real $b>0$ and $m\ge 2b$, 
\be
P\{ {\bf L}(b,b)\ge m\} \le \exp\left(-b I(m/b)\right).
\label{Luppertail}
\ee

{\rm (b)}  There are fixed positive constants
$B_0$, $B_1$, $d_0$, $C_0$ and $C_1$ such that if
$a\ge B_0$ and 
$B_1 a^{4/3}\le hs\le d_0a^2$, then 
$$P\left\{ \mmGa([a],s)> \frac{a^2}{4s}+h\right\} \le 
C_0\exp\left\{ -C_1 \frac{s^3h^3}{a^{4}}\right\}\,.
$$

{\rm (c)}   There are finite positive constants
 $C_0$ and $C_1$ such that for all  $0<a\le s$,
$$P\left\{ \mmGa([a],s)> s\right\} \le 
C_0\exp( -C_1s^2)\,.
$$
\label{gammalm}
\end{lm}

Part (a) was first proved by Kim \cite{Kim}. Sepp\"al\"ainen
\cite{Seldp} proved that 
$I(x)$ is the correct rate function for the  deviations
in (\ref{Luppertail}). Part (b) is a consequence
of Lemma 7.1(iv) in Baik-Deift-Johansson \cite{BDJ}. [See Lemma 5.2 in
\cite{Seperturb} for the conversion of Baik-Deift-Johansson's
lemma into part (b) above.]
Part (c) is a consequence of Lemma 2.2 in Johansson \cite{Joh}. 

Next we use these inequalities to derive estimates tailored
to our needs. Most technical complications arise 
from the need to treat small $t$ that vanish as $n\to\infty$,
in order to get the $t$-uniformity of the theorems. 
When using Lemma \ref{gammalm}, it is often useful 
to note that ${\bf L}(a,b)\stackrel{d}={\bf L}\left(\sqrt{ab},\sqrt{ab\,}\,\right)$
($\stackrel{d}=$ means equality in distribution). This follows
from the invariance of the homogeneous planar Poisson point process
under the maps $(x,y)\mapsto(rx,r^{-1}y)$, $r>0$.

\begin{lm} {\rm (a)} Let $\beta, \tau>0$. Then there exists a 
constant $\alpha\in(0,\infty)$ such that 
$$
\sum_{n\ge 1} P\left\{ 
\mmGa([\alpha nt^{1/2}],nt) \le n\beta \ \mbox{ for some 
$t\in[n^{-2}(\log n)^2,\tau]$} \right\} <\infty.
$$

{\rm (b)} Let $b, \delta, \gamma, \tau>0$. Assume these
restrictions: $\gamma<3/4$ and $\gamma(1+\delta)<1$. 
Then there exists a 
constant $\alpha\in(0,\infty)$ such that 
$$
\sum_{n\ge 1} nP\left\{ 
\mmGa([\alpha nt^{\gamma}],nt) \le b\alpha nt^{1/2} 
 \ \mbox{ for some 
$t\in[n^{-(1+\delta)},\tau]$} \right\} <\infty. 
$$

{\rm (c)} Let $b, \delta, \gamma, \eta, \tau>0$. Assume these
restrictions: $1/2< \gamma<3/4$, $\gamma(1+\delta)<1$, and 
$\delta< (3-4\gamma-2\eta)/(4\gamma-2)$.
Then there exists a 
constant $\alpha\in(0,\infty)$ such that 
$$
\sum_{n\ge 1}  n P\left\{ 
\mmGa([\alpha nt^{\gamma}],nt) \le bn^{1/2+\eta}
 \ \mbox{ for some 
$t\in[n^{-(1+\delta)},\tau]$} \right\} <\infty. 
$$
\label{smalltlm1}
\end{lm}

{\it Proof.} Part (a). Set $t_i=4^in^{-2}(\log n)^2$ for $i\ge 0$. 
Pick $K$ so that $t_{K-1}<\tau\le t_K$. Then $K\le C\log n$
for a constant $C$. Let $\alpha_0=\alpha/2$ to account
for the effect of replacing the integer 
$[\alpha nt^{1/2}]$ by $\alpha nt^{1/2}$. 
Suppose $\alpha_0> 6\beta^{1/2}$. 
Pick $a_1$ so that $I(x)\ge a_1x$ for $x\ge 3$.
Increase $\alpha_0$ further so that $a_1\alpha_0\ge 2$. 
 Use Lemma \ref{gammalm}(a) to bound the probability:
\beas
&&P\left\{ 
\mmGa([\alpha nt^{1/2}],nt) \le n\beta \ \mbox{ for some 
$t\in[n^{-2}(\log n)^2,\tau]$} \right\}\\
&\le& 
P\left\{ 
{\bf L}(n\beta,nt) \ge \alpha_0 nt^{1/2} \ \mbox{ for some 
$t\in[n^{-2}(\log n)^2,\tau]$} \right\}\\
&\le& 
\sum_{i=0}^{K-1} P\left\{ 
{\bf L}\left(n\beta,nt_{i+1}\right) \ge \alpha_0 nt_i^{1/2} \right\}\\
&=& 
\sum_{i=0}^{K-1} P\left\{ 
{\bf L}\left(n(\beta t_{i+1})^{1/2},n(\beta t_{i+1})^{1/2}\right)
 \ge \alpha_0 nt_i^{1/2} \right\}\\
&\le& 
\sum_{i=0}^{K-1} \exp\left\{ -n(\beta t_{i+1})^{1/2}
I\left( \alpha_0(4\beta)^{-1/2} \right) \right\}\\
&\le& 
K \exp\left\{ -a_1\alpha_0\log n\right\}
\le Cn^{-2}\log n. 
\eeas
This bound is summable over $n$. 

Part (b).  Follow a similar  partition argument, with 
$t_i=4^in^{-(1+\delta)}$. 
Take 
$\alpha$ sufficiently large, $\alpha_0=\alpha/2$, and 
use Lemma \ref{gammalm}(a) to get 
the upper bound
\beas
&&P\left\{ 
\mmGa([\alpha nt^{\gamma}],nt) \le b\alpha nt^{1/2} 
 \ \mbox{ for some 
$t\in[n^{-(1+\delta)},\tau]$} \right\}\\
&\le&\sum_{i=0}^{K-1} P\left\{ 
{\bf L}\left(b\alpha nt_{i+1}^{1/2},nt_{i+1}\right) 
\ge \alpha_0 nt_i^{\gamma} \right\}\\
&\le& 
\sum_{i=0}^{K-1} \exp\left\{ -n(b\alpha)^{1/2} t_{i+1}^{3/4}
\,I\left( (1/2)(\alpha/b)^{1/2}4^{-\gamma}t_{i+1}^{\gamma-3/4} \right) \right\}\\
&\le& 
\sum_{i=0}^{K-1} \exp\left\{ -(1/2)na_1\alpha t_i^\gamma\right\}
\le C\log n\exp(-Cn^{1-\gamma(1+\delta)}).  
\eeas
To get the inequality in part (b), 
multiply this bound by $n$
and use the assumption 
$\gamma(1+\delta)<1$. 

Part (c). With the same  partition as in part (b), 
\beas
&&\sum_{i=0}^{K-1} P\left\{ 
{\bf L}\left(b n^{1/2+\eta},nt_{i+1}\right) 
\ge \alpha_0 nt_i^{\gamma} \right\}\\
&\le& 
\sum_{i=0}^{K-1} \exp\left\{ -n^{3/4+\eta/2}(bt_{i+1})^{1/2} 
I\left( (1/2)\alpha_0 b^{-1/2}
n^{1/4-\eta/2} t_{i}^{\gamma-1/2} 
\right) \right\}.
\eeas
Replace $t_i$ by its lower bound $n^{-(1+\delta)}$ inside $I$
to get $n^{1/4-\eta/2} t_{i}^{\gamma-1/2} 
\ge n^{1/4-\eta/2-(1+\delta)(\gamma-1/2)}$. This last
exponent is positive by the assumption on $\delta,\eta,\gamma$.
Now proceed as above. 
\qed

\begin{lm} Let $r,\tau>0$ be positive constants. Then 
there exist finite positive constants $C_0$ 
and $C_1$ such that, for all $n\ge 1$, 
\beas
&&\sum_{m=1}^{[nr]}
P\left\{ 
\mmGa(m,nt)\le \frac{m^2}{4nt}-\frac{m^{4/3}\log n}{4nt}
\ \mbox{ for some $t\in[n^{-2},\tau]$} 
\right\}\\
&&\qquad\le
C_0 n^{5/3}\exp\left(-C_1(\log n)^{3/2}\right).
\eeas
\label{largetlm1}
\end{lm}

{\it Proof.} It suffices to consider $m^{2/3}>\log n$, otherwise the 
probability is 0. 
As in the previous proof, 
 partition the time interval $[n^{-2},\tau]$
by
$n^{-2}=t_0<t_1<\cdots<t_{K-1}<\tau\le t_K$, and bound 
the probability by Lemma \ref{gammalm}(a):
\beas
&&P\left\{ 
\mmGa(m,nt)\le \frac{m^2}{4nt}-\frac{m^{4/3}\log n}{4nt}
\ \mbox{ for some $t\in[n^{-2},\tau]$} 
\right\} \\
&\le& 
P\left\{ 
{\bf L}\left( m^2(4nt)^{-1}(1-m^{-2/3}\log n), nt\right) \ge m 
\ \mbox{ for some $t\in[n^{-2},\tau]$} 
\right\}\\
&\le& 
\sum_{i=0}^{K-1} P\left\{ 
{\bf L}\left( m^2(4nt_i)^{-1}(1-m^{-2/3}\log n), nt_{i+1}\right) 
\ge m    \right\}\\
&\le&  
\sum_{i=0}^{K-1} \exp\left(-b_i I(m/b_i)\right),
\eeas
where we wrote 
$$b_i=\frac{m}2 \left(\frac{t_{i+1}}{t_i}\right)^{1/2}
\left( 1-m^{-2/3}\log n\right)^{1/2}.
$$
We argue separately for two ranges of $m$. 

\hbox{}

{\it Case 1:} $\log n<m^{2/3}\le (1+\delta)\log n$ for 
 a small $\delta\in(0,1/4)$. Define 
the partition by $t_i=4^it_0$. Check that then $m/b_i\ge \delta^{-1/2}$.
Pick a constant $a_0$ so that $I(x)\ge a_0x$ for $x\ge \delta^{-1/2}$.
(This makes sense because  
$\delta^{-1/2}>2$.) The size $K$ of the partition satisfies
$K\le C\log n$ for a constant $C$ that depends on $\tau$. 
Summing the bound $\exp\left(-b_i I(m/b_i)\right)\le \exp(-a_0m)$
 over $m$ and $i$ gives the following upper
bound:
$$
\sum_{m: \log n<m^{2/3}\le (1+\delta)\log n}\;
\sum_{i=0}^{K-1} e^{-a_0 m} \le C_0 (\log n)^{5/2} 
\exp\left(-C_1 (\log n)^{3/2}\right)
$$
for suitable constants $C_0,C_1$. 

\hbox{}

{\it Case 2:} $ (1+\delta)^{3/2}(\log n)^{3/2}\le m\le nr$. 
Set 
$\theta=1+(1/2)m^{-2/3}\log n$, 
 and define the partition
by  $t_i=t_0\theta^i$.
Now $K\le Cm^{2/3}$. Note that
\beas
\frac{m}{b_i}&=&2 \left(\frac{t_{i}}{t_{i+1}}\right)^{1/2}
\left( 1-\frac{\log n}{m^{2/3}}\right)^{-1/2}
\ge 2 \left(\frac{t_{i}}{t_{i+1}}\right)^{1/2}
\left( 1+\frac{\log n}{2m^{2/3}}\right)\\
&=& 2\left( 1+\frac{\log n}{2m^{2/3}}\right)^{1/2}
\ge 2+\frac{\log n}{4m^{2/3}} ,
\eeas
where we used the inequalities $(1-h)^{-1/2}\ge 1+h/2$
and $(1+h)^{1/2}\ge 1+h/4$ that are valid for $0\le h<1$.
Choose $a_1$ so that $I(2+x)\ge a_1x^{3/2}$ for 
$0<x<1/4$. Check that $b_i\ge Cm$ for a constant $C=C(\delta)$. 
Then finally the bound becomes
$$\sum_{i=0}^{K-1}\exp\left(-b_iI(m/b_i)\right)\le 
K\exp\left\{-Cm\left(\frac{\log n}{4m^{2/3}}\right)^{3/2}\right\}
\le Cm^{2/3}\exp\left(-C_1(\log n)^{3/2}\right). 
$$
Summing this over 
$(1+\delta)^{3/2}(\log n)^{3/2}\le m\le nr$, and 
combining with {\it Case 1} above, gives 
the bound in the statement of the lemma.
\qed

\begin{lm} Let $ \tau>0.$  
There exist finite positive constants $M_0$, $n_0$, $C_0$,
and $C_1$ such that, for all $n\ge n_0$ and 
$m\ge M_0(\log n)^{3/2}$, 
\beas
&&P\left\{ 
\mmGa(m,nt)> \frac{m^2}{4nt}+\frac{m^{4/3}\log n}{4nt}
\ \mbox{ for some $t\in[n^{-1},\tau]$} 
\right\}\\
&&\qquad\le
 C_0 m^{2/3}\exp\left(-C_1(\log n)^3\right).
\eeas
\label{largetlm2}
\end{lm}

{\it Proof.}  Fix $n\ge 3,m\ge 1$. Set 
$$\theta=\frac{2m^{2/3}+\log n}{2m^{2/3}} >1.$$
Consider a partition of $[ n^{-1},\tau]$, defined by 
$t_0= n^{-1}$, $t_i=t_0\theta^i$ for $i\ge 1$, and 
choose $K$ so that  
$t_{K-1}<\tau\le t_K$. Then 
$K\le Cm^{2/3}$ for a constant $C=C(\tau)$. 

 Bound the probability in question from above by 
$$
\sum_{i=0}^{K-1} 
P\left\{ 
\mmGa(m,nt_i)> \frac{m^2}{4nt_{i+1}}+\frac{m^{4/3}\log n}{4nt_{i+1}}
\right\}.
$$
To apply Lemma \ref{gammalm}(b) to each of these probabilities, identify
$a=m$, $s=nt_i$, and $h=m^{4/3}\log n (4nt_{i+1})^{-1}
-m^2(4n)^{-1}(t_i^{-1}-t_{i+1}^{-1})$. 
 Check that
there exist $n_0$ and 
 $M_0$ such that, if $n\ge n_0$ and  $m\ge M_0(\log n)^{3/2}$,
then the hypotheses of Lemma \ref{gammalm}(b) are met. Observe
that $hs= (8\theta)^{-1}m^{4/3}\log n$. Applying the 
estimate in Lemma \ref{gammalm}(b) to each term of 
the sum above gives an upper bound of
$$C_0K\exp\left( -C_1 (sh)^3a^{-4} \right)
\le C_0 m^{2/3}\exp(-C_1(\log n)^3),
$$
where the constants $C_i$ are no longer the original ones
from Lemma \ref{gammalm}(b). 
\qed

The range of $m$'s not covered by the last lemma are 
taken care of by the next statement. 

\begin{lm} Let $\tau, \alpha, M$ be positive constants. Then 
there exist finite positive constants $n_0$, $C_0$
and $C_1$ such that, for all $n\ge n_0$, 
\beas
&&P\left\{ 
\mmGa(m,nt)> n^\alpha
\ \mbox{ for some $t\in[n^{-1},\tau]$ and 
$1\le m\le M(\log n)^{3/2 }$} 
\right\}\\
&&\qquad\le
C_0\exp(-C_1n^\alpha).
\eeas
\label{largetlm3}
\end{lm}

{\it Proof.} Since $\mmGa(m,nt)$ is nonincreasing in $t$
and nondecreasing in $m$, the probability is bounded by 
\beas
&&P\left\{ 
\mmGa([M(\log n)^{2/3}],1)> n^\alpha\right\}\\
&=&P\left\{ 
\mmGa([M(\log n)^{2/3}],n^{\alpha/2})>
n^{\alpha/2}\right\}
\le C_0\exp(-C_1n^\alpha).
\eeas
We used the equality in distribution 
${\bf L}(a,b)\stackrel{d}={\bf L}((ab)^{1/2},(ab)^{1/2})$
and Lemma \ref{gammalm}(c).  
\qed

\section{Estimates for the microscopic variational formula}
\label{varcoupsect}

Recall that for the $n$th process the variational coupling equality,
appropriately scaled,  
reads
\be
z_{[nx]}^n(nt)=\inf_{i:i\le [nx]}
\{ z^n_i(0) +\Gamma^{n,i}_{[nx]-i}(nt)\}.
\label{varcoupn}
\ee
Let $i_n(x,t)$ denote the minimal $i$ at which 
the infimum is attained in (\ref{varcoupn}):
\be
i_n(x,t)=\inf\{i: z_{[nx]}^n(nt)=
 z^n_i(0) +\Gamma^{n,i}_{[nx]-i}(nt) \}.
\label{inxtdef}
\ee
 It is proved in \cite{Seejp} 
that under assumption (\ref{unifass}) $i_n(x,t)$ is almost surely 
finite, and it 
 is a nonincreasing function of $t$. For $t=0$ we interpret
$\Gamma^{n,i}_m(0)=\infty$ for $m\ge 1$, and $i_n(x,0)=[nx]$. 

The technical key to benefiting from (\ref{varcoupn}) lies
in estimating how far the $n^{-1}$-scaled random minimizing indices 
of (\ref{varcoupn}) lie from the set  $I(x,t)$ of macroscopic minimizers. 
We start with a crude bound, and successively
refine it.

\begin{lm} 
For $c<a<b$ and  $\tau>0$ define the events
\be
G_n=\left\{ 
\mbox{for some $x\in[a,b]$ and $t\in(0,\tau]$, 
{\rm (\ref{varcoupn})} is minimized by some $i\le cn$} \right\}. 
\label{Gndef}
\ee
For fixed $a<b$ and $\tau$ we can choose  $c<0$ such that this holds: 

{\rm (i)} Under the uniformity assumption {\rm (\ref{unifass})}, 
$\lim_{n\to\infty} P(G_n)=0$. 

{\rm (ii)} Under assumptions {\rm (\ref{assrho}) and (\ref{assloqeq})},
$\sum_{n=1}^\infty P(G_n)<\infty.
$
\label{Gnlm}
\end{lm}

{\it Proof.}  Since $i_n(x,t)$ is nonincreasing in $t$, we can 
express (\ref{Gndef}) as 
$$G_n=\left\{ 
\mbox{for some $x\in[a,b]$ and $t=\tau$, 
{\rm (\ref{varcoupn})} is minimized by some $i\le cn$} \right\}. 
$$
 Let $C_1>0$, and pick $C_0>0$ so that 
$I(x)\ge C_1x$ for $x\ge C_0$ [recall Lemma \ref{gammalm}(a)]. Pick $\e>0$ small
enough so that $C_0^2\tau\e<1$. And then pick $c<0$ so that 
$c<a\left(1-C_0\sqrt{\tau\e\,}\,\right)^{-1}$ and $c<q$ for the $q$ that
satisfies assumption (\ref{unifass}) for $\e$. 
Abbreviate $Y_{n,i}=z^n_{[nb]}(0)-z^n_i(0)$. 
Under these conditions 
$i\le nc$ implies $\e i^2\le ([na]-i)^2/(C_0^2\tau)$, and 
consequently on the event 
$$A_n=\{ \mbox{$i^{-2}nY_{n,i}\le \e$ for all $i\le nc$}\}$$
we have 
\be
\left([na]-i\right)\cdot\left(n\tau Y_{n,i}\right)^{-1/2}\ge C_0
\quad\mbox{ for all $i\le nc$.}
\label{Yineq1}
\ee
Now we can estimate:
\beas
P(G_n)&\le& 
P\left\{ 
\mbox{for some $i\le cn$,
${\bf L}\left((z^n_i(0),0),(z^n_{[nb]}(0), n\tau)\right) \ge [na]-i$
}\right\}\\
&\le&
P(A_n^c)+\sum_{i\le nc} E\left[ {\bf 1}_{A_n}\cdot 
\exp\left\{ -(n\tau Y_{n,i})^{1/2} I\left( ([na]-i)(n\tau 
Y_{n,i})^{-1/2}\right)\right\}
\right]\\
&\le&
P(A_n^c)+\sum_{i\le nc} \exp[-C_1([na]-i) ].
\eeas
The first inequality above comes from 
$z^n_i(0) +\Gamma^{n,i}_{[nx]-i}(n\tau)\le z^n_{[nx]}(0)$ 
which must follow if $i$ is to be a minimizer in (\ref{varcoupn}),
and also from $[na]\le[nx]\le[nb]$. The second inequality
comes from (\ref{Luppertail}), and the last
 from  (\ref{Yineq1}) and $I(x)\ge C_1x$.

Assumption (\ref{unifass}) implies that $P(A_n^c)\to 0$, so 
part (i) of the lemma is 
proved. To get the summability $\sum P(G_n)<\infty$
required in part (ii), we need to check that
 assumptions (\ref{assrho}) and (\ref{assloqeq}) 
imply $\sum P(A_n^c)<\infty$. 

Let $\rho^*(r)=\sup_{r\le x\le b}\rho_0(x)$. Write 
$$P(A_n^c)=P\left(\,\sum_{i=j+1}^{[nb]}\eta^n_i(0)
> \frac{\e j^2}n\  \mbox{for some $j\le nc$}\right).
$$
Since the variables
 $\{\eta^n_i(0): j<i\le [nb]\}$ are
independent exponentials with means bounded by $\rho^*(j/n)$,
they are stochastically
dominated by  $\{ \rho^*(j/n) X_i: j<i\le [nb]\}$ where
the $X_i$'s are i.i.d.\ exponential variables with common 
mean $EX_i=1$. 
Recall that for $s>1$ we have the 
large deviation bound $P(\sum_1^m X_i\ge ms)\le \exp(-m\kappa(s))$
with the rate function $\kappa(s)=s-1-\log s$. Thus
\beas
P(A_n^c)\le \sum_{j\le nc} P\left(\,\sum_{i=j+1}^{[nb]}X_i 
> \frac{\e j^2}{n\rho^*(j/n)}\right)
\le \sum_{j\le nc} \exp\left\{-([nb]-j)
\kappa(s_{n,j})\right\},
\eeas
where 
$$s_{n,j}=
 \frac{\e j^2}{n([nb]-j)\rho^*(j/n)}.
$$
By assumption (\ref{assrho}) we can guarantee $s_{n,j}\ge M$ for an
arbitrarily large  $M$, for all $n$ and $j$, by taking $c<0$ large enough
negative.  Then $\sum P(A_n^c)<\infty$ follows, and the 
lemma is proved. 
\qed

\begin{lm}  Let $a<b$, $\alpha>0$, $\tau>0$, $\gamma\in(1/2,3/4)$.
Suppose 
$\delta>0$ satisfies $\gamma(1+\delta)<1$ and 
$\delta< (3-4\gamma)/(4\gamma-2)$.  
Define the events
\beas
H_{n,0}&=&\left\{\mbox{for some $x\in[a,b]$ and
 $t\in[n^{-(1+\delta)},\tau]$,}\right.\\ 
&&\qquad\left.\mbox{{\rm (\ref{varcoupn})} is minimized  
by some $i< [nx]-\alpha nt^{\gamma}$}\right\}
\eeas
and 
\beas
H_{n,1}&=&\left\{\mbox{for some $x\in[a,b]$ and
 $t\in(0,n^{-(1+\delta)}]$,}\right.\\ 
&&\qquad\left.\mbox{{\rm (\ref{varcoupn})} is minimized  
by some $i< [nx]-\alpha n^{(1-\delta)/2}$}\right\}.
\eeas
If  $\alpha$ is chosen large enough, the following
is true: 

{\rm (i)} Under assumptions  {\rm (\ref{ass1d})}  and  
 {\rm (\ref{unifass})}, 
$\lim_{n\to\infty} P(H_{n,0}\cup H_{n,1})=0$. 

{\rm (ii)} Under assumptions {\rm (\ref{assrho})} and  
 {\rm (\ref{assloqeq})},
$\sum_{n=1}^\infty P(H_{n,0}\cup H_{n,1})<\infty.
$
\label{Hn01lm}
\end{lm}

{\it Proof.} To see that $H_{n,0}$ is a measurable event, 
let $T$ be a countable dense subset of $[n^{-(1+\delta)},\tau]$
that contains $\tau$. Then almost surely 
$$H_{n,0}=\bigcup_{k=[na]}^{[nb]}\bigcup_{t\in T}
\bigcup_{i<k-\alpha nt^\gamma}\left\{ z^n_k(nt)=z^n_i(0)+
\Gamma^{n,i}_{k-i}(nt)\right\}.
$$
To see why it suffices to consider only $t\in T$ in the union
above, note that as functions of $t$,
$z^n_k(nt)$ and $\Gamma^{n,i}_{k-i}(nt)$ are
right-continuous jump processes whose jumps do not accumulate
with probability 1. So almost surely, for any $t$ there exists
a (random) $\e>0$ such that these processes do not jump in $(t,t+\e)$.
A similar argument works for $H_{n,1}$ also. 

 We prove statements (i) and (ii)
first for $H_{n,0}$. The challenge here is in the 
small values of $t$ that vanish as $n\to\infty$. 
 The proof will be achieved in two rounds.
First we rule out minimizers $i\le nx-\alpha nt^{1/2}$,
and then in the second step we rule out 
$i\le nx-\alpha nt^{\gamma}$.
 By conditioning on the event $G_n^c$
of Lemma \ref{Gnlm}, it suffices to consider $i\ge nc$.  

Suppose some $i\in [nc, nx-\alpha nt^{1/2}]$
minimizes (\ref{varcoupn}) for some  
$x\in[a,b]$ and 
 $t\in[n^{-(1+\delta)},\tau]$. Since $[nx]$ is among the 
indices over which the infimum is taken in (\ref{varcoupn}), 
it must follow that 
 $$z^n_i(0)+\Gamma^{n,i}_{[nx]-i}(nt)\le z^n_{[nx]}(0)\le z^n_{[nb]}(0).$$
Bound the left-hand side from below:
$$
z^n_i(0)+\Gamma^{n,i}_{[nx]-i}(nt)\ge 
z^n_{[nc]}(0)+\Gamma^{n,[nc]}_{[nx]-i}(nt)\ge  
z^n_{[nc]}(0)+\Gamma^{n,[nc]}_{[\alpha nt^{1/2}]}(nt).
$$
The consequence is that for some $t\in[n^{-(1+\delta)},\tau]$,
$$\Gamma^{n,[nc]}_{[\alpha nt^{1/2}]}(nt)
\le z^n_{[nb]}(0)-z^n_{[nc]}(0).$$
Let $\beta=u_0(b)-u_0(c)+1$, and $G_n$ be the event in 
Lemma \ref{Gnlm}. Then the previous reasoning gives
\beas
&&P\left\{\mbox{for some $x\in[a,b]$ and
 $t\in[n^{-(1+\delta)},\tau]$,}\right.\\ 
&&\left.\qquad\qquad\qquad\qquad\mbox{{\rm (\ref{varcoupn})} is minimized  
by some $i\le nx-\alpha nt^{1/2}$}\right\}\\
&\le& P(G_n)+ P\left\{ 
\mmGa\left([\alpha nt^{1/2}], nt\right) \le n\beta
\ \mbox{for some $t\in[n^{-(1+\delta)},\tau]$} \right\}
\\
&&\qquad +P\left\{ z^n_{[nb]}(0)-z^n_{[nc]}(0)> n\beta \right\}.
\eeas
Note that $\gamma(1+\delta)<1$ forces $\delta<1$,
and then $n^{-(1+\delta)}>n^{-2}(\log n)^2$ for large $n$. 
Thus Lemma \ref{smalltlm1}(a) applies, and we can conclude
that the second probability after the inequality above is 
summable over $n\ge 1$ if $\alpha$ is chosen large enough.    The last probability
 converges to zero in Case (i), 
and is summable over $n$ in Case (ii) of the lemma. 

Now condition on the event that all minimizers 
satisfy $i\ge nx-\alpha nt^{1/2}$, for $(x,t)$ in the 
range under consideration. Under this condition, 
\beas
H_{n,0}&\Longrightarrow& 
\mbox{for some $x\in[a,b]$, 
 $t\in[n^{-(1+\delta)},\tau]$,   
$i\in [nx-\alpha nt^{1/2},nx-\alpha nt^{\gamma}]$:}\\
&&\qquad\qquad\qquad\qquad
 z^n_i(0)+\Gamma^{n,i}_{[nx]-i}(nt)\le z^n_{[nx]}(0)\\
&\Longrightarrow& \mbox{for some  $x\in[a,b]$, 
$t\in[n^{-(1+\delta)},\tau]$,}\\
&&\qquad\qquad\qquad\qquad
\Gamma^{n,[nx-\alpha nt^{1/2}]}_{[\alpha nt^{\gamma}]}(nt)
\le z^n_{[nx]}(0)-z^n_{[nx-\alpha nt^{1/2}]}(0).
\eeas
By the definition of $\zeta_n$, we can write
\beas
&&z^n_{[nx]}(0)-z^n_{[nx-\alpha nt^{1/2}]}(0)\\
&=&n\left(u_0(x)-u_0(x-\alpha t^{1/2})\right)+n^{1/2}
\left(\zeta_n(x,0)-\zeta_n(x-\alpha t^{1/2},0)\right)\\
&\le& 
C\alpha n t^{1/2}+2n^{1/2}\cdot\sup_{x\in [c,b]}\zeta_n(x,0),
\eeas
where $[c,b]$ is an interval that contains 
$[x-\alpha t^{1/2},x]$ for all $(x,t)$ under consideration, 
and $C$ is the Lipschitz constant for 
$u_0$ on the interval $[c,b]$.
By assumption $\delta< (3-4\gamma)/(4\gamma-2)$, so we can 
pick a small $\eta>0$ so that $\delta< (3-4\gamma-2\eta)/(4\gamma-2)$. 
Now summarize everything in this upper bound:
\beas
&&P(H_{n,0})\le P(G_n) + P\left\{ 
\mmGa\left([\alpha nt^{1/2}], nt\right) \le n\beta
\ \mbox{for some $t\in[n^{-(1+\delta)},\tau]$} \right\}
\\
&&\qquad +P\left\{ z^n_{[nb]}(0)-z^n_{[nc]}(0)> n\beta \right\}\\
&&+P\left\{ \Gamma^{n,[nx-\alpha nt^{1/2}]}_{[\alpha nt^{\gamma}]}(nt)
\le C\alpha n t^{1/2}+2n^{1/2+\eta} \ \mbox{for some  $x\in[a,b]$, 
$t\in[n^{-(1+\delta)},\tau]$} \right\}\\
&&\qquad +P\left\{ \sup_{x\in [c,b]}\zeta_n(x,0) > n^\eta\right\}\\
&&\equiv P(G_n) + p_{n,1}+p_{n,2}+p_{n,3}+p_{n,4}.
\eeas
By Lemma  \ref{smalltlm1}(a)--(c), $\sum_n p_{n,j}<\infty$ for $j=1,3$. 
Note that for $p_{n,3}$ we have to sum over 
the superscript $[nx-\alpha nt^{1/2}]$ as $x$ varies over
$[a,b]$. This gives $O(n)$ terms, which is why the 
probabilities in Lemma \ref{smalltlm1}(b)--(c) are multiplied by $n$. 
Under assumption (\ref{ass1d}), 
$\lim_{n\to\infty} p_{n,k}=0$ for $k=2,4$. Under 
assumption (\ref{assloqeq}) and local boundedness of $\rho_0$, 
 $\sum_n p_{n,k}<\infty$ for $k=2,4$.
This proves the lemma for the event $H_{n,0}$.

Repeat the first step for $H_{n,1}$.
Notice  that $\Gamma^{n,[nc]}_{[\alpha n^{(1-\delta)/2}]}(nt)$
is nonincreasing in $t$, so we can replace $t$ by its upper
bound $n^{-(1+\delta)}$. Then we get 
\beas
P(H_{n,1})&\le& P(G_n) + P\left\{ 
\mmGa\left([\alpha n^{(1-\delta)/2}], n^{-\delta}\right) \le n\beta
 \right\}
\\
&&\qquad +P\left\{ z^n_{[nb]}(0)-z^n_{[nc]}(0)> n\beta \right\}.
\eeas
These probabilities are handled as above. Note that the next
to last probability is the special case $t=n^{-(1+\delta)}$
of the event in Lemma \ref{smalltlm1}(a). 
\qed

As usual, the distance between a  point $x$ and a set $A$ is 
denoted by $\dist(x,A)=\inf\{|x-y|: y\in A\}$.

\begin{lm} Let $A\subseteq\RR\times[0,\infty)$ be a compact set. Assume
that $A$ satisfies either  assumption (a) or (b):

(a) $A$ is a finite set; or 

(b) there are no shocks in $A$, in other words $y^-(x,t)=y^+(x,t)$
for all $(x,t)\in A$.
  
 For   $\delta>0$, define the events $H_n=H_n(\delta)$ by
\beas
H_n&=&\left\{\mbox{for some $(x,t)\in A$, {\rm (\ref{varcoupn})} is minimized  
by some $i$}\right.\\
&&\left.\qquad\mbox{ such that $\dist\left(n^{-1}i, I(x,t)\right)>
\delta$}\right\}.
\eeas

{\rm (i)} Under assumptions   
 {\rm (\ref{ass1d})} and {\rm (\ref{unifass})}, 
$\lim_{n\to\infty} P(H_n)=0$. 

{\rm (ii)} Under assumptions {\rm (\ref{assrho})} and  
 {\rm (\ref{assloqeq})},
$\sum_{n=1}^\infty P(H_n)<\infty.
$
\label{Hnlm}
\end{lm}

{\it Proof.}  Measurability of $H_n$ is obvious for a finite $A$.
We prove the measurability of $H_n$ for the other case in 
the appendix.  

Fix finite $a<b$ and $\tau>0$ so that 
$A\subseteq [a,b]\times[0,\tau]$. For small enough $\sigma>0$, 
$I(x,t)\subseteq [x-\delta/2,x]$ for all $x\in[a,b]$ and 
$t\in(0,\sigma]$ by Lemma \ref{Ixtlm}, 
so Lemma \ref{Hn01lm} gives the conclusion for 
$0< t\le \sigma$. Thus for the proof we can assume that 
$A\subseteq [a,b]\times[\sigma,\tau]$
where $0<\sigma<\tau$. The important point here is bounding 
$t$ away from 0 because the estimation gets harder if $t\to 0$ as
$n\to\infty$.  

Choose $c<0$, $c<y^-(a,\tau)\le a$
 so that Lemma \ref{Gnlm} is satisfied. 
By that  Lemma  we only need to consider minimizers 
in the range $[nc, nb]$.
Let 
$$I(x,t)^{(\delta)}=\{q: \mbox{$|q-y|<\delta$ for some $y\in I(x,t)$}\}
$$
be the $\delta$-neighborhood of $I(x,t)$. Set 
\beas
\e&=&\frac15\cdot\inf\{ u_0(y)+tg((x-y)/t)-u(x,t):\\
 &&\qquad\qquad\qquad (x,t)\in A, y\in [c,x]\setminus
I(x,t)^{(\delta)} \}. 
\eeas

We claim that $\e$ is a positive quantity if $A$ satisfies one of the 
two assumptions (a) or (b) in the statement of the lemma.
This is clear if $A$ is finite. Suppose next that 
$y^\pm(x,t)=y(x,t)$ for all $(x,t)\in A$, but $\e=0$.   
Pick a sequence $(x_j,t_j)$ in $A$ and $y_j\in[c,x_j]\setminus I(x_j,t_j)^{(\delta)}$
 so that 
$$
u_0(y_j)+t_jg((x_j-y_j)/t_j)-u(x_j,t_j)\to 0\qquad\mbox{as $j\to\infty$}.
$$
Pass to a convergent subsequence $(x_j,t_j,y_j)\to (\xbar,\tbar,\ybar)$ with 
$(\xbar,\tbar)\in A$. By 
continuity, 
$$
u_0(\ybar)+\tbar g((\xbar-\ybar)/\tbar)-u(\xbar,\tbar)=0
$$
which implies that $\ybar$ must be the Hopf-Lax minimizer for
$(\xbar,\tbar)$, in other words $\ybar=y(\xbar,\tbar)$. 
On the other hand, we also 
have $y(x_j,t_j)\to y(\xbar,\tbar)$ and $|y_j-y(x_j,t_j)|\ge \delta$,
so in the $j\to\infty$ limit $|\ybar-y(\xbar,\tbar)|\ge\delta$. 
This contradiction shows that $\e>0$. 

Note that we cannot
make this argument in case $y^-(\xbar,\tbar)\ne y^+(\xbar,\tbar)$, because
then it is perfectly possible that $y_j\to \ybar=y^-(\xbar,\tbar)$ 
while
$y(x_j,t_j)\to y^+(\xbar,\tbar)$ without contradicting  
$|\ybar-y^+(\xbar,\tbar)|\ge\delta$. This is the step where the proof of
a uniform limit fails for a compact set with shocks. 
Of course, Remark \ref{unifremark} already showed that we cannot hope to
prove a uniform limit for such a set.

For each 
$x\in[a,b]$, $t\in[\sigma,\tau]$, choose 
 finitely many points $a_{k}=a_{k}(x,t)$
and $b_{k}=b_{k}(x,t)$, $1\le k\le K=K(x,t)$, so that 
$$ [c,x]\setminus
I(x,t)^{(\delta)} =\bigcup_{k=1}^{K} [a_{k},b_{k}]
$$
and 
$$\left| tg\left(\frac{x-b_{k}}t\right)-
tg\left(\frac{x-a_{k}}t\right)\right| \le \e
$$
for each $k=1,\ldots,K$. 
To do this,
 pick $\delta_1\in(0,\delta)$ so that 
$|tg(q/t)-tg(r/t)|<\e$ for all $q,r\in[0,b-c]$, $t\in[\sigma,\tau]$,
such that $|q-r|\le \delta_1$. Then pick a partition 
$c=y_0<y_1<\cdots<y_m=b$ with mesh $\max(y_{i+1}-y_i)<\delta_1$.
Every connected component of $I(x,t)^{(\delta)}$ 
is an open interval of length at least $2\delta$, so
each $\{[y_i,y_{i+1}]\cap[c,x]\}\setminus I(x,t)^{(\delta)}$
is either empty or a closed interval. Let
$\{[a_{k},b_{k}]: 1\le k\le K\}$ be the collection 
of the nonempty ones among the intervals 
$\{([y_i,y_{i+1}]\cap[c,x])\setminus I(x,t)^{(\delta)}:0\le i\le m-1\}$. 

For each $(x,t)\in A$, choose a point 
 $y_{x,t}\in I(x,t)$. Reason as follows:
\beas
&&\mbox{for some $(x,t)\in A$, 
(\ref{varcoupn}) is minimized  
by some $i$}\\
 &&\qquad\qquad
\mbox{such that $n^{-1}i\in [c,x]\setminus I(x,t)^{(\delta)}$}\\
&\Longrightarrow&
\mbox{for some $(x,t)\in A$, 
(\ref{varcoupn}) is minimized  
by some $i$}\\
&&\qquad\qquad\mbox{such that $n^{-1}i\in [a_{k},b_{k}]$ for some 
$1\le k\le K$}\\
&\Longrightarrow&
\mbox{for some $(x,t)\in A$  and $1\le k\le K$,}\\
&&\qquad\qquad 
\mbox{$z^n_{[na_{k}]}(0)+\Gamma^{n,[na_{k}]}_{[nx]-[nb_{k}]}(nt) 
\le 
z^n_{[ny_{x,t}]}(0)+\Gamma^{n,[ny_{x,t}]}_{[nx]-[ny_{x,t}]}(nt)$}
\\
&\Longrightarrow&
\mbox{for some $(x,t)\in A$  and $1\le k\le K$,}\\
&&\qquad\qquad\mbox{either $z^n_{[na_{k}]}(0)<nu_0(a_{k})-n\e$,}\\
&&\qquad\qquad\mbox{or 
$\Gamma^{n,[na_{k}]}_{[nx]-[nb_{k}]}(nt)<ntg((x-b_{k})/t)-n\e$,}\\
&&\qquad\qquad\mbox{or
$z^n_{[ny_{x,t}]}(0)>nu_0(y_{x,t})+n\e$,}\\
&&\qquad\qquad\mbox{or  
$\Gamma^{n,[ny_{x,t}]}_{[nx]-[ny_{x,t}]}(nt)>ntg((x-y_{x,t})/t)+n\e$}\\
&\Longrightarrow&
\mbox{for some $y\in[c,b]$, $|z^n_{[ny]}(0)-nu_0(y)|>n\e$,
         or for  some  }\\
&&\qquad\qquad\mbox{$x\in[a,b]$, $t\in[\sigma,\tau]$,   $y\in[c,x]$,
$\left|\Gamma^{n,[ny]}_{[nx]-[ny]}(nt)-ntg((x-y)/t)\right|\ge n\e$}.
\eeas
The next to last implication above
 followed from the choice of $\e$, because
\beas
u_0(a_{k})+tg\left(\frac{x-b_{k}}t\right)&\ge& 
u_0(a_{k})+tg\left(\frac{x-a_{k}}t\right)-\e\\
&\ge& 
u_0(y_{x,t})+tg\left(\frac{x-y_{x,t}}t\right)   +4\e.
\eeas

The entire argument can be summarized in this bound: 
\beas
&&P(H_n)\le P(G_n) + P(H_{n,0})  + P(H_{n,1}) \\
&&\quad+ P\left(\mbox{ $|z^n_{[ny]}(0)-nu_0(y)|>n\e$
for some $y\in[c,b]$ }\right)\\
&+& P\left( \mbox{ 
$\left|\Gamma^{n,[ny]}_{[nx]-[ny]}(nt)-ntg((x-y)/t)\right|\ge n\e$
for some $x\in[a,b]$,  $y\in[c,x]$, $t\in[\sigma,\tau]$ } \right).
\eeas
Apply the assumptions and previous lemmas to treat the 
terms on the right-hand side above. The probabilities
of $\Gamma^{n,[ny]}_{[nx]-[ny]}(nt)$ are handled by Lemmas
\ref{largetlm1} and \ref{largetlm2}.
Assumption (\ref{ass1d}) of weak convergence of
$n^{-1/2}\{z^n_{[ny]}(0)-nu_0(y)\}$ to a $y$-continuous
 process in the topology of uniform convergence on compact
sets of $y$'s guarantees that 
$$\lim_{n\to\infty} P\left(\mbox{ $|z^n_{[ny]}(0)-nu_0(y)|>n\e$
for some $y\in[c,b]$ }\right) =0.
$$
Under assumptions (\ref{assrho}) and (\ref{assloqeq}) 
 use elementary large deviation
estimates after a partitioning: if $c=b_0<b_1<\cdots<b_k=b$
is a fine enough partition, monotonicity of both 
$z^n_{[ny]}(0)$ and $nu_0(y)$, and the Lipschitz
continuity of  $u_0(y)$, give
\beas
&&P\left(\mbox{ $|z^n_{[ny]}(0)-nu_0(y)|>n\e$
for some $y\in[c,b]$ }\right)\\
&\le& \sum_{j=0}^k P\left(\,|z^n_{[nb_j]}(0)-nu_0(b_j)|>n\e/2\right).
\eeas
These probabilities are summable over $n$, by large 
deviation bounds for exponential random variables. 
\qed

\section{Proof of Theorem \ref{probthm}}
\label{pfofprobthm}

\subsection{Proof of Theorem \ref{probthm}(i)}
\label{pfofprobthmi}

\begin{lm} Suppose $X$ is a measurable function
defined on some measurable space $(\Omega,{\cal F})$,
and 
$C$ is a compact subset of $\RR$. Then there exists
a  measurable function $Y$ such that, for all 
$\omega\in\Omega$,  $Y(\omega)\in C$ and $\dist(X(\omega),C)
=|X(\omega)-Y(\omega)|$.
\label{distrvlm}
\end{lm}

{\it Proof.} The function $g(x)=\inf\{y\in C: \dist(x,C)=|x-y|\}$
is nondecreasing, hence Borel measurable. Set 
$Y(\omega)=g(X(\omega))$. 
\qed

Recall the definition (\ref{inxtdef})
of $i_n(x,t)$. By Lemma \ref{distrvlm}, we may choose a 
random $y_n(x,t)\in I(x,t)$ such that 
$$|n^{-1}i_n(x,t) -y_n(x,t)|=
\dist\left(n^{-1}i_n(x,t), I(x,t)\right).
$$ 
To 
prove the limit (\ref{problim}) in Theorem \ref{probthm},
we 
 bound $\zeta_n(x,t)-\inf_{y\in I(x,t)}\zeta_n(y,0)$ 
 from below and from above, uniformly over
$(x,t)\in A$, with four separate arguments
for different ranges of $t$. Let $[a,b]\times[0,\tau]$ 
be a compact rectangle that contains $A$.


\subsubsection{ Lower Bound, Case 1} 
\label{lbcase1}
Consider $t\in(0,n^{-(1+\delta)}]$ for a 
small $\delta>0$. By Lemma \ref{Hn01lm}(i) we may
condition on the event $H_{n,1}^c$, and thereby assume  
 that
 $i_n(x,t)\ge nx-\alpha n^{(1-\delta)/2}$ for all 
$(x,t)\in [a,b]\times (0,n^{-(1+\delta)}]$.
 Since the $\Gamma$-term 
is always nonnegative, 
$$z^n_{[nx]}(nt)\ge z^n_{i_n(x,t)}(0)\ge 
z^n_{[nx-\alpha n^{(1-\delta)/2}]}(0). $$
 Furthermore, by Lemma \ref{Ixtlm}
and by the local Lipschitz property of $u_0$, there 
exists a constant $C$ such that 
$$u_0(y)\ge u_0(x)-Ct\ge u_0(x)-Cn^{-(1+\delta)}$$
for all $x\in[a,b]$, $t\in(0,n^{-(1+\delta)}]$, and $y\in I(x,t)$. 
Monotonicity in time and space give
$u(x,t)\le u_0(x)$, and  $z^n_{[ny]}(0)\le z^n_{[nx]}(0)$
whenever $y\in I(x,t)$. We get the following
lower bound, valid on the event  $H_{n,1}^c$ for $y\in I(x,t)$: 
\beas
&&z^n_{[nx]}(nt)-nu(x,t)-\{z^n_{[ny]}(0)-nu_0(y)\}\\
&\ge& -\{ z^n_{[nx]}(0)-z^n_{[nx-\alpha n^{(1-\delta)/2}]}(0)\}
-Cn^{-\delta}.
\eeas
Add and subtract the term 
$nu_0(x)-nu_0(x-\alpha n^{-(1+\delta)/2})$,
which is of order $O(n^{1/2-\delta/2})$ uniformly over $x\in[a,b]$
by the local Lipschitz property of $u_0$. 
Multiply through  by $n^{-1/2}$ and uniformize over 
 $(x,t)$:
\beas
&&\inf\{ \zeta_n(x,t)-\zeta_n(y,0): (x,t)\in [a,b]\times 
(0,n^{-(1+\delta)}],\,
y\in I(x,t)\}\\
&\ge& -\sup\{ \zeta_n(x,0)- \zeta_n(y,0): x\in [a,b],\,|x-y|\le  
 2\alpha n^{-(1+\delta)/2}\}
-C n^{-\delta/2}. 
\eeas
This bound is valid on the event  $H_{n,1}^c$,
hence by Lemma \ref{Hn01lm} with probability $1-\e$ if $n$ is large
enough. The lower bound converges to 0 in probability 
by assumption (\ref{ass1d}). 
The constant $\alpha$ was replaced by $2\alpha$ to account
for the effects of integer parts.

\subsubsection{ Lower Bound, Case 2} 
\label{lbcase2}

Now $t\in[n^{-(1+\delta)},\tau]$.
\beas
&&z^n_{[nx]}(nt)-nu(x,t)\\
&=&z^n_{i_n(x,t)}(0)+
\Gamma^{n,i_n(x,t)}_{[nx]-i_n(x,t)}(nt)
-nu(x,t)\\
&=&\left\{ z^n_{[ny_n(x,t)]}(0) -nu_0(y_n(x,t))\right\} 
+\left\{ \Gamma^{n,i_n(x,t)}_{[nx]-i_n(x,t)}(nt)
- (nx-i_n(x,t))^2/4tn \right\}\\
&&\qquad+\left\{z^n_{i_n(x,t)}(0)- z^n_{[ny_n(x,t)]}(0)
-nu_0(n^{-1}i_n(x,t)) +nu_0(y_n(x,t))\right\} \\
&&\qquad+n\left\{ \Phi(x,n^{-1}i_n(x,t))-
\Phi\left(x,y_n(x,t)\right) \right\}.
\eeas
Above  we used the notation  
$$\Phi(x,y)=u_0(y)+tg\left(\frac{x-y}t\right)$$
for the function minimized in the Hopf-Lax formula (\ref{hopflax}). 
Since $y_n(x,t)\in I(x,t)$ minimizes $\Phi(x,\cdot)$, the term  
$\Phi\left(x,n^{-1}i_n(x,t)\right)-\Phi\left(x,y_n(x,t)\right)$
is nonnegative and can be discarded. 
Recalling the definition (\ref{zetanxtdef}) of $\zeta_n$, we get 
\bea
&&\zeta_n(x,t)=n^{-1/2}\left\{ z^n_{[nx]}(nt)-nu(x,t)\right\}\nn\\
&\ge& 
\inf_{y\in I(x,t)} \zeta_n(y,0)  \; +\; 
n^{-1/2}\left\{ \Gamma^{n,i_n(x,t)}_{[nx]-i_n(x,t)}(nt)
- (nx-i_n(x,t))^2/4tn \right\}\label{interm1}\\
 &+&\;
n^{-1/2}\left\{z^n_{i_n(x,t)}(0)-nu_0(n^{-1}i_n(x,t))\right\}
- n^{-1/2}\left\{z^n_{[ny_n(x,t)]}(0) -nu_0(y_n(x,t))\right\}.
\nn
\eea

Recall the definitions of the events $H_{n,0}$ and $H_n(\delta)$
 in  Lemmas \ref{Hn01lm} and \ref{Hnlm}.
By Lemma  \ref{Hnlm}, $\lim_{n\to\infty} P(H_n(\delta))=0$ for 
any fixed $\delta>0$. Then it is possible
to find  a sequence $\delta_n\searrow 0$ 
such that  $\lim_{n\to\infty} P(H_n(\delta_n))=0$. 
Now condition on the event  $H_{n,0}^c\cap H_n(\delta_n)^c$,
the complement  of 
these events. Then 
 for all $(x,t)\in A$ such that  
$t\in [n^{-(1+\delta)},\tau]$, 
\be
\mbox{$[nx]-i_n(x,t)\le \alpha nt^\gamma$ and 
$|n^{-1}i_n(x,t)-y_n(x,t)|\le\delta_n$.}
\label{caseineq1}
\ee
Consequently we get, on the event  $H_{n,0}^c\cap H_n(\delta_n)^c$,
 for all $(x,t)\in A$ such that 
$t\in [n^{-(1+\delta)},\tau]$, 
$$
\zeta_n(x,t)- \inf_{y\in I(x,t)} \zeta_n(y,0) 
\ge  R_{n,1} + R_{n,2}
$$
where we abbreviated
\beas
R_{n,1}&=&\inf\left\{
n^{-1/2}\left(\Gamma^{n,i}_{m}(nt)
-{m^2}/{(4tn)}\right) :
[nc]\le i\le [nb], \right.\\ 
&&\qquad\qquad\qquad\qquad\qquad\qquad\qquad
 \left. 0\le m\le \alpha nt^\gamma,
t\in [n^{-(1+\delta)},\tau]\right\}
\eeas
and 
$$R_{n,2}=\inf
\left\{ \zeta_n(r,0)-\zeta_n(s,0) :
\mbox{ $|r-s|\le\delta_n$ and $r,s\in[c,d]$}\right\}.
$$
In the definition of $R_{n,1}$ and  $R_{n,2}$
 we picked $c<d$ depending on $\alpha$, $\gamma$,
and $\delta_n$ in (\ref{caseineq1}) to ensure that
$y_n(x,t)$ and 
$n^{-1}i_n(x,t)\in[c,d]$ for all $(x,t)\in A$ 
and for all $n$.

\begin{lm} For any $\e>0$,  
$\sum_{n=1}^{\infty} P\left( R_{n,1} \le - \e\right)<\infty. $ 
\label{Rn1lm}
\end{lm}

{\it Proof.} The case $m=0$ in the definition of $R_{n,1}$
can be ignored.  For $1\le m\le \alpha nt^\gamma$ and 
 $t\ge n^{-(1+\delta)}$, 
$$
\frac{ m^{4/3}\log n}{4nt}\le Cn^{1/3-(1+\delta)(4\gamma/3-1)}\log n.
$$
This is less than $\e n^{1/2}$ for large $n$, if we choose 
$\gamma$ close enough to $3/4$ and $\delta$ small enough. This
can be done while satisfying the hypotheses of Lemma \ref{Hn01lm}. 
Consequently the estimate in Lemma \ref{largetlm1}
 is valid for the 
entire  range of $m$-values in the definition of $R_{n,1}$,  
and gives 
\beas
&&\sum_{n=1}^{\infty} P\left( R_{n,1} \le -\e
\right) \\
 &\le&C\sum_{n=1}^{\infty} n \sum_{m=1}^{[nr]}
P\left\{\, \mmGa(m,nt)
\le \frac{m^2}{4tn}  - \frac{ m^{4/3}\log n}{4nt}\  
\mbox{for some
$t\in[n^{-1-\delta},\tau]$} \,\right\}\\
&<&\infty.
\eeas
The factor $Cn$ came from summing over 
$[nc]\le i\le [nb]$ as required by the definition of
$R_{n,1}$, and $r$ was chosen sufficiently large so that 
$nr\ge  \alpha nt^\gamma$ for all $t$ in the range.
\qed

\begin{lm} $\lim_{n\to\infty}
|R_{n,2}|= 0$ in probability. 
\label{Rn2lm}
\end{lm}

{\it Proof.} Recall the definition of
$D_u(\RR)$ as the space of RCLL functions with the locally
uniform metric.  For a fixed $\beta>0$ 
define the continuous  function  $\phi_\beta$ on $D_u(\RR)$
by  
\be\phi_\beta(f)=
\sup
\left\{\left| f(r)-f(s)\right|:
\mbox{ $|r-s|\le\beta$ and $r,s\in[c,d]$}\right\}.
\label{fibetadef}
\ee
For large enough $n$ so that $\delta_n<\beta$, $|R_{n,2}|\le \phi_\beta(\zeta_n(\cdot,0))$.
The set $\{f\in D_u(\RR): 
 \phi_\beta(f)\ge \e\}$ is closed, so  
by the weak convergence $\zeta_n(\cdot,0)\to \zeta_0$ on
 $D_u(\RR)$, 
$$
\limsup_{n\to\infty} P\left( |R_{n,2}|\ge \e\right)
\le
\limsup_{n\to\infty} P\left( \phi_\beta(\zeta_n(\cdot,0))\ge \e\right)
\le 
P\left( \phi_\beta(\zeta_0)\ge \e\right).
$$
Since $\zeta_0$ has continuous paths by assumption (\ref{ass1d}), 
the events $ \phi_\beta(\zeta^0)\ge \e$ decrease to the null event
as $\beta\searrow 0$ and $\e,c,d$ are held fixed.
\qed

We now have for {\bf Case 2}
\beas
&&\lim_{n\to\infty} P\left( \inf_{(x,t)\in A\,,\, t\in[n^{-(1+\delta)},\tau]}
\left\{ \zeta_n(x,t)-\inf_{y\in I(x,t)}\zeta_n(y,0)\right\} \le -\e\right)\\
&\le&\lim_{n\to\infty}
\left\{ P(H_{n,1})+P(H_{n}(\delta_n))+P\left( R_{n,1} \le - \e/2\right)
+P\left( R_{n,2}\le -\e/2\right) \right\}=0
\eeas
for any $\e>0$. {\bf Cases 1} and {\bf 2} together give
$$
\lim_{n\to\infty} P\left( \inf_{(x,t)\in A}
\left\{ \zeta_n(x,t)-\inf_{y\in I(x,t)}\zeta_n(y,0)\right\} \le -\e\right)
=0. 
$$
This completes the proof of the lower bound.
Next we bound 
$\zeta_n(x,t)-\inf_{y\in I(x,t)}\zeta_n(y,0)$ from above.
 
\subsubsection{ Upper Bound, Case 1} 
\label{ubcase1}

Consider  $0<t\le  n^{-1}$.  
Let $C$ be a finite constant such that $u(x,t)\ge u_0(x)-Ct$ 
for all $a\le x\le b$, $0\le t\le \tau$. 
Since $z^n_{[nx]}(nt)\le z^n_{[nx]}(0)$, we can write
\beas
&&z^n_{[nx]}(nt)-nu(x,t)-\{z^n_{[ny]}(0)-nu_0(y)\}\\
&\le& z^n_{[nx]}(0)-nu_0(x)-\{z^n_{[ny]}(0)-nu_0(y)\}+Ctn.
\eeas
Estimate this uniformly over $a\le x\le b$, $0<t\le  n^{-1}$. 
By Lemma \ref{Ixtlm}, there is a constant
$\gamma$ such that $y\in I(x,t)$ implies $|x-y|\le \gamma t\le 
\gamma n^{-1} $, 
for all $(x,t)$ in this range. 
Dividing by $\sqrt{n}$ above gives, for $n$ large enough
to have $x-\gamma  n^{-1}\ge a-1$, 
\beas
&&\sup\{ \zeta_n(x,t)-\zeta_n(y,0): a\le x\le b,\, 0<t\le  n^{-1},\,
y\in I(x,t)\}\\
&\le& \sup\{ \zeta_n(x,0)- \zeta_n(y,0): x,y\in [a-1,b],\,|x-y|\le  
\gamma n^{-1}\}
+C n^{-1/2}.
\eeas 
The last quantity converges to 0 in probability 
by assumption (\ref{ass1d}).

\subsubsection{ Upper Bound, Case 2} 
\label{ubcase2}

Lastly consider $  n^{-1}\le t\le \tau$.
Use the fact that $u(x,t)=u_0(y)+(x-y)^2/(4t)$ for 
any $y\in I(x,t)$. 
\beas
&&z^n_{[nx]}(nt)-nu(x,t) 
=\inf_{y\le x}\{z^n_{[ny]}(0)+
\Gamma^{n,[ny]}_{[nx]-[ny]}(nt)
-nu(x,t)\}\\
&\le& 
\inf_{y\in I(x,t)}\left\{ z^n_{[ny]}(0) -nu_0(y)\right\} 
+ R_{n,3} +C,
\eeas
where 
$$
R_{n,3}=\sup_{x\in[a,b], t\in[ n^{-1},\tau]}
\sup_{y\in I(x,t)}\left\{ \Gamma^{n,[ny]}_{[nx]-[ny]}(nt)
-\frac{([nx]-[ny])^2}{4nt}\right\}, 
$$
and the constant $C$ accounts for replacing $(x-y)^2/(4t)$
with $([nx]-[ny])^2/(4nt)$. 
Consequently  
$$
\sup_{ (x,t)\in A,\,   n^{-1}\le t\le\tau }
\left\{ \zeta_n(x,t)-\inf_{y\in I(x,t)}\zeta_n(y,0)
\right\}
\le  n^{-1/2}R_{n,3} +Cn^{-1/2}. 
$$
The required upper bound follows by
taking $\alpha\in(1/3,1/2)$ in the next lemma. 
 
\begin{lm} For any $\alpha>1/3$, 
$\sum_{n=1}^\infty P\left( R_{n,3}>n^\alpha\right) <\infty.
$
\label{Rn3lm}
\end{lm}

{\it Proof.} $R_{n,3}>n^\alpha$ implies that
\be
\mbox{for some $y\in[c,b]$, $t\in[ n^{-1},\tau]$, and 
$1\le m\le n\gamma t+1$,} \ \Gamma^{n,[ny]}_{m}(nt)
>\frac{m^2}{4nt}+n^\alpha.
\label{Rn3event1}
\ee
To see this, choose $\gamma$ according to Lemma \ref{Ixtlm}
 so that for $y\in I(x,t)$, 
$m=[nx]-[ny]\le n\gamma t+1$. Choose $c\le a-\gamma\tau$
so that  the range of possible 
$y$-values is contained in $[c,b]$.   The case $m=0$ is empty
because $\Gamma^{n,j}_0(nt)\equiv 0$. 

By the assumption $\alpha>1/3$, 
 the inequality
$n^\alpha\ge m^{4/3}(4tn)^{-1}\log n$ is valid  
for the range $1\le m\le n\gamma t+1$, for large enough $n$.
Pick $r>0$ so that $n\gamma \tau+1\le nr$, and let 
$M_0$ be the constant that appeared in Lemma \ref{largetlm2}. 
We can then assert that the event in (\ref{Rn3event1}) is contained in the
union of 
$$
\bigcup_{ \scriptstyle\begin{array}{cc} nc\le j\le na\\ M_0(\log n)^{3/2}\le m\le nr
\end{array} } 
 \left\{ \Gamma^{n,j}_m(nt)> 
\frac{m^2}{4nt} +\frac{m^{4/3}\log n}{4tn} \ \mbox{for some
 $t\in[ n^{-1},\tau]$} \right\}
$$
and 
$$
\bigcup_{nc\le j\le na}  
 \left\{ \Gamma^{n,j}_m(nt)> n^\alpha 
\ \mbox{for some
 $t\in[ n^{-1},\tau]$ and $1\le m\le M_0(\log n)^{3/2}$} \right\}.
$$
The conclusion now follows from the estimates in Lemmas 
\ref{largetlm2} and \ref{largetlm3},
because for any fixed $(n,j)$, $\Gamma^{n,j}_m(nt)$ has the 
same distribution as $\mmGa(m,nt)$. 
\qed

Combining Cases 1 and 2,  we have bounded 
$$
\sup_{(x,t)\in A}\left\{ \zeta_n(x,t)-\inf_{y\in I(x,t)} 
\zeta_n(y,0)\right\} 
$$
above by a random variable that 
  vanishes in
probability as $n\to\infty$. Together with the lower
bound, this completes the proof
of part (i) of  Theorem \ref{probthm}.

\subsection{Proof of Theorem \ref{probthm}(ii)}
\label{pfofprobthmii}
First a lemma
whose proof is partly a repetition of the above argument.  

\begin{lm} Let $-\infty<a<b<\infty$ and $0<\tau<\infty$, and set
$$
M_n=\sup_{x\in[a,b], t\in[0,\tau]}\left| \zeta_n(x,t)-\inf_{y\in I(x,t)}
\zeta_n(y,0)\right|.
$$
Then for every $\e>0$ there exists a finite constant $C$ such that
$P(M_n\le C)\ge 1-\e$ for all $n$. 
\label{suplm}
\end{lm}

{\it Proof.} We prove the measurability of 
$M_n$ in the Appendix.
 Let $A_0=[a,b]\times[0,\tau]$. If $A_0$ were an admissible compact set 
for the proof of Theorem \ref{probthm}(i) just completed, there would
be nothing more to prove. But $A_0$ might have shocks. The only
step where this makes a difference in the above proof is 
Section \ref{lbcase2} 
{\bf Lower Bound,  Case 2} because this case appealed to Lemma \ref{Hnlm}.
Without Lemma \ref{Hnlm} the argument  still gives a lower bound. Condition
on the event $G_n^c$ so that $i_n(x,t)\ge nc$  for $(x,t)\in A_0$. 
Then inequality (\ref{interm1}) gives 
\bea
&&\inf_{(x,t)\in A_0}\left\{\zeta_n(x,t)-
\inf_{y\in I(x,t)}\zeta_n(y,0)\right\} 
\ge-2\cdot\sup_{y\in [c,b]}|\zeta_n(y,0)|\nn\\
&&\quad
+\inf_{1\le m\le n(b-c)\,,\, t\in[n^{-(1+\delta)},\tau]\,,
\, nc\le j\le nb}\left\{ 
\Gamma_{m}^{n,j}(nt)-\frac{m^2}{4nt}\right\}.
\label{interm2}
\eea
This  bound is valid on the event $G_n^c$, hence with probability
$1-\e$ for large enough $n$. 
Combined with the other cases proved in Section
\ref{pfofprobthmi} it 
proves the lemma.
\qed

Fix $-\infty<a<b<\infty$ and $\tau<\infty$. The goal is to prove the
limit in probability (\ref{problim2}). Let $\e>0$. 
Let
$$
Y_n=\sup_{(x,t)\in[a,b]\times[0,\tau]}\left| \zeta_n(x,t)-\inf_{y\in I(x,t)}\zeta_n(y,0)
\right|^p.
$$
By Lemma \ref{suplm} we can find a constant $C\in(0,\infty)$ such that 
$P(Y_n> C)\le \e$ for all $n$. Define  the event $D_n=\{ Y_n\le C\}$. 
By Proposition \ref{Gprop} we can find an open set $G\subset\RR\times(0,\infty)$
such that $G$ contains all the shocks in $[a,b]\times[0,\tau]$, and 
its $t$-section has 1-dimensional Lebesgue measure $|G_t|<\e/(3C)$ for all $t$. 
(Recall that by definition there are no shocks on the $t=0$ line.)

Let $A=[a,b]\times[0,\tau]\setminus G$. $A$ is a compact set 
with no shocks, so by Theorem \ref{probthm}(i) 
\be
X_n\equiv \sup_{(x,t)\in A}\left| \zeta_n(x,t)-\inf_{y\in I(x,t)}\zeta_n(y,0)
\right|^p \to 0 \qquad\mbox{in probability.}
\label{Xnlim}
\ee
Let $A_t=\{x: (x,t)\in A\}$ be the $t$-section of $A$. 
On the event $D_n$ we can now bound
\beas
&&\sup_{0\le t\le\tau} \int_a^b \left| \zeta_n(x,t)-\inf_{y\in I(x,t)}\zeta_n(y,0)
\right|^p dx\\
&\le& \sup_{0\le t\le\tau} \left\{ \int_{A_t} X_n\, dx 
+\int_{[a,b]\setminus A_t} Y_n\, dx\right\}\\
&\le& (b-a)X_n + Y_n \e/(3C) \le (b-a)X_n+ \e/3. 
\eeas
Thus by (\ref{Xnlim})
$$
\limsup_{n\to\infty} P\left\{ \sup_{0\le t\le\tau} \int_a^b 
\left| \zeta_n(x,t)-\inf_{y\in I(x,t)}\zeta_n(y,0)
\right|^p dx \ge \e\right\} \le P(D_n^c) \le \e.
$$ 
This proves (\ref{problim2}). 

\section{Proof of the weak limit and the linearized equation} 
\label{sectpfweaklimitthm}

\subsection{Proof of Theorem \ref{weaklimitthm}}

For part (i), take $A=\{(x_1,t_1),\ldots,(x_k,t_k)\}$ in  (\ref{problim}).
Note that the mapping 
$$h\mapsto\left( \inf_{y\in I(x_1,t_1)}h(y),\ldots, \inf_{y\in I(x_k,t_k)}h(y)\right)$$
from $D_u(\RR)$ into $\RR^k$ is continuous.  Then 
use the assumption (\ref{ass1d}) of weak convergence at time zero, and the
continuous mapping theorem \cite[p.\ 30]{billingsley}. 

Part (ii) goes by the same general principle. Let us abbreviate 
\be
\sigma_n(x,t)=\inf_{y\in I(x,t)}\zeta_n(y,0).
\label{sigmadef}
\ee
We need to check that $\zeta_n$ defines a random element of
$D\left([0,\infty),\Lploc\right)$ and that $\sigma_n$ and $\zeta$ 
define random elements of
$C\left([0,\infty),\Lploc\right)$. 

For $f\in D_u(\RR)$ let $Gf(x,t)=\inf_{y\in I(x,t)}f(y)$. 

\begin{lm} $G$ is a continuous map 
from $D_u(\RR)$ into $C\left([0,\infty),\Lploc\right)$, when 
we interpret $Gf$ as the path $t\mapsto Gf(\cdot,t)\in\Lploc$.
\label{GLplm}
\end{lm}

{\it Proof.} 
For $0\le t\le T$ and $a\le x\le b$, $I(x,t)\subseteq [y^-(a,T),b]$, and so 
 $Gf$ is locally bounded as a function of $(x,t)$.  
Consequently for a fixed $t$, $Gf(\cdot,t)$
is in $\Lploc$. 

Secondly, we need to argue that
 as $s\to t$, $Gf(\cdot,s)\to Gf(\cdot,t)$ in $\Lploc$. By the local
boundedness and dominated convergence, we need only
 show   
\be
\mbox{$Gf(x,s)\to Gf(x,t)$ for a.e.\ $x$. }
\label{Gaeconv}
\ee
Recall from Section \ref{charsection} that if $y_1\in I(x,s_1)$ and $y_2\in I(x,s_2)$ for
$s_1<s_2$, then $y_2\le y_1$. Consider first $s\nearrow t$. Fix $x$ so that
$(x,t)$ is not a shock. Then for any choice $y_s\in I(x,s)$, 
$y_s\searrow y(x,t)$ $=$ the unique Hopf-Lax minimizer for $(x,t)$.
By right-continuity $f(y_s)\to f(y(x,t))=Gf(x,t)$, and since we can 
let $f(y_s)$ be arbitrarily close to $Gf(x,s)$, we have (\ref{Gaeconv}). 

Now suppose $s\searrow t$. $f$ has at most countably many 
discontinuities, so we still have a.e.\ $x$ if we exclude all $x$
such that $(x,t)$ is a shock, and all points $x=w^\pm(\ybar,t)$ 
for discontinuities $\ybar$ of $f$. Suppose $x$ is not one of the 
excluded points.  Then $(x,t)$ has a unique minimizer $y(x,t)$,
and the previous paragraph shows again $f(y_s)\to f(y(x,t))$ if
$f$ is continuous at $y(x,t)$. But suppose $\ybar=y(x,t)$ is a discontinuity
for $f$. Then it must be that $w^-(\ybar,t)<x<w^+(\ybar,t)$. 
[Justification: (\ref{xyequiv}) forces  $w^-(\ybar,t)\le x\le w^+(\ybar,t)$, 
but $x\in\{w^\pm(\ybar,t)\}$ cannot happen because 
$x$ is not among the excluded points.] Since forward characteristics
are continuous, it follows that for $s>t$ but close enough to $t$,
 $w^-(\ybar,s)<x<w^+(\ybar,s)$. This implies $y^\pm(x,s)=\ybar$
which in turn says $Gf(x,s)=Gf(x,t)$. Again (\ref{Gaeconv}) checks. 

We have now shown   that the function $t\mapsto Gf(\cdot,t)\in\Lploc$ 
is continuous. Finally, we check that the map 
$G:D_u(\RR)\to C\left([0,\infty),\Lploc\right)$ is continuous.
This is a consequence of having the locally uniform topology 
on $D_u(\RR)$. For by the observation made in the beginning of the proof, 
$$
\sup_{0\le t\le T} \int_{a}^b | Gf(x,t)-Gg(x,t)|^pdx
\le \sup_{y^-(a,T)\le x\le b}|f(x)-g(x)|.  \qquad\qed
$$

By definitions (\ref{zetaxtdef}) and (\ref{sigmadef}), the processes
$\zeta$ and $\sigma_n$ are obtained by applying the mapping $G$
to the $D_u(\RR)$-valued random functions $\zeta_0$ and $\zeta_n(\cdot,0)$. 
This checks that $\sigma_n$ and $\zeta$
define random elements of
$C\left([0,\infty),\Lploc\right)$. Also, by assumption (\ref{ass1d})
and the continuous mapping theorem, $\sigma_n\stackrel{d}\to\zeta$
in the space $C\left([0,\infty),\Lploc\right)$. 

Now consider $\zeta_n$ defined by (\ref{zetanxtdef}). It is jointly measurable
in $(x,t,\omega)$, where $\omega$ is a sample point of the underlying 
probability space $\Omega$ (see the Appendix).
 It is also locally bounded in $(x,t)$
so local $L^p$-integrability is not  a problem.
  Fix $t$. First we argue that 
\be
\mbox{$\omega\mapsto\zeta_n(\cdot, t;\omega)$ 
is a measurable map from $\Omega$ into $\Lploc$,}
\label{Lpmeasurab}
\ee
 where $\Lploc$ is
endowed with its Borel $\sigma$-algebra defined by the metric $d_p$ 
of (\ref{dpdef}). 
By Fubini's theorem, for any $f\in \Lploc$, 
$d_p(f, \zeta_n(\cdot,t;\omega))$ is a measurable function of $\omega$. 
Hence for any open $d_p$-ball $B(f,r)$ in $\Lploc$, the inverse image
$\{\omega: \zeta_n(\cdot,t;\omega)\in B(f,r)\}$ is measurable. 
By the separability of $\Lploc$, (\ref{Lpmeasurab}) follows. 

By convention interacting systems are constructed to be right-continuous 
in time, so $t\mapsto \zeta_n(x,t;\omega)$ is right-continuous. 
By local boundedness and dominated convergence, 
the map $t\mapsto \zeta_n(\cdot,t;\omega)$ from $[0,\infty)$
into $\Lploc$ is right-continuous. The measurability 
of $\zeta_n(\cdot,\cdot;\omega)$ as a $D\left([0,\infty),\Lploc\right)$-valued
random element follows because in $D$-space measurability is 
equivalent to measurability of the time-coordinate projections
$\zeta_n(\cdot,t;\omega)$.

The weak convergence
$\zeta_n\stackrel{d}\to\zeta$ now follows readily. For any finite
time-horizon $T$, 
$\sup_{0\le t\le T} d_p(\zeta_n(t),\sigma_n(t))\to 0$ in probability
    by (\ref{problim2}). Uniform in time is stronger than the
Skorokhod topology. So it follows that, if we let $d_D$ denote
the Skorokhod metric on $D\left([0,\infty),\Lploc\right)$, 
$d_D(\zeta_n,\sigma_n)\to 0$ in probability also. This together with
$\sigma_n\stackrel{d}\to\zeta$ implies $\zeta_n\stackrel{d}\to\zeta$.
We have proved Theorem \ref{weaklimitthm}. 


\subsection{ Proof of Theorem \ref{lineqnthm}}


\begin{lm} Let $F,G$ be right-continuous functions,
$F$ locally BV and $G$  nondecreasing. Let  
$H^-$ the left-continuous inverse of $G$ defined by
$$
H^-(y)=\sup\{x: G(x)<y\}=\inf\{x: G(x)\ge y\}.
$$
Then for all continuous functions $\varphi$ for which the integrals
exist, 
$$
\int \varphi(H^-(y))dF(y)=\int \varphi(x) d(F\circ G)(x).
$$
\label{FGHlm}
\end{lm}

{\it Proof.} It suffices to take $\varphi={\bf 1}_{(a,b]}$, the indicator
function of a left-open right-closed interval. Check that 
$
\{y: a<H^-(y)\le b\}=(G(a),G(b)].
$
Then 
\beas
\int {\bf 1}_{(a,b]}(H^-(y))dF(y)&=&\int {\bf 1}_{(G(a),G(b)]}(y)dF(y)=
F(G(b))-F(G(a))\\
&=&\int {\bf 1}_{(a,b]}d(F\circ G). \qquad\qed
\eeas

This lemma will be applied below to the pair $G(a)=w^+(a,t)$,  
$H^-(b)=y^-(b,t)$. 

Fix a test function $\phi\in C^\infty_c(\RR\times[0,\infty))$. 
Let $(A,B)\times[0,T)$ contain the support of $\phi$. 
Let $\Phi(x,t)=\int_{-\infty}^x\phi(y,t)dy$. By Theorem 
\ref{odethm}, for any $q\in\RR$ we have the formula 
\bea
-\Phi(q,0)&=&\int_0^T \frac{d}{dt}\Phi(w^+(q,t),t)dt\nn\\
&=&\int_0^T \Phi_t(w^+(q,t),t)dt + \int_0^T \phi(w^+(q,t),t)h(w^+(q,t),t)\,dt.
\label{Phieqn}
\eea
Now we calculate, beginning with the leftmost 
term of (\ref{weaktranseqn}), with $v$ in place of
$\zetabar$. Note that 
$v(x,t)=v_0(y^-(x,t))$ a.e.\ so in this first integral these two
are interchangeable.  
\beas
\int_0^T dt\int  v(x,t)\phi_t(x,t) dx
&=& \int_0^T dt\int v_0(y^-(x,t)) d[\Phi_t(\cdot,t)](x)\\
&=& \int_0^T dt\int v_0(q) d[\Phi_t(w^+(\cdot,t),t)](q).
\eeas
There is a fixed compact interval $[a,b]$ on which the 
Lebesgue-Stieltjes   measure \break
 $d[\Phi_t(w^+(\cdot,t),t)](q)$ is
supported for all $t\in[0,T]$. Let $a=q_0<q_1<\cdots<q_m=b$ be a partition of this
interval with mesh $\Delta=\max(q_i-q_{i-1})$. 
We can choose the partitions so that $w(q_i,t)=w^\pm(q_i,t)$ 
for all $i$, because by Lemma \ref{charlm1}(c) we only need 
to pick the $q_i$'s outside a certain Lebesgue null set. 
 The integrand $v_0$ is
continuous by assumption, hence the $q$-integral can be 
written as a limit, and 
the last line above equals
$$
=\int_0^T dt \lim_{\Delta\to 0}\sum_i v_0(q_i) 
\{\Phi_t(w(q_i,t),t) - \Phi_t(w(q_{i-1},t),t)\}.
$$
The function inside the $t$-integral is bounded by a constant, uniformly
over $t\in[0,T]$ and over partitions of $[a,b]$. Hence we can take the limit
outside, apply (\ref{Phieqn}), and then put the limit back inside, to get
\bea
&=& \lim_{\Delta\to 0}  \sum_i v_0(q_i)
\int_0^T \{\Phi_t(w(q_i,t),t) - \Phi_t(w(q_{i-1},t),t)\} dt\nn\\
&=& 
\lim_{\Delta\to 0} \left\{ -\sum_i v_0(q_i)
\int_0^T \{\phi(w(q_i,t),t)h(w(q_i,t),t)  \right.\nn\\
&&\qquad \left. -\phi(w(q_{i-1},t),t)h(w(q_{i-1},t),t)\} dt
-\sum_i v_0(q_i) [\Phi(q_i,0)-\Phi(q_{i-1},0)]\right\}
\nn\\
&=&-\int_0^T \lim_{\Delta\to 0}  \sum_i v_0(q_i)
\{\phi(w(q_i,t),t)h(w(q_i,t),t)  \nn\\
&&\qquad -\phi(w(q_{i-1},t),t)h(w(q_{i-1},t),t)\} dt
-\int v_0(q)\phi(q,0)dq.
\label{linebeforeR}
\eea
At this point we  replace $h(w(q_i,t),t)$ by $f'(\rho^+(w(q_i,t),t))$
and write $R_\Delta$ for the error term. Then the last line above
equals
\beas
&=&-\int_0^T \lim_{\Delta\to 0}  \sum_i v_0(q_i)
\{\phi(w(q_i,t),t)f'(\rho^+ (w(q_i,t),t))  \\
&& -\phi(w(q_{i-1},t),t)f'(\rho^+ (w(q_{i-1},t),t))\} dt
-\int v_0(q)\phi(q,0)dq +\lim_{\Delta\to 0}R_\Delta.
\eeas 
Ignoring the term $\lim_{\Delta\to 0}R_\Delta$ for the moment, 
take the $\Delta\to 0$ limit in the first sum to get
again a Lebesgue-Stieltjes integral.  
After another application of Lemma \ref{FGHlm}, we get 
this intermediate equation:
\bea
&&\int_0^T dt\int  v(x,t)\phi_t(x,t) dx\nn\\
&=& -\int_0^T dt \int v_0(q) d[\phi(w^+(\cdot,t),t)
 f'(\rho^+ (w^+(\cdot,t),t))](q)
-\int v_0(q)\phi(q,0)dq\nn\\
&&\qquad\qquad  -\lim_{\Delta\to 0}R_\Delta\nn\\
&=&-\int_0^T dt \int v_0(y^-(x,t)) d[\phi(\cdot,t) f'(\rho^+ (\cdot,t))](x)
-\int v_0(x)\phi(x,0)dx\nn\\
&&\qquad\qquad  -\lim_{\Delta\to 0}R_\Delta.
\label{insertRhere}
\eea
It remains to take care of  $\lim_{\Delta\to 0}R_\Delta$. 
Notice that on line (\ref{linebeforeR}), 
at the stage where $R_\Delta$ was
introduced, the summation can be restricted to $i$ such that
$w(q_{i-1},t)<w(q_i,t)$ because otherwise $w(q_{i-1},t)=w(q_i,t)$
and the expression in braces $\{ \}$ equals zero. Thus
we can write $R_\Delta$ as follows, and sum by parts: 
\beas
R_\Delta&=&\int_0^T dt\sum_{i: w(q_{i-1},t)<w(q_i,t)}
v_0(q_i)\left[ \phi(w(q_i,t),t) \left\{
h(w(q_i,t),t)- f'(\rho^+(w(q_i,t),t)) \right\} \right.\\
&& \left. -\phi(w(q_{i-1},t),t) \left\{
h(w(q_{i-1},t),t)- f'(\rho^+(w(q_{i-1},t),t)) \right\} \right]\\
&=&\int_0^T dt \sum_i \phi(w(q_i,t),t)
\left\{ h(w(q_i,t),t)- f'(\rho^+(w(q_i,t),t)) \right\}\\
&&\quad \times \left[ v_0(q_i){\bf 1}\{ w(q_{i-1},t)<w(q_i,t)\}
- v_0(q_{i+1}){\bf 1}\{ w(q_{i},t)<w(q_{i+1},t)\}\right].
\eeas
Now note that the last sum can be restricted to $i$ such 
that $(w(q_i,t),t)$ is a shock because 
$h(x,t)-f'(\rho^+(x,t))=0$ unless $(x,t)$ is a shock. 
Supposing that $(x,t)$ is a shock, observe that 
if $y^-(x,t)\le q_i<q_{i+1}\le y^+(x,t)$ then 
$w(q_{i},t)=w(q_{i+1},t)=x$. [In general $w^+(y^+(x,t),t)$ 
could be strictly larger than $x$, but then 
$w^-(y^+(x,t),t)< w^+(y^+(x,t),t)$, which we have prevented by
assuming 
$w^-(q_i,t)=w^+(q_i,t)$.] Consequently, for the shock 
$(x,t)$, 
\beas
&&\sum_{i:w(q_{i},t)=x}
\left[ v_0(q_i){\bf 1}\{ w(q_{i-1},t)<w(q_i,t)\}
- v_0(q_{i+1}){\bf 1}\{ w(q_{i},t)<w(q_{i+1},t)\}\right]\\
&=&v_0\left(\min\{q_i: q_i\ge y^-(x,t)\}\right)
-v_0\left(\min\{q_i: q_i> y^+(x,t)\}\right).
\eeas
To express this in a single function, write
\beas
L_\Delta(x,t)&=&
{\bf 1}\left[
\mbox{$(x,t)$ is a shock, and $x\in\{w(q_i,t):0\le i\le m\}$}\right]\\
&&\cdot \left[ v_0\left(\min\{q_i: q_i> y^+(x,t)\}\right)
-v_0\left(\min\{q_i: q_i\ge y^-(x,t)\}\right) \right].
\eeas
The subscript $\Delta$ expresses the dependence of $L_\Delta$
on the partition. For any shock $(x,t)$,  
some $q_i$ lies in $(y^-(x,t),y^+(x,t))$ when $\Delta$ is
small enough, and then $x=w(q_i,0)$. 
By the continuity of $v_0$ we have the convergence 
\be
\lim_{\Delta\to 0} L_\Delta(x,t)= v_0(y^+(x,t))-v_0(y^-(x,t)),
\label{LDlim}
\ee
which happens boundedly and at all $(x,t)$. Now write
\beas
R_\Delta
&=&
\int_0^T dt \sum_{a\le x\le b} \phi(x,t)
\left\{   f'(\rho^+(x,t),t) -h(x,t) \right\}
L_\Delta(x,t)\\
&=&\int_0^T dt \int \theta(x,t) L_\Delta(x,t)
d[\phi(\cdot,t)f'(\rho^+(\cdot,t))](x)
\eeas
where we recognized that 
$$
 f'(\rho^+(x,t),t) -h(x,t) = \theta(x,t)[ f'(\rho^+(x,t),t)
-  f'(\rho^-(x,t),t)].
$$
The $x$-integral lives entirely on the countable set 
of shocks because $L_\Delta$ vanishes elsewhere.
This is why we can slip the continuous $\phi(x,t)$ factor
into the integrator. Taking the limit gives
\beas
&&
\lim_{\Delta\to 0} R_\Delta\\
&=&\int_0^T dt \int \theta(x,t) \left[v_0(y^+(x,t))-v_0(y^-(x,t))\right]
d[\phi(\cdot,t)f'(\rho^+(\cdot,t)](x). 
\eeas
Substituting this on line (\ref{insertRhere}) above
completes the proof that $v(x,t)$ satisfies 
(\ref{weaktranseqn}).  We have proved Theorem \ref{lineqnthm}.

\section{Proof of Theorem \ref{thmloqeq}} 
\label{sectpfthmloqeq}
We begin by realizing the initial configurations 
$(z^n_i(0):i\in\ZZ)$ with Skorokhod's 
representation. Let $(\Omega, {\cal F}, P)$ be a  probability space on which are defined
a two-sided Brownian motion $B(\cdot)$, and 
independently of it a  
space-time Poisson point process for constructing
the Hammersley dynamics. Recall that $B(\cdot)$ is defined by
$B(s)=B_1(s)$ for $s\ge 0$ and $B(s)=-B_2(-s)$ for $s< 0$, where
$B_1(\cdot)$, $B_2(\cdot)$ are two independent standard 1-dimensional Brownian 
motions defined on $[0,\infty)$. 
For each $n$, define a two-sided Brownian motion $B_n(\cdot)$
by $B_n(s)=n^{1/2}B(s/n)$ for $s\in\RR$. 

Fix $n$. Construct  the Skorokhod representation for the 
independent mean zero random variables $\eta^n_i-E[\eta^n_i]$ whose
distribution is defined in assumption (\ref{assloqeq}).  The usual
construction (see e.g.
Section 7.6 in \cite{Dur}) is applied  to $B_1$ for $i> 0$
and to $B_2$ for $i\le 0$. This gives  random variables 
$$
\cdots\le T_{n,-2}\le T_{n,-1}\le 0=T_{n,0}\le T_{n,1}\le T_{n,2}\le 
\cdots
$$
such that the variables
 $\{\tau_{n,i}= T_{n,i}-T_{n,i-1}:i\in\ZZ\}$ are mutually
independent,  we have the equality in distribution 
of the processes
$$
\{B_n( T_{n,i})-B_n(T_{n,i-1}):i\in\ZZ\}
\stackrel{d}= \{\eta^n_i(0)-
E[\eta^n_i(0)]: i\in\ZZ\}
$$
and for each $i$ 
\be
E[\tau_{n,i}]=\Var[\eta^n_i(0)]=E[\eta^n_i(0)]^2=\left( n\int^{i/n}_{(i-1)/n}
\rho_0(s)ds\right)^2.
\label{Etau}
\ee
Note that the assumption of exponentially distributed 
$\eta^n_i(0)$ was used here. 

Now we take this construction as the {\it definition} of
the initial interface: 
\be
z^n_i(0)= nu_0(i/n)+B_n(T_{n,i}) =nu_0(i/n)+n^{1/2}B(n^{-1}T_{n,i}).
\label{zconstrB}
\ee
The initial process $\zeta_n(y,0)$ defined by (\ref{zetanxtdef})
is now given by
\be
\zeta_n(y,0)=B(n^{-1}T_{n,[ny]})+n^{1/2}\left(u_0([ny]/n)-u_0(y)\right).
\label{zetaB}
\ee

\begin{lm} For any $-\infty<a<b<\infty$,
$$\lim_{n\to\infty} \sup_{y\in[a,b]}
\left| \frac{T_{n,[ny]}}n - \int_0^y \rho^2_0(s)ds
\right| =0 \qquad\mbox{almost surely.}
$$
\label{Tlm}
\end{lm}

{\it Proof.} Suppose $0\le a<b$. The other cases are
handled with similar arguments. 
Let us first check 
\be
\lim_{n\to\infty} \sup_{y\in[a,b]}
\left| \frac1n{ET_{n,[ny]}} - \int_0^y \rho^2_0(s)ds
\right| =0.
\label{ETlim}
\ee
By (\ref{Etau}), 
\beas
\frac1n{ET_{n,[ny]}}&=&\frac1n\sum_{i=1}^{[ny]} 
\left( n\int^{i/n}_{(i-1)/n}
\rho_0(s)ds\right)^2\\
&=&\int_0^{[ny]/n}
\sum_{i=1}^{[ny]} 
\left( n\int^{i/n}_{(i-1)/n}
\rho_0(s)ds\right)^2
{\bf 1}_{\left[\frac{i-1}n,\frac{i}n\right)}(r) dr.
\eeas
The integrand is bounded by the assumption
$\rho_0\in L^\infty_{\rm loc}(\RR)$, 
and converges to $\rho^2_0(r)$
at every Lebesgue point $r$ of $\rho_0$. So the required
convergence in (\ref{ETlim}) holds for each fixed $y$.  To get uniformity
over $y$, 
\beas
&&\sup_{y\in[a,b]}
\left| \frac1n{ET_{n,[ny]}} - \int_0^y \rho^2_0(s)ds
\right|\\
&\le&\sup_{y\in[a,b]}
\left| \int_0^{[ny]/n} \left\{
\sum_{i=1}^{[ny]} 
\left( n\int^{i/n}_{(i-1)/n}
\rho_0(s)ds\right)^2
{\bf 1}_{\left[\frac{i-1}n,\frac{i}n\right)}(r) 
-\rho^2_0(r)\right\} dr \right|+\frac{C}{n}\\
&\le&\sup_{y\in[a,b]}
 \int_0^{[ny]/n} \left| 
\sum_{i=1}^{[nb]} 
\left( n\int^{i/n}_{(i-1)/n}
\rho_0(s)ds\right)^2
{\bf 1}_{\left[\frac{i-1}n,\frac{i}n\right)}(r) 
-\rho^2_0(r)\right| dr +\frac{C}{n}\\
&&\qquad\qquad\qquad\qquad\qquad\mbox{(note that the integrand no longer
 depends on $y$)}\\
&\le& 
 \int_0^{b} \left| 
\sum_{i=1}^{[nb]} 
\left( n\int^{i/n}_{(i-1)/n}
\rho_0(s)ds\right)^2
{\bf 1}_{\left[\frac{i-1}n,\frac{i}n\right)}(r) 
-\rho^2_0(r)\right| dr +\frac{C}{n}.
\eeas
The error $C/n$ accounts for the effect of switching
between $y$ and 
$[ny]/n$ ($b$ and $[nb]/n$) as the upper limit of integration. 
Now (\ref{ETlim}) follows by dominated convergence. 

Next we show that, for a fixed $y$, 
\be
\lim_{n\to\infty}
\left| n^{-1}{T_{n,[ny]}} -  n^{-1}{ET_{n,[ny]}}
\right| =0\qquad\mbox{a.s.}
\label{Tnylim}
\ee
The moments of the waiting times $\tau_{n,i}$ satisfy 
$$E[(\tau_{n,i})^k]\le C_k E\left[\left\{\eta^n_i(0)-
E[\eta^n_i(0)]\right\}^{2k}\right]
\le C'_k <\infty
$$
for all $0\le i\le nb$, for constants $C_k, C'_k$.
The first inequality follows from the Burkholder-Davis-Gundy
inequalities, and the second from the local boundedness
of $\rho_0$. 
 Thus
\beas
&&P\left( | T_{n,[ny]} -  ET_{n,[ny]} |\ge n\e\right) 
=P\left( \left| \sum_{i=1}^{[ny]}(\tau_{n,i} -  E\tau_{n,i})
 \right|\ge n\e\right) \\
&&\qquad\qquad\le \frac1{n^4\e^4}
E\left[ \left( \sum_{i=1}^{[ny]}(\tau_{n,i} -  E\tau_{n,i})
 \right)^4\,\right] \le \frac{C}{n^2\e^4},
\eeas
and now Borel-Cantelli gives (\ref{Tnylim}). 
Finally, given $\e>0$, pick $\delta>0$ so that 
$$\int_y^x \rho^2_0(s)ds <\e
\qquad
\mbox{for any $a\le y<x<y+\delta\le b$.}
$$
Pick a partition $a=a_0<a_1<\cdots<a_m=b$ such that 
$a_{k+1}-a_k<\delta$. Apply (\ref{Tnylim}) for the values $y=a_k$, 
and between the partition points estimate by 
$$
\inf_{[na_k]\le i\le [na_{k+1}]}
\frac{T_{n,i}-ET_{n,i}}n\ge \frac{T_{n,[na_k]}-ET_{n,[na_k]}}n
+\frac{ET_{n,[na_k]}-ET_{n,[na_{k+1}]}}n\,,
$$
to get, by (\ref{ETlim}) and (\ref{Tnylim}), 
$$\liminf_{n\to\infty} \inf_{[na]\le i\le [nb]}
\frac{T_{n,i}-ET_{n,i}}n\ge
-\max_{0\le k\le m}\int_{a_k}^{a_{k+1}} \rho^2_0(s)ds\ge -\e.
$$
Repeat the argument for the upper bound, and let $\e\searrow 0$
 to get 
$$\lim_{n\to\infty} \sup_{[na]\le i\le [nb]}
\left|n^{-1}T_{n,i}-n^{-1}ET_{n,i} \right|=0 \qquad\mbox{a.s.}
$$
Combine this with (\ref{ETlim}) to get the conclusion of the lemma. 
\qed

By definition (\ref{zconstrB})
and the path-continuity of Brownian motion, Lemma
\ref{Tlm} is sufficient
for proving that 
\be
\lim_{n\to\infty} \sup_{y\in[a,b]}
\left| \zeta_n(y,0)
- B\left(\int_0^y \rho^2_0(s)ds\right)
\right| =0 \qquad\mbox{almost surely.}
\label{limBt0}
\ee
Thus to prove limit  (\ref{limznxtB}) in Theorem
\ref{thmloqeq} it suffices to show 
\be
\lim_{n\to\infty} \sup_{(x,t)\in A}
\left| \zeta_n(x,t)-\inf_{y\in I(x,t)} 
\zeta_n(y,0) \right| =0
\qquad\mbox{almost surely.}
\ee
In other words, we need to strengthen (\ref{problim}) 
to a.s.\ convergence.
The proof follows the  case-by-case reasoning 
in Section \ref{pfofprobthm} for (\ref{problim}).  The 
error terms $R_{n,j}$, $j=1,2,3$, are the same 
as there.  We check that in
each case the stronger assumptions (\ref{assrho})
and  (\ref{assloqeq}) give almost
sure convergence.

\hbox{}

{\it Lower Bound, Case 1.} By the argument for 
the case  $t\in(0,n^{-(1+\delta)}]$
in Section \ref{lbcase1}, 
\beas
&&\sum_{n\ge 1}P\left( \inf_{(x,t)\in A\,,\, t\in(0,n^{-(1+\delta)}]}
\left[ z^n_{[nx]}(nt)-nu(x,t)  \right.\right.\\
&&\qquad\qquad\qquad \left.\left.
-\inf_{y\in I(x,t)}
\{z^n_{[ny]}(0)-nu_0(y)\}\right] \le -2\e n^{1/2}\right)\\
&\le& \sum_{n\ge 1} P(H_{n,1})+\sum_{n\ge 1}
P\left( \sup_{x\in[a,b]} \left\{ 
z^n_{[nx]}(0)-z^n_{[nx-\alpha n^{(1-\delta)/2}]}(0) \right\}
\ge \e n^{1/2} \right).
\eeas
$\sum P(H_{n,1})<\infty$ follows from Lemma \ref{Hn01lm}(ii).
To show that the last sum is finite, write
$$
z^n_{[nx]}(0)-z^n_{[nx-\alpha n^{(1-\delta)/2}]}(0)=
\sum_{i=[nx-\alpha n^{(1-\delta)/2}]+1}^{[nx]} \eta^n_i(0)
$$
in terms of increment, or stick  variables. 
By the local boundedness of $\rho_0$ and by assumption
 (\ref{assloqeq}),
the $\eta^n_i(0)$ are stochastically
dominated by i.i.d.\ exponential variables of 
finite mean. Now apply standard large deviation
estimates.

\hbox{}

{\it Lower Bound, Case 2.} Now $t\in[n^{-(1+\delta)},\tau]$.
 The argument for this case in
Section \ref{lbcase2} showed that on the event 
$H_{n,0}^c\cap H_{n}(\delta_n)^c$, 
$$
 \inf_{(x,t)\in A\,,\, t\in[n^{-(1+\delta)},\tau]}
\left\{ \zeta_n(x,t)-\inf_{y\in I(x,t)}\zeta_n(y,0)\right\}
\ge R_{n,1} + R_{n,2}.
$$
By Lemma \ref{Hnlm}(ii), $\sum P(H_{n}(\delta))<\infty$
for any fixed $\delta>0$. Then it is possible to find a 
sequence $\delta_n\searrow 0$ such that
 $\sum P(H_{n}(\delta_n))<\infty$. For example, 
for $j>0$ find $n_0(j)\nearrow\infty$ such that 
$\sum_{n\ge n_0(j)} P(H_{n}(j^{-1}))< 2^{-j}$. And then 
set $\delta_n=j^{-1}$ for $n_0(j)\le n<n_0(j+1)$. 
Use also 
Lemmas \ref{Hn01lm}  and \ref{Rn1lm} 
to see that 
$$\sum_{n=1}^\infty 
\left\{ P(H_{n,0})+ P(H_{n}(\delta_n))+
P\left( R_{n,1} \le - \e\right) \right\}
<\infty.
$$
It remains to show 
$
\liminf_{n\to\infty} R_{n,2}\ge 0
$ a.s. Recall the definition of the function 
$\phi_\beta$ in (\ref{fibetadef}). 
We get the desired conclusion by showing that 
$\phi_{\delta_n}(\zeta_n(\cdot,0))\to 0$ a.s. 
By (\ref{zetaB}), 
\beas
\phi_{\delta_n}(\zeta_n(\cdot,0))&=&
\sup
\left\{\,\left| B\left( n^{-1}T_{n,[nr]}\right) -  
 B\left( n^{-1}T_{n,[nq]}\right) \,\right|: \right.\\
&&\qquad \left. \mbox{ $|r-q|\le\delta_n$ and $r,q\in[c,d]$}\right\} +O(n^{-1/2}).
\eeas

Let $\e\in(0,1)$. 
Choose $\beta>0$ so that 
$\left|\int_q^r \rho^2_0(s)ds\right|<\e$ for $q,r\in[c,d]$ such that 
$|q-r|\le\beta$. Let $U_m$ be the event such  that
$$\left| n^{-1}T_{n,[nq]}-\int_0^q \rho^2_0(s)ds\right|\le \e 
\quad\mbox{ for all $q\in[c,d]$ and  $n\ge m$.}
$$
 Let $[g,h]$ be a compact interval 
such that $\int_0^q\rho^2_0(s)ds\in [g+1,h-1]$ for all $q\in[c,d]$. 
Then on the event $U_m$ we  also have 
$n^{-1}T_{n,[nq]}\in [g,h]$ for all  $q\in[c,d]$, and 
\be
\left| n^{-1}T_{n,[nq]}-n^{-1}T_{n,[nr]}\right|\le 3\e 
\label{Tnqrineq}
\ee
for all $q,r\in[c,d]$ such that 
$|q-r|\le\beta$, again for all $n\ge m$. 
Since $\delta_n<\beta$ for large $n$, we have on the event $U_m$
\be
\limsup_{n\to\infty} \phi_{\delta_n}(\zeta_n(\cdot,0))
\le \sup_{v,s\in[g,h], |v-s|\le 3\e}
\left| B(v)-B(s)\right|. 
\label{fiB}
\ee
By Lemma \ref{Tlm}
$\lim_{m\to\infty} P(U_m)=1$, so (\ref{fiB}) holds almost 
surely. Letting   $\e\searrow 0$ turns  (\ref{fiB}) into
$\lim_{n\to\infty} \phi_{\delta_n}(\zeta_n(\cdot,0))=0$ by the
path-continuity of $B(\cdot)$.  

This completes the proof of the lower bound: we now have 
\be
\liminf_{n\to\infty} \inf_{(x,t)\in A}
\left\{ \zeta_n(x,t)-\inf_{y\in I(x,t)} 
\zeta_n(y,0) \right\}\ge 0
\qquad\mbox{almost surely.}
\ee

\hbox{}

{\it Upper Bound, Case 1.} This is handled by large 
deviation estimates as was done in {\it Lower Bound, Case 1}
above. 

\hbox{}

{\it Upper Bound, Case 2.} Lemma \ref{Rn3lm} already gives almost sure
convergence for the error $R_{n,3}$. 

\hbox{}

This completes the proof of part (i) of Theorem \ref{thmloqeq}. 

Following Lemma \ref{suplm}, one can show that 
$M\equiv\sup_n\sup_{(x,t)\in A_0}|\zeta_n(x,t)|$ is a.s.\ finite
for an arbitrary compact set $A_0$. This 
furnishes the a.s.\ bound needed to turn the proof of 
part (ii) of Theorem \ref{probthm} (Section \ref{pfofprobthmii}) into
a proof for  part (ii) of Theorem \ref{thmloqeq}.

\section{The  distribution-valued processes}
\label{distrsect}

First we discuss the general setting of Theorem \ref{stickthm}.
An element $F\in\cal D'$ is determined by its actions on
$C^\infty_c(K)$'s for compact sets $K\subseteq\RR$, so
$R(F,G)$ is clearly a metric on $\Hloc$. The separability
and completeness of this metric follow from the separability
and completeness of $H^{-1}(\RR)$. 

Let us first check that $\xi_n(t)\in\Hloc$, or in other words
that $\chi\xi_n(t)\in H^{-1}(\RR)$ for any $\chi\in C^\infty_c(\RR)$.
Let $\varphi\in H^1(\RR)$. Set
$$
\ubar_n(x,t)=n^{1/2}\left( u(x,t)-u([nx]/n,t)\right).
$$
 By a summation
by parts, 
\beas
&&\chi\xi_n(t,\varphi)=\xi_n(t,\chi\varphi)\\
&=&n^{-1/2}\sum_{i\in\ZZ}\chi(i/n)\varphi(i/n)
\{ z^n_i(nt)-nu(i/n,t)-z^n_{i-1}(nt)+nu((i-1)/n,t)\}\\
&=&
-\sum_{i\in\ZZ}\left[\chi((i+1)/n)\varphi((i+1)/n)-\chi(i/n)\varphi(i/n)\right]
n^{-1/2}\{z^n_i(nt)-nu(i/n,t)\}\\
&=& -\sum_{i\in\ZZ}\int_{i/n}^{(i+1)/n}(\chi\varphi)'(x) dx \cdot
n^{-1/2}\{z^n_i(nt)-nu(i/n,t)\}
\\
&=&-\int (\chi\varphi)'(x) \left( \zeta_n(x,t)+\ubar_n(x,t)\right) dx.  
\eeas
Use $(\chi\varphi)'=\varphi\chi'+\varphi'\chi$ and the
Schwarz inequality. Let
$[a,b]$ contain the support of $\chi$. By the local Lipschitz
property of $u$ (Lemma \ref{Ixtlm}(c)), 
$\ubar_n=O(n^{-1/2})$ on any compact set.  We get
\beas
|\chi\xi_n(t,\varphi)|
&\le&
\|\varphi\|_{L^2(\RR)}\left( \int |\chi'(x)|^2 
\left\{ \zeta_n(x,t)+\ubar_n(x,t)\right\}^2 dx\right)^{1/2}\\
&&\quad +
\|\varphi'\|_{L^2(\RR)}\left( \int |\chi(x)|^2
\left\{ \zeta_n(x,t)+\ubar_n(x,t)\right\}^2 dx\right)^{1/2}\\
&\le& C\cdot \|\varphi\|_{H^1(\RR)}\left\{\left( \int_a^b  \zeta_n(x,t)^2 dx 
\right)^{1/2}+\frac1{\sqrt{n}}\right\}.
\eeas
This verifies that $\chi\xi_n(t)\in H^{-1}(\RR)$, because for any fixed
$\omega$, the process $\zeta_n(x,t)$ is locally bounded in $(x,t)$.
In other words, $\xi_n(t)\in\Hloc$. 

Next we check that $\xi_n(t)$ is measurable as a random element of
$\Hloc$. For fixed $\chi$ and $\varphi$, the function 
$\omega\mapsto \xi_n(t,\chi\varphi;\omega)$ is measurable. 
For a fixed $F\in\Hloc$, the metric $R(F,\xi_n(t;\omega))$ is a measurable
function of $\omega$, because the supremum  in
$$
\|\chi_kF-\chi_k\xi_n(t)\|_{H^{-1}(\RR)}
=\sup_{\|\varphi\|_{H^1(\RR)}\le 1}| F(\chi_k\varphi)-\xi_n(t,\chi_k\varphi)|
$$
can be restricted to countably many $\varphi$'s due to the separability of 
$H^1(\RR)$. Now the Borel measurability of the $\Hloc$-valued 
function $\omega\mapsto \xi_n(t;\omega)$ follows from the separability 
of $\Hloc$. 

Right-continuity of the path $t\mapsto \xi_n(t;\omega)\in\Hloc$
follows from the right-continuity and local boundedness
of the process $\zeta_n$, via a calculation that resembles the
one performed above. 

To summarize,  $\xi_n(\cdot)$ is a random 
element of the Skorokhod space $D([0,\infty),\Hloc)$. We leave
it to the reader to show that $\xi(\cdot)$ is a random 
element of $C([0,\infty),\Hloc)$. 

We move to the main point, 
to prove the strong law $\sup_{0\le t\le \tau} R(\xi_n(t),\xi(t))\to 0$
a.s. Note that in the definition (\ref{metricdef}) of $R(F,G)$
the quantities $\|\chi_kF-\chi_kG\|_{H^{-1}(\RR)}$ are nondecreasing
in $k$. Consequently, for any $k$, 
$$R(F,G)\le \|\chi_kF-\chi_kG\|_{H^{-1}(\RR)} + 2^{-k}. $$
So it suffices to show that for any $k$,  almost surely
$$
\lim_{n\to\infty}\sup_{0\le t\le \tau}\sup_{\|\varphi\|_{H^1(\RR)}\le 1}
|\xi_n(t,\chi_k\varphi)-\xi(t,\chi_k\varphi)|=0.
$$
Following the earlier calculation and by the definition (\ref{xidef}), we get
\beas
&&|\xi_n(t,\chi_k\varphi)-\xi(t,\chi_k\varphi)|\\
&\le&
\left|-\int (\chi_k\varphi)'(x)\zeta_n(x,t)dx+
\int (\chi_k\varphi)'(x)\zeta(x,t)dx\right| +
\left|\int (\chi_k\varphi)'(x)\ubar_n(x,t)dx\right|\\
&\le&
C_k\cdot 
\|\varphi\|_{H^1(\RR)}\left\{\left( \int_a^b  |\zeta_n(x,t)-\zeta(x,t)|^2 dx
\right)^{1/2}+\frac1{\sqrt{n}}\right\}.
\eeas
The constant $C_k$ is determined by $\chi_k$. 
Consequently
\beas
&&\sup_{0\le t\le \tau}\sup_{\|\varphi\|_{H^1(\RR)}\le 1}
|\xi_n(t,\chi_k\varphi)-\xi(t,\chi_k\varphi)|\\
&\le& 
C_k\cdot \left\{\left( \sup_{0\le t\le \tau}\int_a^b  |\zeta_n(x,t)-\zeta(x,t)|^2 dx
\right)^{1/2}+\frac1{\sqrt{n}}\right\}
\eeas
which converges a.s.\ to 0 by Theorem \ref{thmloqeq}(ii). 
This completes the proof of Theorem \ref{stickthm}.

\hbox{}

As the final item of this section, we verify the correlation formula
(\ref{corr1}), which requires us to show that, for 
$\phi,\psi\in C^\infty_c(\RR)$, 
\beas
&&  E\left[\int_\RR\int_\RR \phi'(x)\psi'(z)B\left(\int_0^{y^+(x,t)}\rho^2_0(r)dr\right)
B\left(\int_0^{y^+(z,s)}\rho^2_0(r)dr\right) \,dx\,dz\right]\\
&=& \int_\RR \phi(w^+(r,t))\psi(w^+(r,s))\rho^2_0(r)dr.
\eeas
Recall that in the above integrals $y^\pm(x,t)$ are interchangeable
because they differ only on Lebesgue null sets, and the same for
$w^\pm(r,t)$. In the sequel we manipulate Lebesgue-Stieltjes integrals
of the form $\int_a^b f\,dG$. Then we always use the right-continuous version of
$G$, with the measure defined by $\mu(a,b]=G(b)-G(a)$, and in case jumps
make a difference,  the integral is taken over the set $(a,b]$.

   We shall make use of the following integration by
parts formula. Its proof follows from standard integration by parts
\cite[Theorem 3.36]{folland} and Lemma \ref{FGHlm}. Let 
 $f\in L^1(\RR)$ and $\varphi$ be  compactly supported
and differentiable. Then 
\bea
\int_a^b \varphi'(x)\int_{-\infty}^{y^+(x,t)}f(r)dr\,dx
&=&\varphi(b)\int_{-\infty}^{y^+(b,t)}f(r)dr -
\varphi(a)\int_{-\infty}^{y^+(a,t)}f(r)dr \nn\\
&&\quad -
\int_{y^+(a,t)}^{y^+(b,t)}\varphi(w^+(r,t))f(r)dr.
\label{ibpeqn1}
\eea

To begin, recall that $B(\cdot)$ is a two-sided Brownian motion with independent
halves, so if $rq<0$ then $E[B(r)B(q)]=0$. Assume $s\le t$ without loss of generality.
   Then 
\beas
&&  E\left[\int_\RR\int_\RR \phi'(x)\psi'(z)B\left(\int_0^{y^+(x,t)}\rho^2_0(r)dr\right)
B\left(\int_0^{y^+(z,s)}\rho^2_0(r)dr\right) \,dx\,dz\right]\\
&=&
\int_\RR\int_\RR \phi'(x)\psi'(z) \left\{\int_0^{y^+(x,t)\wedge y^+(z,s)}\rho^2_0(r)dr\right\}
{\bf 1}\{y^+(x,t)\wedge y^+(z,s)>0\}\,dx\,dz \\
&+&
\int_\RR\int_\RR \phi'(x)\psi'(z) \left\{\int_{y^+(x,t)\vee y^+(z,s)}^0\rho^2_0(r)dr\right\}
{\bf 1}\{y^+(x,t)\vee y^+(z,s)<0\}\,dx\,dz \\
&\equiv& A_++A_-,
\eeas
where the last line defines the abbreviations $A_\pm$. We do the calculation
for $A_+$ and leave the similar steps for $A_-$ to the reader. 
\beas
A_+&=& \int_{w^+(0,t)}^\infty dx\int_{w^+(0,s)}^\infty dz\, \phi'(x)\psi'(z)
 \left\{\int_0^{y^+(x,t)\wedge y^+(z,s)}\rho^2_0(r)dr\right\}\\
&=& \int_{w^+(0,s)}^\infty dz\,\psi'(z)
 \left\{\int_0^{ y^+(z,s)}\rho^2_0(r)dr\right\}
\int_{w(z;s,t)}^\infty dx\, \phi'(x)\\
&&\qquad +
\int_{w^+(0,s)}^\infty dz\, \psi'(z) \int_{w^+(0,t)}^{w(z;s,t)} dx\,
\phi'(x) \int_0^{ y^+(x,t)}\rho^2_0(r)dr\\
&=&
-\int_{w^+(0,s)}^\infty dz\, \psi'(z)\phi(w(z;s,t))\int_0^{ y^+(z,s)}\rho^2_0(r)dr\\
&&\qquad+
\int_{w^+(0,s)}^\infty dz\, \psi'(z)\phi(w(z;s,t))\int_0^{ y^+(w(z;s,t),t)}\rho^2_0(r)dr\\
&&\qquad -
\int_{w^+(0,s)}^\infty dz\, \psi'(z)\phi(w^+(0,t))\int_0^{ y^+(w^+(0,t),t)}\rho^2_0(r)dr \\
&&\qquad-
\int_{w^+(0,s)}^\infty dz\, \psi'(z) \int_{ y^+(w^+(0,t),t)}^{y^+(w(z;s,t),t)}
\phi(w^+(r,t))\rho^2_0(r)dr, 
\eeas
where the last equality came from applying (\ref{ibpeqn1}). Observe that
\beas
&&\mbox{$w(z;s,t)=w^+(r,t)$ for $y^+(z,s)<r<y^+(w(z;s,t),t)$, and}\\
&&\mbox{$w^+(0,t)=w^+(r,t)$ for $0<r<y^+(w^+(0,t),t)$. }
\eeas
Then the terms above add up to give
$$
A_+=-\int_{w^+(0,s)}^\infty dz\, \psi'(z) \int_0^{y^+(z,s)} \phi(w^+(r,t))\rho^2_0(r)dr
$$
which after an application of (\ref{ibpeqn1}) is
\beas
&=&\psi(w^+(0,s))\int_0^{y^+(w^+(0,s),s)}\phi(w^+(r,t))\rho^2_0(r)dr\\
&&\qquad 
+\int_{y^+(w^+(0,s),s)}^\infty\psi(w^+(r,s))\phi(w^+(r,t))\rho^2_0(r)dr.
\eeas
Observe that $w^+(0,s)=w^+(r,s)$ for $0<r<y^+(w^+(0,s),s)$, to turn this into
$$
A_+=\int_0^\infty\psi(w^+(r,s))\phi(w^+(r,t))\rho^2_0(r)dr.
$$
Similar arguments show that $A_-$ is the complementary integral $\int_{-\infty}^0$. Equation (\ref{corr1}) is proved.

\section{Appendix: Some technical issues}
\label{measurabilitysection}

The assumptions in this section are the same as those of
Theorem \ref{probthm}. 
We first check the measurability of certain functions and 
sets. When convenient we add the sample point $\omega$ as an argument to
random quantities.   At the end of this section
 we indicate how Proposition \ref{kellprop}
follows from the estimates of Sections \ref{gammasect} and 
\ref{varcoupsect}. 

\begin{prop} For each $n$, the function 
$\zeta_n(x,t,\omega)-\inf_{y\in I(x,t)}\zeta_n(y,0,\omega)$
is jointly measurable in $(x,t,\omega)$. 
\label{mbleprop1}
\end{prop}

{\it Proof.} The term $\zeta_n(x,t,\omega)=
n^{-1/2}\{z^n_{[nx]}(nt,\omega)-nu(x,t)\}$ needs no
special argument, as for each $k$ the variable 
$z^n_k(t,\omega)$ is right-continuous in $t$ and hence
progressively measurable in $(t,\omega)$.  

Fix $n$ and consider the term $\sigma(x,t,\omega)=
\inf_{y\in I(x,t)}\zeta_n(y,0,\omega)$. For integers $m$ and
$(x,t,y)\in \RR\times[0,\infty)\times\RR$ let
$$h_m(x,t,y)=
\left\{ 
\begin{array}{rl}
 0 &\mbox{if $ y\in \bigcup_{q\in I(x,t)}[q,q+1/m]$,}\\
\infty &\mbox{otherwise.}
\end{array}
\right.
$$
Check that $h_m$ is lower semicontinuous, and hence Borel
measurable in $(x,t,y)$.  Consequently
$$
\sigma^{(m)}(x,t,\omega)=\inf_{y\in\QQ}
\{ \zeta_n(y,0,\omega)+h_m(x,t,y)\}
$$
is measurable in $(x,t,\omega)$. [The infimum is over 
rational $y$.] 

It remains to check $\sigma^{(m)}(x,t,\omega)\to \sigma(x,t,\omega)$
as $m\to\infty$. We leave this to the reader. 
\qed

In several places in the paper we needed
 the measurability of the function 
$$Z(\omega)= \sup_{(x,t)\in A}
\left| \zeta_n(x,t,\omega)-\inf_{y\in I(x,t)} 
\zeta_n(y,0,\omega) \right|
$$
where $A$ is one of three types of compact sets: 
(i) $A=[a,b]\times[0,\tau]$, (ii) $A$ has no shocks, or 
(iii) $A$ is
finite. Finite is of course trivial. We give the proof
for type (i)  and leave the (simpler) type (ii) to the 
reader. 

\begin{prop} Let  $A=[a,b]\times[0,\tau]$ in the definition 
of $Z(\omega)$.  
Then there exists a countable set ${\cal S}\subseteq 
[a,b]\times(0,\tau]$ such that 
\be
Z= \sup_{(x,t)\in\cal S}
\left| \zeta_n(x,t)-\inf_{y\in I(x,t)} 
\zeta_n(y,0) \right|. 
\label{ZSdef}
\ee
Consequently, by Proposition \ref{mbleprop1},
 $Z$ is a measurable function on the
probability space of the process $z^n(\cdot)$. 
\label{mbleprop2}
\end{prop}

{\it Proof.} 
 Let $\cal S$ be a countable subset of
$[a,b]\times(0,\tau]$ that satisfies these requirements:

(i) $\cal S$ contains a subset ${\cal S}'\subseteq\cal S$
such that 
(a) 
each $(x,t)\in\cal S'$ has $y^\pm(x,t)=y(x,t)$,
in other words, points in $\cal S'$ are not
shocks; and 
(b)  $\cal S'$ is dense in $[a,b]\times(0,\tau]$,
and  dense in the boundary line segments 
$\{a\}\times(0,\tau]$, $\{b\}\times(0,\tau]$, and 
$[a,b]\times\{\tau\}$. 
 
(ii) 
$\cal S$ contains all the shocks on the boundary line
segments, and on the 
vertical segments $\{k/n\}\times(0,\tau]$
where $k$ ranges over integers such that $k/n\in [a,b]$. 
Also, $\cal S$ contains the points $(x,\tau)$ for 
$x=a$, $x=k/n$ for any integer $k$ such that $k/n\in[a,b]$,
and for $x=b$. 

(iii) For each integer $\ell$, let $V_\ell$ be the 
set of points $(b,t)$, $0<t\le \tau$, such that
$(b,t)$ is not a shock and $y(b,t)=\ell/n$. (We already
included shocks $(b,t)$ in $\cal S$ in step (ii).)
Include in $\cal S$ a dense countable  subset  of $V_\ell$
so that each point of $V_\ell$ can be approached from
above by a point of $\cal S$. 

(iv) $\cal S$ contains all the $\cal U$-points in 
$[a,b]\times(0,\tau]$. 

Requirements (i)--(ii) can be satisfied because there 
are no more than countably many shocks  on any horizontal
or vertical line segment. 
Requirement (iv) can be satisfied by Theorem \ref{Uthm}.

Suppose now that $(x,t)$ is an arbitrary point 
of $[a,b]\times(0,\tau]$ outside $\cal S$.
 Since  
$\cal U$-points are in $\cal S$, it follows 
that $I(x,t)=\{y^\pm(x,t)\}$. We first show that 
we can find a sequence 
$(x_j,t_j)$ in $\cal S$ such that $y^\pm(x_j,t_j)=y_j$
and 
\be
\lim_{j\to\infty} \left| \zeta_n(x_j,t_j)-\zeta_n(y_j,0)\right|
= \left| \zeta_n(x,t)-\zeta_n(y^+(x,t),0)\right|.
\label{interm4}
\ee

Start with the case $x<b$. 
Find $(x_j,t_j)\in\cal S$ 
 so that
$I(x_j,t_j)=\{y_j\}$,  
 $x_j\searrow x$, $t_j\searrow t$,
and  so that 
$x_j> w(x;t,t_j)$ [possible 
because $w(x;t,t_j)\searrow x$ as $t_j\searrow t$]. 
If $t=\tau$ we can choose $t_j=\tau$ for all $j$. 
Notice that since $x_j$ approaches $x$ from the right, 
$[nx_j]=[nx]$ for large enough $j$, and since the 
jump processes are right-continuous in time, 
$z^n_{[nx_j]}(nt_j)=z^n_{[nx]}(nt)$ for large enough $j$. 
From $x_j> w(x;t,t_j)$ we have   $y_j\searrow y^+(x,t)$, so
$z^n_{[ny_j]}(0)=z^n_{[ny^+(x,t)]}(0)$ for large $j$. 
Since $u(x,t)$ is continuous, we get 
 (\ref{interm4}) for $x<b$. 

Now let $(x,t)=(b,t)\notin\cal S$. Then 
$(b,t)$ cannot be a shock so $y^\pm(b,t)=y(b,t)$. 
Pick $(b,t_j)\in\cal S$ so that $t_j\searrow t$. 
Exactly as above, $\zeta_n(b,t_j)\to \zeta_n(b,t)$
as $j\to\infty$. 
Now $y_j=y(b,t_j)$ satisfies $y_j\nearrow y(b,t)$.
We still have $\zeta_n(y_j,0)\to \zeta_n(y(b,t),0)$,
except possibly in 
 the case where $ny(b,t)$ is an integer, when it may happen
that $[ny_j]=[ny(b,t)]-1$ for all large enough $j$.
But by part (iii) of the definition of $\cal S$, then 
we can choose  $(b,t_j)\in\cal S$ so that 
$y_j=y(b,t)$. 

We have checked (\ref{interm4}) for all $(x,t)\in [a,b]\times(0,\tau]$.
Repeating (\ref{interm4}) with $y^-(x,t)$ in place of
$y^+(x,t)$ is trickier, so we exclude all cases where
 the $\zeta_n(y^-(x,t),0)$ alternative is irrelevant for the
value of $Z$. 
First, for  the $\zeta_n(y^-(x,t),0)$ alternative to matter
we must have 
\be
\zeta_n(y^-(x,t),0)=\min \{\zeta_n(y^-(x,t),0), 
\zeta_n(y^+(x,t),0)\} < \zeta_n(y^+(x,t),0).
\label{interm4a}
\ee
Secondly, the quantity $|\zeta_n(x,t)-\zeta_n(y^-(x,t),0)|$ will
not influence the value of $Z$ unless at least 
$$
 \left| \zeta_n(x,t)-\zeta_n(y^-(x,t),0)\right|
> \lim_{j\to\infty} \left| \zeta_n(x_j,t_j)-\zeta_n(y_j,0)\right|
$$
where $(x_j,t_j)$ is the sequence from $\cal S$ that 
appears in (\ref{interm4}). This implies
$$
 \left| \zeta_n(x,t)-\zeta_n(y^-(x,t),0)\right|
> \left| \zeta_n(x,t)-\zeta_n(y^+(x,t),0)\right|
$$
which together with (\ref{interm4a}) implies that 
$$
 \zeta_n(x,t)-\zeta_n(y^-(x,t),0) >0.
$$
Thus to show that $\left| \zeta_n(x,t)-\zeta_n(y^-(x,t),0)\right|$
is no larger than the supremum over $\cal S$
in (\ref{ZSdef}), it suffices to 
 find a sequence 
$(x_j,t_j)$ in $\cal S$ such that $y^\pm(x_j,t_j)=y_j$
and 
\be
\limsup_{j\to\infty} \left\{ \zeta_n(x_j,t_j)-\zeta_n(y_j,0)\right\}
\ge  \zeta_n(x,t)-\zeta_n(y^-(x,t),0).
\label{interm4b}
\ee

Consider first $x\in\{a,\{k/n\}\}$. Then we choose
$(x_j,t_j)=(x,t_j)$ so that $t_j\searrow t$. (If 
$t=\tau$ we would not be able to approximate $t$ from above; this is
why the point $(x,\tau)$ was included explicitly in $\cal S$.)
Now $y_j=y(x_j,t_j)$ satisfies $y_j\nearrow y^-(x,t)$.
Depending on whether $ny^-(x,t)$ is an integer or not
and whether $y_j<y^-(x,t)$ or equal, 
$z^n_{[ny_j]}(0)$ converges to $z^n_{[ny^-(x,t)]}(0)$  
or to  $z^n_{[ny^-(x,t)]-1}(0)$. In either case
$\lim_{j\to\infty}\zeta_n(y_j,0)\le \zeta_n(y^-(x,t),0)$.
Right $t$-continuity of the random dynamics and the
continuity of $u$ give 
$\zeta_n(x,t_j)\to \zeta_n(x,t)$. Thus (\ref{interm4b}) holds
in this case. 

It remains to consider $x\notin\{a,\{k/n\}\}$
 in (\ref{interm4b}). Now pick $t_j\searrow t$
and $x_j\nearrow x$. (Again if $t=\tau$ we can take $t_j=\tau$.)
This forces $y_j\nearrow y^-(x,t)$. 
Since $nx$ is not an integer, $[nx_j]=[nx]$ for large enough
$j$, and the right $t$-continuity of the random dynamics 
together with the continuity of $u$ gives $\zeta_n(x_j,t_j)\to \zeta_n(x,t)$.
The argument of the last paragraph again gives
$\lim_{j\to\infty}\zeta_n(y_j,0)\le \zeta_n(y^-(x,t),0)$
and (\ref{interm4b}) holds. 
 \qed

As the last measurability issue we show that the event
$H_n$ in Lemma \ref{Hnlm} is a measurable subset of the 
underlying probability space. 

\begin{prop} Fix $n$ and $\delta>0$. 
Let $A$ be a compact subset of
$\RR\times[0,\infty)$. Let 
\beas
H&=&\{\omega:\, z^n_{[nx]}(nt,\omega)=z^n_i(0,\omega)+
\Gamma^{n,i}_{[nx]-i}(nt,\omega)
        \ \mbox{for some $(x,t)\in A$ }\\
            &&\qquad \mbox{and some $i$ such that 
$\dist\left(i/n, I(x,t)\right)>\delta$}\}
\eeas
Then $H$ is a measurable event on the probability space of the 
process $z^n(\cdot)$. 
\label{mbleprop3}
\end{prop}

{\it Proof.} Fix a left-closed right-open bounded
 rectangle $[a,b)\times[0,\tau)$ that
contains $A$, and such that $an$, $bn$ and $\tau n$ are integers. 
Let $p$ be a positive integer. For integers $v,w$
such that  $2^p na+1\le v\le 2^pnb$ and 
$1\le w\le 2^pn\tau$, let 
$$
K^p_{v,w}=\left[ \frac{v-1}{2^pn},  \frac{v}{2^pn}\right)
\times 
\left[ \frac{w-1}{2^pn},  \frac{w}{2^pn}\right)
$$
be a tiling of  $[a,b)\times[0,\tau)$  with small rectangles. 
The size goes by multiples of $2^{-p}$ so that if 
$p'>p$ then each $K^{p'}_{v',w'}$ lies inside a unique
$K^p_{v,w}$. $\overline{K}^p_{v,w}$ is the closure. 
Let 
$$
I^p(v,w)=\bigcup_{(x,t)\in \overline{K}^p_{v,w}} I(x,t).
$$
For integers $L<0$ put  
$$
J^L_{v,w}=\{i\in\ZZ:  Ln\le i\le[2^{-p}v]\,,\, 
 \dist(i/n, I^p(v,w))>\delta\}.
$$
Define the measurable event
$$
V^L_{v,w}=\bigcup_{i\in J^L_{v,w}}
\{ \omega: z^n_{[2^{-p}v]}(2^{-p}w,\omega)=z^n_i(0,\omega)+
\Gamma^{n,i}_{[2^{-p}v]-i}(2^{-p}w,\omega)\}. 
$$
 Let ${\cal I}^p$ be the set of indices $(v,w)$ such that
$K^p_{v,w}$ intersects $A$. 
Let 
$$
U^L_p=\bigcup_{(v,w)\in{\cal I}^p}V^L_{v,w}.
$$

Our goal is now to show that 
\be
H=\bigcup_{L<0}\;\bigcup_{m\ge 1}\;\bigcap_{p\ge m} U^L_p
\quad\mbox{ a.s.}
\label{HisU}
\ee
The set on the right-hand side is evidently measurable,
and we conclude that so is $H$. 

Fix $L$ and suppose $\omega\in U^L_p$ for all large enough $p$. 
Then it is possible to choose a subsequence of $p$'s
and $(v_p,w_p)\in{\cal I}^p$ such that 
 $\omega\in V^L_{v_p,w_p}$ and the squares
$K^p_{v_p,w_p}$ are nested decreasing. 
Since $2^{-p}v\le nb$, there are in general only
finitely many choices for the index $i\in J^L_{v,w}$.  
 Thus by passing to
an even further subsequence we may assume that 
there is a fixed $i$ that satisfies $i\in J^L_{v_p,w_p}$
and 
\be
z^n_{[2^{-p}v_p]}(2^{-p}w_p)=z^n_i(0)+
\Gamma^{n,i}_{[2^{-p}v_p]-i}(2^{-p}w_p)
\label{interm5}
\ee
for all the $p$ in the subsequence. Since the squares
$K^p_{v_p,w_p}$ are nested, there is a fixed $k$ such that
$K^p_{v_p,w_p}\subseteq [k/n,(k+1)/n)\times [0,\tau)$. 
Now the treatment splits into two cases. 

{\it Case 1.} Suppose $2^{-p}n^{-1}v<(k+1)/n$ for some $p$ in the 
relevant subsequence. Then  $2^{-p}n^{-1}v$ is bounded away
from $(k+1)/n$ for all large
enough $p$ because the nesting of the $K^p_{v_p,w_p}$'s
forces $2^{-p}n^{-1}v$ to be nonincreasing.
Pass to the $p\to\infty$ limit along the relevant
subsequence. By the nesting and compactness, there
exists a point $(x,t)\in A$ such that $2^{-p}n^{-1}v
\searrow x$ and $2^{-p}n^{-1}w
\searrow t$.  Since the convergence comes from the right,
we can pass to the limit in (\ref{interm5}) to get
\be
z^n_{[nx]}(nt)=z^n_i(0)+
\Gamma^{n,i}_{[nx]-i}(nt). 
\label{interm6}
\ee
Note also that $(x,t)\in \overline{K}^p_{v_p,w_p}$
and  $i\in J^L_{v_p,w_p}$ imply $\dist(i/n, I(x,t))>\delta$.
Thus (\ref{interm6}) says that $\omega\in H$.

{\it Case 2.} Suppose $2^{-p}n^{-1}v=(k+1)/n$  for all $p$ in the 
relevant subsequence. Then after passing to the limit
$p\to\infty$ we have $x=(k+1)/n$. Again  (\ref{interm5})
gives  (\ref{interm6}) with the consequence  $\omega\in H$.

Conversely, we now show that $H$ lies a.s.\ in the event
on the right-hand side of (\ref{HisU}). 
 $H$ is a.s.\ the union of  the sets 
$$
H_L=\{\omega\in H: \mbox{$i_n(x,t,\omega)\ge Ln$ 
for $(x,t)\in[a,b]\times[0,\tau]$}
\}
$$
over $L<0$ [recall that $n$ is fixed now],
 so it suffices to consider $\omega\in H_L$ for
a  fixed $L$. Fix $(x,t)\in A$ and $i\in[Ln,[nx]]$ such that
$\dist(i/n,I(x,t))>\delta$ and 
(\ref{interm6}) holds.

 For each $p$ let $(v_p,w_p)$ be the
index such that $(x,t)\in K^p_{v_p,w_p}$. Pick $\beta>0$ so that 
$\dist(i/n,I(x,t)) >\delta+\beta$. For all $(x',t')$ close
enough to $(x,t)$, $I(x',t')$ is contained in the 
 $\beta$-neighborhood around $I(x,t)$. 
 Thus for large enough $p$, $I^p(v_p,w_p)$ lies in this 
$\beta$-neighborhood, and consequently $i\in J^L_{v_p,w_p}$. 
Increase $p$ so that $[nx]=[2^{-p}v]$ and so that neither
$z^n_{[nx]}(\cdot)$ nor $\Gamma^{n,i}_{[nx]-i}(\cdot)$ jumps
in the time interval $(nt, 2^{-p}w_p]$. Then for these large
enough $p$'s, (\ref{interm6}) implies (\ref{interm5}), which
says that $\omega\in V^L_{v_p,w_p}\subseteq U^L_p$.  
This completes the proof. \qed

Finally, we indicate briefly how to deduce Proposition
\ref{kellprop} from the estimates. Fix 
$\mu\in(2/3,1)$. The task is to show that 
$\lim_{n\to\infty}Z_{r,n}= 0$ in probability for $r=1,2$
where
$$
Z_{1,n}=\sup_{an\le k\le nb\,,\,0\le t\le n^{-\mu}}n^{-1/2}
\{ z^n_{k+\ell}(nt)-z^n_k(nt)\}
$$
and 
$$
Z_{2,n}=\sup_{an\le k\le nb\,,\, n^{-\mu} \le t\le \tau}
n^{-1/2}\{ z^n_{k+\ell}(nt)-z^n_k(nt)\}. 
$$
Let $i(k)$ be a minimizer for $z^n_k(nt)$ in the 
variational formula (\ref{varcoupn}). Bound $Z_{1,n}$ above by
$$
\sup_{an\le k\le nb\,,\,0\le t\le n^{-\mu}}n^{-1/2}
\{ z^n_{k+\ell}(0)-z^n_{i(k)}(0)\}\,,
$$
use Lemma \ref{Hn01lm} to bound $i(k)$ from below, and
appeal to assumption (\ref{ass1d}).

Bound $Z_{2,n}$ above by 
\beas
&&\sup_{an\le k\le nb\,,\, n^{-\mu}\le t\le \tau }n^{-1/2}
\left\{ \Gamma^{n,i(k)}_{k+\ell-i(k)}(nt)-\Gamma^{n,i(k)}_{k-i(k)}(nt)
\right\}\\
&\le& 
\sup_{an\le k\le nb\,,\, n^{-\mu}\le t\le \tau  }n^{-1/2} \left[
\cdot \Gamma^{n,i(k)}_{k+\ell-i(k)}(nt) \cdot
{\bf 1}\{ k-i(k)\le M_0(\log n)^{3/2}\} \right.\\
&&\qquad \left. +
\left( \Gamma^{n,i(k)}_{k+\ell-i(k)}(nt)-\Gamma^{n,i(k)}_{k-i(k)}(nt)\right)
\cdot
{\bf 1}\{ k-i(k)> M_0(\log n)^{3/2}\} \right],
\eeas
where $M_0$ is the constant appearing in Lemma \ref{largetlm2}.
Now apply Lemmas \ref{largetlm1} and \ref{largetlm2}.

Under the stronger assumptions of local equilibrium 
these estimations can be made summable in $n$ and a.s.\ convergence
follows by Borel-Cantelli.

\bibliographystyle{plain}

\enddocument